\documentclass[12 pt]{article}
\usepackage[top=2cm, bottom=2cm, left=2cm, right=2cm]{geometry}
\usepackage[utf8]{inputenc}
\usepackage[english]{babel}
\usepackage{polski}

\usepackage{amsmath}
\usepackage{parskip}
\usepackage{amssymb}
\usepackage{color}
\usepackage[all]{xy}  
\usepackage{tikz-cd}
\usepackage{float}
\usepackage{setspace}
\usepackage{hyperref}
\newtheorem{theorem}{Theorem}

\newtheorem{proposition}[theorem]{Proposition}

\newtheorem{definition}{Definition}

\newtheorem{remark}{Remark}

\def\d{\mathrm d}

\title{Reviewing the Geometric Hamilton-Jacobi Theory concerning Jacobi and Leibniz identities}
\author{O. Esen, M. de Le\'{o}n, M. Lainz, C. Sard\'{o}n, M. Zaj\k{a}c}
\date{January 2022}

\begin{document}

\maketitle

\begin{abstract}
    In this survey, we review the classical Hamilton–Jacobi theory from a geometric point of view in different geometric backgrounds. We propose a Hamilton-Jacobi equation for different geometric structures attending to one particular characterization: whether they fulfill the Jacobi and Leibniz identities simultaneously, or if at least they satisfy one of them.
    In this regard, we review the case of time-dependent ($t$-dependent in the sequel) and dissipative physical systems as systems that fulfill the Jacobi identity but not the Leibnitz identity. Furthermore, we review the contact-evolution Hamilton-Jacobi theory as a split off the regular contact geometry, and that actually satisfies the Leibniz rule instead of Jacobi.
    Furthermore, we include a novel result, which is the Hamilton-Jacobi equation for conformal Hamiltonian vector fields as a generalization of the well-known Hamilton-Jacobi on a symplectic manifold, that is retrieved in the case of a zero conformal factor.
    
The interest of a geometric Hamilton–Jacobi equation is the primordial observation that if a Hamiltonian vector field $X_H$ can be projected
into a configuration manifold by means of a $1$-form $dW$, then the integral curves of the projected vector field $X_{H}^{dW}$ can be transformed into 
integral curves of $X_H$ provided that $W$ is a solution of the Hamilton-Jacobi equation. Geometrically, the solution of the Hamilton-Jacobi equation plays the role of a Lagrangian submanifold of a certain bundle. Exploiting these features in different geometric scenarios we propose a geometric theory for multiple physical systems depending on the fundamental identities that their dynamic satisfies. Different examples are pictured to reflect the results provided, being all of them new, except for one that is reassessment of a previously considered example. 
\end{abstract}
\tableofcontents
\setlength{\parskip}{4mm}
\setlength{\parindent}{0pt}
\onehalfspacing
\section{Introduction}

The Hamilton-Jacobi theory is an alternative formulation of classical mechanics, equivalent to other classical theories as the Lagrangian or Hamiltonian formulation, as well as Newtonian mechanics. The Hamilton-Jacobi theory (henceforth HJ theory) is particularly useful in the identification of conserved quantities, even if the system cannot be resolved completely. In this formulation, one is able to represent a particle as a wave, and the wave equation obeyed by the mechanical systems is similar but not identical to a Schr\"{o}dinger's equation, this is why the HJ is considered the closest approach of classical mechanics to quantum mechanics \cite{Goldstein-book}. From the point of view of mathematics, the HJ equation makes an appearance in calculus of variations and it is a special case of the HJ-Bellman equation in dynamics programming \cite{Bellman}.

Furthermore, the HJ theory provides important physical examples connecting first order partial differential equations and systems of first order differential equations. For example, one obtains
  HJ-type equations whenever we consider a short-wave approximation for the solutions of wave-type equations, i.e., hyperbolic-type equations. In framework of wave mechanics, a complete solution of the HJ equation allows us to reconstruct an approximate solution of the Schr\"{o}dinger equation via the Van Vleck determinant of the complete solution itself \cite{CariHJ,carinena2006geometric,EMS}.
  
\textbf{Some Historical Remarks on Hamilton-Jacobi Theory.}  The HJ theory is rooted in the works of Sir William Rowan Hamilton (1805-1865) and Carl Gustav Jacob Jacobi (1804-1851). Hamilton
  carried out one of the earliest studies of geometrical optics in an arbitrary medium with varying index of refraction. He found a powerful expression of the topic in a characteristic function, which is the optical path length of a ray that satisfies the eikonal equation, regarded as a function of initial and final positions and times of the ray. This was first introduced in his Memoirs on Optics presented at the Royal Irish Academy. These functions satisfy partial differential equations, and directly determine infinite families of rays. Following an analogy between rays and trajectories of a mechanical system, Hamilton soon extended his concepts to mechanics, incorporating ideas of Lagrange and others concerning generalized coordinates. The resulting Hamiltonian mechanics, notable for its invariance under coordinate transformations, is a cornerstone of theoretical physics. Jacobi sharpened Hamilton's formulation, by clarifying some mathematical issues. The resultant Hamilton-Jacobi theory and later developments are presented in several famous texts: \cite{Arnold-book,BornWolf,Caratheodory,CourantHilbert,Gantmacher,Lanczos,landau2000mechanics}. For studies using modern PDE theory see \cite{Benton,Evans,Lions}.
  
  In a view broader than that of the original work, a solution of the Hamilton-Jacobi equation is the generator of a canonical transformation, a symplectic change of variables intended to simplify the equations of motion. In this framework (as applied to mechanics) there are solutions of a type different from that of Hamilton, which determine not only orbits but also invariant tori in phase space on which the orbits lie. These solutions, which are known to exist only under special circumstances, are the subject of the celebrated work of Kolmogorov, Arnol'd, and Moser; see \cite{Gallavotti}. Even approximate invariants, constructed by approximate solutions of the Hamilton-Jacobi equation, have implications for stability of motion over finite times \cite{Nehorosev,WarnockRuth}. Approximate invariants also find applications in the Einstein-Brillouin-Keller quantization of semi-classical quantum theory \cite{Keller}. Various forms and generalizations of the Hamilton-Jacobi equation occur widely in contemporary applied mathematics, for instance in optimal control theory \cite{FlemingRishel}. 
  
  Even if the Hamilton-Jacobi theory was proposed two centuries ago, there is still a lot of work remaining around this equation. We have mentioned that through the years, scientist have kept busy finding applications of the HJ theory in dynamics, stability, semiclassical quantum theory. Indeed, the geometric interpretation of the theory, in which this survey is based, is only a few years old \cite{CariHJ,Narciso,Cari16,Vaqtesis}. There are multiple recent publications on the geometry of the HJ equation, varying through the different author's wide range of interests. In the case of Cariñena, Román-Roy et al., they find their interest in the Hamilton–Jacobi theory was generated by the existence of bi-Hamiltonian descriptions for completely integrable dynamical systems and their desire to unveil and understand the quantum counterpart of bi- Hamiltonian systems, as well as nonholonomic systems \cite{CariHJ}. De León et al. have deeply inspected the geometry of the HJ theory on different geometric backgrounds, proposing an explicit expression for such equation: Nambu geometry, systems with dissipation terms, explicitly time-dependent problems, implicit Hamiltonians, multisymplectic backgrounds, field theories and conformal fields among many others \cite{LeonSar2,LeonSar1,LeonSar3, LeonZaj1,EsenLeonSar1,EsenLeonSar2,EsenLeonSarZaj1,EsJiLeSa19}.
  
  Other authors have found their interest in the HJ theory rooted in their necessity to solve discrete differential equations. Indeed, there are plenty of optimal control problems proposed in terms of a discrete Hamiltonian system. The first inklings of discrete mechanics appeared in the realm of Lagrangian mechanics \cite{MarsWest}, and the lack of a corresponding Hamiltonian theory lead to the development of discrete Hamiltonian mechanics. Since then, some works appeared on the discretization of Lagrangian and Hamiltonian systems on tangent and cotangent bundles, what lead
to variational principles for dynamical systems and principles of critical action on both the tangent and cotangent bundle \cite{GuiBloch2, GuiBloch}. 
This gave rise to analogies between discrete and continuous symplectic forms, Legendre transformations, momentum maps,  Noether's theorem and the discrete Hamilton-Jacobi theory \cite{OhsawaBlochLeok,OhsawaBlochLeok2}. 

In this review we are interested in presenting the HJ theory formulated on different geometric backgrounds attending to two fundamental geometric identities: the Jacobi and Leibniz identity. So, we dedicate different sections to structures that either fulfill the Jacobi identity, the Leibniz identity, or the two of them at the same time. It is our intention to generalize these results to geometric structures that do not satisfy these identities, as well as their discrete version.
Let us first introduce the HJ equation.

\subsection{Hamilton-Jacobi Equation}
Consider a Hamiltonian system detemined by the symplectic manifold $(T^*Q,\Omega_Q)$ and a Hamiltonian function $H$ on $T^*Q$. The HJ theory consists on finding a time-dependent {principal function}  $S(t,q^i)$, that fulfils
\begin{equation}\label{tdepHJ}
 \frac{\partial S}{\partial t}+H\left(q^i,\frac{\partial S}{\partial q^i}\right)=0
\end{equation}
where $\{q^i,p_i\}$ is the Darboux coordinates on $T^*Q$. Equation \eqref{tdepHJ} is referred to as the {Hamilton-Jacobi equation.}
 If we set the principal function to be separable in the
time variable $S=W(q^1,\dots,q^n)-Et$,
where $E$ is the total energy of the system, then \eqref{tdepHJ} will now read \cite{AbrahamMarsden,Goldstein-book}
\begin{equation}\label{HJeq1}
 H\left({q}^i,\frac{\partial W}{\partial {q}^i}\right)=E.
\end{equation}
The above equation is a very useful instrument to solve the Hamilton equations for $H$. When the Hamilton equations cannot be straightforwardly solved, one can give \eqref{HJeq1} a try and solve an equivalent problem.
Indeed, if we find a solution $W$ of \eqref{HJeq1}, then any solution of the Hamilton
equations is retrieved by taking ${p}_i=\partial W/\partial {q}^i$.

\textbf{Complete Solutions.}
On a broader sense, the most important result of the HJ theory is Jacobi's theorem, that states that a complete solution of the equation
\begin{equation}\label{tdepHJcomp}
 \frac{\partial S(t,q^i,\alpha_i)}{\partial t}+H\left(q^i,\frac{\partial S(t,q^i,\alpha_i)}{\partial q^i}\right)=0
\end{equation}
is given by the function $S$ that is generator of a canonical transformation with a symplectic change of variables intended to simplify the equations of motion (we can define this transformation at least locally using the implicit function theorem) from $(t,q^i,p_i)$ to $(t,\alpha^i,\beta_i)$. This implies that the complete solution will depend certain parameters $(\alpha_1,\dots,\alpha_n)$, provided that 
$$\text{det}\left(\frac{\partial^2 S}{\partial q^i\alpha_j}\right)\neq 0$$
and satisfies the equations
\begin{equation}
    \frac{\partial S(t,q^i,\alpha_i)}{\partial q^i}=p_i,\quad \frac{\partial S(t,q^i,\alpha_i)}{\partial \alpha^i}=-\beta_i.
\end{equation}
A simple computation shows that the system in these new coordinates reaches the dynamical equilibrium given by the equations:
\begin{equation}
    \frac{\partial \alpha_i}{\partial t}=0,\quad \frac{\partial \beta_i}{\partial t}=0. \end{equation}

The application of Jacobi’s theorem in the integration of Hamiltonian systems is usually based on the method of separation of variables. In this regard, we quote V.I. Arnol'd \cite{Arnold-book}:

\medskip

{\it ``The technique of generating functions for canonical transformations, developed by Hamilton and Jacobi, is the most powerful method available for integrating the differential equations of dynamics".}

\textbf{Geometric Hamilton-Jacobi Theorem.}
This survey is based on the geometric Hamilton-Jacobi approach, which is basically the translation of the classical HJ theory into a geometric language. Therefore, let us state the geometric reformulation of the HJ theory, which is rooted in the following assertion \cite{carinena2006geometric}: {\sl If a Hamiltonian
vector field $X_{H}$ can be projected into the configuration manifold by means of a 1-form $dW$, then the integral curves of the projected
vector field $X_{H}^{dW}$can be transformed into integral curves of $X_{H}$ provided that $W$ is a solution of \eqref{HJeq1}}.

This explanation can be represented by the following diagram:

\begin{equation}\label{geomdiagram-}
\xymatrix{ T^{*}Q
\ar[dd]^{\pi} \ar[rrr]^{X_H}&   & &TT^{*}Q\ar[dd]^{T\pi}\\
  &  & &\\
 Q\ar@/^2pc/[uu]^{dW}\ar[rrr]^{X_H^{dW}}&  & & TQ}
\end{equation}

\bigskip

This implies that $(dW)^{*}H=E$, with $dW$ being a section of the cotangent bundle. In other words, we are looking for a section $\alpha$ of $T^{*}Q$ 
such that $\alpha^{*}H=E$. As it is well-known, the image of a one-form is a Lagrangian submanifold of $(T^{*}Q, \Omega_Q)$ if and only if $d\alpha=0$.
That is, $\alpha$ is locally exact, say $\alpha=dW$ on an open subset around each point.
  
This is the geometric Hamilton--Jacobi theory enunciated on a symplectic manifold that is the cotangent bundle. Depending on the kind of manifold we are working on, the Lagrangian submanifold will be generalized to other types of submanifolds, as we shall point out.

\subsection{Outline and Content}
As we have discussed in the previous paragraphs, the geometric Hamilton-Jacobi theory has been generalized and applied in various geometries and to numerous physical systems. In this review paper we limit ourselves to recent advances on the Hamilton-Jacobi theory for classical systems generated by brackets of functions. In other words, we focus on the geometric Hamilton-Jacobi theory for Hamiltonian dynamics on Poisson, Jacobi and almost Poisson manifolds. These geometric investigations will find applications ranging from time-dependent Hamiltonian flows to dissipative dynamical systems, contact flows and nonholonomic systems. Before going into more details of the content, let us now briefly present how one can relax the algebraic conditions of a Poisson bracket to arrive at its generalizations.   

\textbf{Hierarchy of Brackets.}
A Poisson manifold admits a Poisson bracket that satisfies both the Leibniz and Jacobi identities \cite{BhaskaraViswanath,Dufour,Vaisman94,weinstein1983local,Weinstein98}. The Leibniz identity manifests the existence of Hamiltonian vector fields, whilst the Jacobi identity implies integrability.  We state two possible generalizations of Poisson manifolds. On one hand, one can relax the Jacobi identity: this gives rise to almost Poisson manifolds, where there is a skew-symmetric bracket satisfying only the Leibniz identity.

In the almost Poisson case \cite{CannasWeinstein}, one can further relax the requirement of being skew-symmetric, which reads Leibniz manifolds \cite{Ortega-Leibniz}. 
On the other hand, in a Poisson setting, instead of relaxing the Jacobi identity, one could relax the Leibniz identity. This gives Jacobi manifolds, see \cite{Marle-Jacobi,vaisman2002jacobi} and references therein. 
The following diagram represents these hierarchical relationships:  
\begin{center}
	\begin{tikzcd}
	\begin{matrix}
		\text{Leibniz Manifold} \end{matrix} \arrow[dd,"\text{$\{\bullet,\bullet\}$ skew-sym}",swap] \\\\
	\text{Almost Poisson Manifold} \arrow[ddr,"\mathfrak{J}=0",swap] && \text{Jacobi Manifold}\arrow[ddl,"Z=0"] \\\\
		&\text{Poisson Manifold}\arrow[dd,"\text{$\Lambda$ is nondegenerate}",swap] \\\\
		&\text{Symplectic Manifold}
		\end{tikzcd}
		\end{center}
See that, on the left upper hand, one starts with a bracket $\{\bullet,\bullet\}$ of smooth functions satisfying only the Leibniz identity. Then, by assuming the skew-symmetry of this bracket, one defines almost Poisson manifolds. Adding the requirement that the Jacobiator 
\begin{equation}\label{Jacobiator-}
\mathfrak{J}:\mathcal{F}(P)\times \mathcal{F}(P) \times \mathcal{F}(P)\longrightarrow \mathcal{F}(P), \quad (F,G,H)\mapsto ~ \circlearrowright \{F,\{H,G\}\},
\end{equation}
identically vanishes, one arrives at a Poisson manifold. Here, $\circlearrowright$ stands for the cyclic sum. On the other hand, by starting from the right upper hand, one has a bivector field $\Lambda$ and a vector field $Z$ satisfying \begin{equation}\label{ident-Jac-intro}
        [\Lambda,\Lambda] = 2 Z \wedge \Lambda, \qquad 
       [Z,\Lambda] = 0,
\end{equation}
    where $[\bullet,\bullet ]$ is the Schouten-Nijenhuis bracket \cite{Perelomov-SN,Marle-SN}. The case $Z=0$ determines a Poisson manifold. In this case, $\Lambda$ turns out to be a Poisson bivector. A non-degenerate Poisson bivector defines a symplectic two-form.  
    
\textbf{The content.}    Following this reasoning, we are going to start developing a HJ theory for these geometric structures and subcases of them. 
This survey consists of 3 main parts. Let us summarize them one by one.


		
\textbf{Section \ref{Sec-Poisson}: HJ on Poisson Manifolds.}
	In this section, we shall introduce the geometric Hamilton theory for conformal Hamiltonian dynamics and Hamiltonian dynamics on cosymplectic manifolds. The first case will be examined as a particular case of a HJ for dynamical systems subject to external forces, whereas cosymplectic will be the framework for $t$- dependent Hamiltonian motion. 

\textbf{Section \ref{Sec-Jacobi}: HJ on Jacobi Manifolds.} 	Contact manifolds and locally conformally symplectic (LCS) manifolds are two examples of Jacobi manifolds that admit a Hamilton-Jacobi formulation. This section contains a brief review of these theories. LCS manifolds look locally like symplectic manifolds, but not globally. In this sense, they are proper for gluing together all of the locally conformally symplectic patches to give rise to global Hamiltonian dynamics. Contact manifolds are odd-dimensional manifolds that are commonly used to describe certain physical phenomena, specially in thermodynamics.

\textbf{Section \ref{Sec-Almost}: HJ on Almost Poisson Manifolds.} We shall start this section by establishing a HJ theory for linear almost Poisson manifolds. This will be achieved on the dual of an almost Lie algebroid. Later, we shall determine a HJ for Hamiltonian dynamics with nonholonomic constraints. This section finishes with the study of almost Poisson Hamiltonian dynamics on contact manifolds. This approach is new in the literature, it is called evolution contact dynamics.

\textbf{References.} This review paper is based on the following publications of the authors:

\begin{itemize}

\item O. Esen, M. de León, C. Sardón, M. Zajac, Hamilton-Jacobi Formalism on Locally Conformally Symplectic Manifolds, {\it J. Math. Phys.} {\bf 62}, 033506 (2021).

\item O. Esen, V. Jiménez, M. de León, C. Sardón,
Reduction of a Hamilton — Jacobi Equation for Nonholonomic Systems,
{\it Regular and Chaotic Dynamics} {\bf 24}, 525- 559 (2019).

\item M. de León, C. Sardón,
Cosymplectic and contact structures for time-dependent and dissipative Hamiltonian systems, {\it J. Phys. A: Math. and Theor.} {\bf 50}, 255205 (2017).
	
\item M. de León, D.M. de Diego, M. Vaquero,
A Hamilton-Jacobi theory on Poisson manifolds,
{\it J. Geom. Mech.} {\bf 6}, 121 (2014).

\item M. de León, M. Lainz, Á. Muñiz-Brea,
The Hamilton–Jacobi Theory for Contact Hamiltonian Systems, in: {\sl  Special Issue New Trends in Hamilton-Jacobi Theory: Conservative and Dissipative Dynamics}
{\it Mathematics} 2021, {\bf 9}, (2021).

\item M. de León, M, Lainz Valcázar, Contact Hamiltonian systems,
{\it J. Math. Phys.} {\bf 60}, 102902 (2019).

\item O. Esen, M. Lainz Valcázar, M. de León, J.C. Marrero,
Contact Dynamics: Legendrian and Lagrangian Submanifolds,
{\it Mathematics} {\bf 9}, 2704 (2021).

\item O. Esen, M. Lainz Valcázar, M. de León, C. Sardón, Implicit Contact Dynamics and Hamilton-Jacobi Theory. arXiv preprint arXiv:2109.14921, (2021).

\end{itemize}

\section{HJ on Poisson Manifolds}\label{Sec-Poisson}
In this section we introduce the fundamentals for the main geometric structures in which we review our HJ theory. To reach this objective, we 
are interested in the most classical formalisms included in the Poisson category namely, symplectic, conformally symplectic, and cosymplectic Hamiltonian dynamics. We first summarize the HJ theory for the most common case, which is the case of symplectic dynamics.
Here, we shall review the time-dependent and time-independent geometry for Hamiltonian dynamics, construct the classical HJ theory and depict the complete solution of the classical symplectic Hamilton--Jacobi theory. 
Then we shall present HJ for conformal symplectic dynamics. Later, we shall depict HJ theory for cosymplectic Hamiltonian dynamics.

\subsection{Poisson Manifolds and Dynamics}\label{Sec-Poisson-Man}
A manifold $P$ equipped with a bracket defined on the space $C^\infty(P) $ of smooth functions \cite{BhaskaraViswanath,Dufour,Vaisman94,weinstein1983local,Weinstein98}   \begin{equation}\label{PoissonBracket}
  \{\bullet,\bullet\}:C^\infty(P) \times C^\infty(P)\longrightarrow C^\infty(P)
\end{equation}
is called Poisson if the bracket is skew-symmetric, satisfies the Leibnitz identity that is, for all $F$, $H$ and $G$ in $C^\infty(P) $
\begin{equation}\label{LeibId}
 \{F,H\cdot G\}= \{F,H\}\cdot G+H\cdot \{F,G\}
\end{equation}
and the Jacobi identity that is, for all $F$, $H$ and $G$ in $C^\infty(P) $
\begin{equation}\label{Jac-ident-Poiss}
\{F,\{H,G\}\}+\{H,\{G,F\}\}+\{G,\{F,H\}\}=0.
\end{equation}


\textbf{Hamiltonian Vector Fields.}
Given a function $H$, the Hamiltonian vector field $X_H$ is defined to be, for all $F$ in $\mathcal{F}(P)$,
\begin{equation}\label{Hamvf-}
X_H(F):=\{F,H\}.
\end{equation}
A function is called a Casimir function if it commutes with all the other functions (wrt. the Poisson bracket \eqref{PoissonBracket}). So, a function $C$ is a Casimir function if $\{C,F\}$ vanishes for all $F$ in $\mathcal{F}(P)$. Notice that the Hamiltonian vector field for a Casimir function is the zero vector field. The almost Hamiltonian dynamics generated by a Hamiltonian function $H$ is then computed to be 
\begin{equation} \label{HamEq}
\dot{z}=\{z,H\},
\end{equation} 
for $z$ in $P$.
The skew-symmetry of the bracket manifests that the Hamiltonian function $H$ is conserved all along the motion. 
Notice that, for a Poisson manifold, we have the identity that 
\begin{equation} \label{JL-iden}
[X_H,X_F]=-X_{\{H,F\}},
\end{equation}
where the bracket on the left hand side is the Lie bracket of vector fields. We denote the algebra of Hamiltonian vector fields on $P$ as $\mathfrak{X}_{ham}(P)$.

\textbf{Poisson Bivector Field.}
Due to the Leibniz identity, we identify an almost Poisson bracket $\{\bullet,\bullet\}$ with a bivector field $\Lambda$ on $P$ according to the definition 
\begin{equation} \label{bivec-PoissonBra}
\Lambda(dF,dH):=\{F,H\}
\end{equation}
for all $F$ and $H$ in $\mathcal{F}(P)$. Here, $dF$ and $dH$ are de Rham exterior derivatives of the smooth functions $F$ and $H$, respectively. So that we may represent a Poisson manifold by a pair $(P,\Lambda)$ as 
well. This geometry enables us to define a musical mappings $\Lambda^\sharp$  induced by the bivector field $\Lambda$ as
\begin{equation} \label{sharp}
\Lambda^\sharp: \Gamma^1(P)\longrightarrow \mathfrak{X}(P), \qquad \langle \beta, \Lambda^\sharp(\alpha)\rangle = \Lambda(\alpha, \beta),
\end{equation}
where the pairing on the left hand side is the one between the space $\Gamma^1(P)$ of one-form sections and the space $\mathfrak{X}(P)$ of vector fields. In terms of $\Lambda^\sharp$, the Hamiltonian vector field $X_H$ in \eqref{Hamvf} can be defined as 
\begin{equation}
\Lambda^\sharp(dH)= X_H.
\end{equation}

Recall the bivector field defined in \eqref{bivec-PoissonBra}. A direct calculation proves that  
\begin{equation} \label{SN-bra-JI}
\frac{1}{2}[\Lambda,\Lambda](dF,dH,dG)=~\circlearrowright \{F,\{H,G\}\}
=0,
\end{equation}
where the bracket on the left hand side is the Schouten-Nijenhuis bracket and $\circlearrowright $ denotes the cyclic sum. It is immediate to see from \eqref{SN-bra-JI} that $\{\bullet,\bullet\}$ is a Poisson bracket if and only if the bivector field $\Lambda$ commutes with itself under the Schouten-Nijenhuis bracket, that is, 
\begin{equation} \label{Poisson-cond}
[\Lambda,\Lambda]=0.
\end{equation} 
So, we can define a Poisson manifold by a pair $(P,\Lambda)$ with  $\Lambda$ satisfying the Jacobi identity \eqref{Poisson-cond}. A bivector field satisfying \eqref{Poisson-cond} is called a Poisson bivector. 
In terms of bivectors, we say that a differentiable map $\phi$ is called a Poisson map if $\Lambda$ and $\Lambda'$ are $\phi$-related.

\textbf{Local Picture of Poisson Manifolds.}
Consider a finite dimensional manifold $P$ and choose a local coordinate system $\{z^\alpha\}$. Here, $\alpha$ runs from $1$ to the dimension of the Poisson manifold. Then, a Poisson bracket is of the form 
\begin{equation}
\{F,H\}(z)=\Lambda^{\alpha \beta}(z)\frac{\partial F}{\partial z^\alpha}
\frac{\partial H}{\partial z^\beta},
\end{equation} 
where $\Lambda^{\alpha \beta}$ is a family of functions satisfying the skew-symmetry property, that is, $\Lambda^{\alpha \beta}=-\Lambda^{ \beta\alpha}$. The Poisson bivector field defined in \eqref{bivec-PoissonBra} turns out to be
\begin{equation}
\Lambda=\Lambda^{\alpha \beta}(z)\frac{\partial }{\partial z^\alpha}\wedge
\frac{\partial }{\partial z^\beta}.
\end{equation}
 In this local picture, the Jacobi identity \eqref{Poisson-cond} turns out to be 
\begin{equation}
\Lambda^{[\gamma \epsilon}_{,\beta} \Lambda^{\alpha] \beta}=0,
\end{equation}
where the subindex with comma denotes the partial derivative with respect to $z^\beta$, and the bracket denotes the cyclic permutation of the indices $\alpha$, $\gamma$ and $\epsilon$. If this holds, then $\Lambda$ is called a Poisson bivector. The Hamiltonian vector field $X_H$ and the Hamilton equations are computed to be
\begin{equation}
X_H= \Lambda^{\alpha\beta}\frac{\partial H}{\partial z^\alpha}\frac{\partial }{\partial z^\beta},
\qquad 
\dot{z}^\alpha=\Lambda^{\alpha\beta}\frac{\partial H}{\partial z^\beta},
\end{equation}
respectively.

\textbf{Characteristic Distribution.}
The characteristic distribution $\mathfrak{C}$ for a  Poisson manifold $(P,\Lambda)$ is a subbundle of the tangent bundle of $TP$ spanned by the values of the Hamiltonian vector fields  pointwisely as
\begin{equation}
\mathfrak{C}_p=\Lambda^\sharp(T^*_pP).
\end{equation}
The identity \eqref{JL-iden} reads 
that the characteristic distribution, spanned by Hamiltonian vector fields, of a Poisson bracket is integrable since the Jacobi-Lie bracket of two Hamiltonian vector field is another Hamiltonian vector field. A Poisson manifold is called transitive if its characteristic distribution is precisely the tangent bundle of it. In this case, $\Lambda^\sharp$ is invertible and the nondegenerate Poisson bivector $\Lambda$ determines a symplectic two-form $\Omega$ through 
\begin{equation}\label{sf-pt-}
\Omega^\flat:=-(\Lambda^\sharp)^{-1}
\end{equation}
where $\Omega^\flat$ is the musical map associated with $\Omega$ (c.f. \eqref{bemol}). 
In general, as a manifestation of the integrability of the characteristic distribution, a Poisson manifold defines a foliation of symplectic leaves. 

\textbf{Darboux Coordinates.}
There is a distinguished coordinate system at every point of a Poisson manifold $P$. It is called the Darboux-Weinstein coordinate system and provides the local picture of Poisson geometry. 
Assume that $P$ is a $(2m+k)$-dimensional Poisson manifold. Then at any point, there are coordinates $(q^i,p_i,y^a)$ where $i$ runs from $1$ to $m$ whereas $a$ runs from $1$ to $k$ such that the Poisson bivector field is written as
\begin{equation} \label{loc-L-bf}
\Lambda=\frac{\partial }{\partial q^i}\wedge
\frac{\partial }{\partial p_i} + \frac{1}{2}\lambda
^{ab}(z)\frac{\partial }{\partial y^a}\wedge
\frac{\partial }{\partial y^b}
\end{equation}
with $\lambda^{ab}(0)$ equal zero. For $k=0$ the Poisson manifold turns out to be even dimensional and admits a nondegenerate bivector field. 
 This is the case where a Poisson manifold is a symplectic manifold. 
 Another extreme case is $m=0$, which reads a degenerate Poisson bivector at the origin. A large class of Poisson manifolds are of this form, namely Lie-Poisson structures on the duals of Lie algebras \cite{Marsden1999}.

\textbf{Non-degeneracy.} If we can speak of non-degeneracy as admitting an isomorphism from the space of vector fields to the space of one-form sections, then for even dimensions one arrives at symplectic manifolds while for odd dimensions one arrives at cosymplectic manifolds. These particular instances of Poisson manifolds are determined through differential forms instead of bivector fields. 
We present the following table summarizing the geometric constructions addressed in this section. The table is summarizing the tensorial objects determining the geometries as well as the brackets. It is also exhibiting Hamiltonian vector fields and characteristic distributions. Later, we shall explicitly write on these realizations and associated HJ formalisms one by one.

\begin{table}[H]\label{Table-1}{\footnotesize
  \noindent
\caption{{\small {\bf Poisson Manifolds}. $M$ is a manifold. $\Lambda$ is a bivector field  and $[\Lambda,\Lambda]$ is the Schouten bracket. $\Lambda^\sharp$ is the musical mapping induced from $\Lambda$. $\Omega$ is a two-form, $\Omega^n$ is the $n$-th wedge power of it. $\Omega^\sharp$ is the musical mapping induced from $\Omega$. 
The flat mapping  $\flat:\mathfrak{X}(M)\rightarrow \Lambda^{1}(M)$ is a mapping between $C^{\infty}$ modules}.} 
\label{table1}
\medskip
\noindent\hfill
\resizebox{\textwidth}{!}{\begin{minipage}{\textwidth}
\begin{tabular}{ l l l l }
 
\hline
 &&\\[-1.5ex]
 Structure&  Characterization &Bracket and h.v.f.& Induced structure \\[+1.0ex]
\hline
 &  & \\[0.5ex]
{\bf Poisson} & $[\Lambda,\Lambda]=0$  &$X_H:= \Lambda^\sharp(dH)$ &    $\mathfrak{C}_p=\Lambda^\sharp(T_p^{*}M)\subset T_pM$  \\[+1.0ex] $(M,\Lambda)$ &  &$\{F,G\}=\Lambda(dF,dG)$ & Symplectic Leaves \boxed{\ker{\Lambda}=0} 
\\[0.5ex]
\hline\\[0.5ex]
{\bf Symplectic} & $\Omega^n\neq 0$ & $X_H:=\Omega^\sharp(dH)$ &   $\mathfrak{C}_p=\Omega^\sharp(T_p^{*}M)=T_pM$  \\[+1.0ex] $(M,\Omega)$ & $d\Omega=0$  &$\{F,G\}=\Omega(X_F,X_G)$ &   \\[+1.0ex] even $2n-$dim&   & & 
\\ [0.5ex] \hline
\\[0.5ex]
 {\bf Cosymplectic} &  $\Omega^n\wedge \eta\neq 0$& $\flat\circ \operatorname{grad} H=dH$ & $\flat(X)=\iota_X\Omega + \eta(X)\eta$,  \\[+1.0ex] $(M,\Omega,\eta)$   & $d\eta=0$  &  $\{F,H\}:=\Omega(\operatorname{grad} F,\operatorname{grad} H)$ & $\mathcal{R}=\sharp \eta$ where $\sharp=\flat^{-1}$ \\[+1.0ex] odd $(2n+1)-$dim & $d\Omega=0$& $\iota_{X_H}\Omega =dH- \mathcal{R}(H)\eta,\quad \iota_{X_H}\eta=0$  &$\mathfrak{C}_p=\sharp(T_p^{*}M) $    \\ [0.5ex]
\hline\\[0.5ex]
\end{tabular}
  \end{minipage}}
\hfill}
\end{table}

\subsection{Symplectic Manifolds and Dynamics}\label{Sec-Sym-Man}
In this section and following subsections we recall some basics on Hamiltonian dynamics on symplectic manifolds referring to an incomplete list \cite{AbrahamMarsden,Arnold-book,leon89,holm2009geometric,Liber87,Marsden1999}.

An  almost symplectic manifold is a pair $(M,\Omega)$, where $\Omega$ a non-degenerate two-form, i.e., an almost symplectic two-form. If additionally the two-form $\Omega$ is closed, then, $M$ is a symplectic manifold with a symplectic two-form $\Omega$.  If $\Omega$ is the exterior derivative  of a one-form $\Theta$, then it is said to be an exact symplectic two-form.

\textbf{Musical Mappings.}  For a nondegenerate two-form $\Omega _{%
M}$, the mapping 
\begin{equation}
\Omega^{\flat }:\mathfrak{X}\left( M\right)
\longrightarrow \Omega^{1} ( M),\qquad X\mapsto \iota_{X}\Omega
_{M},\label{bemol}
\end{equation}%
where $\iota_X$ is the interior derivative, is an isomorphism. 
The fiberwise inverse of $\Omega^{\flat }$ is 
\begin{equation}\label{diyez}
\Omega^{\sharp }:\Omega ^{1} ( M ) \longrightarrow 
\mathfrak{X}\left( M\right).
\end{equation}
In literature, the mappings  $\Omega^{\flat }$
and $\Omega^{\sharp }$ are called musical isomorphisms induced by the symplectic two-form $\Omega$.
Notice that we can relate the symplectic two-form $\Omega$ and the musical mapping $\Omega^{\sharp }$ through
\begin{equation}
\alpha(X)=\Omega(\Omega^{\sharp }(\alpha),X).
\end{equation}

\textbf{Symplectic Orthogonal.} Now, let $N$ be a submanifold of $M$. We define the symplectic orthogonal complement
of the tangent bundle $TN$ as the collection of tangent spaces
\begin{equation}
 T_zN^{\bot}=\{u\in T_zM |\ \Omega(u,v)(z)=0, \forall v\in T_zN\}.
\end{equation}
Here, $\Omega(u,v)(z)$ is the pointwise evaluation of the two-form with two vectors.
Generally, although it is not possible to express the tangent bundle $TM$ as the direct sum of $TN$  and $TN^{\bot}$, the dimension of $TM$ is the sum of the dimensions of $TN$ and its symplectic orthogonal complement $TN^{\bot}$ that is $\dim TM= \dim TN + \dim TN^{\bot}$.

Consider now a distribution $F\subset TM$ on $M$. The symplectic orthogonal $F^\bot$ is defined as 
\begin{equation*}
F_z^\bot=\{u\in T_zM: \Omega(u,v)=0,\forall v\in F\}. 
\end{equation*}
\noindent
Note that $F^\bot$ is also a distribution over $M$. We denote the annihilator of $F$ by $F^0$, consisting of one-forms on $M$ that vanish when restricted to $F$. In this case, $F^0$ is called a codistribution. Notice that, the musical mappings \eqref{bemol} and \eqref{diyez} satisfy the following identities
 \begin{equation} \label{flat-sharp}
\Omega^\flat(F)=(F^\bot)^0, \qquad F^{\bot}=\Omega^\sharp(F^{0}).
\end{equation}

\textbf{Submanifolds.} 
We have the following characterization of the submanifolds of a symplectic manifold $(M,\Omega)$:
 \begin{itemize}
  \item $N$ is called an isotropic submanifold if $TN\subset TN^{\bot}$. In this case, the dimension of $N$ is less or equal to the half of the dimension of $M$.
  \item $N$ is called a coisotropic submanifold if $TN^{\bot}\subset TN$. In this case, the dimension of $N$ is greater or equal to the half of the dimension of $M$.
  \item $N$ is called a Lagrangian submanifold   if $N$ is a maximal isotropic subspace of $(TM,\Omega)$. That is, if $TN=TN^{\bot}$. In this case, the dimension of $N$ is equal to the half of the dimension of $M$. 
  \item $TN$ is symplectic if $TN\cap TN^{\bot}=0$. In this case, $(N,\left.\Omega \right \vert_N)$ is a symplectic manifold. Here, $\left.\Omega \right \vert_N$ denotes the restriction of $\Omega$ to $N$.
 \end{itemize}

\textbf{Symplectic Diffeomorphisms.} Let $\left( M_{1},\Omega _{1}\right) $ and $\left( M_{2},\Omega _{2}\right) $ be two symplectic manifolds and $\psi:\mathcal{%
M}_{1}\rightarrow M_{2}$ be a diffeomorphism. This diffeo $\psi$ is called
symplectic (or canonical) diffeomorphism (aka. symplectomorphism)  if the pull back of the symplectic two-form $\Omega _{2}$ is precisely $\Omega _{1}$, that is,
\begin{equation}\label{symp-map}
\psi^{\ast }\Omega _{2}=\Omega _{1}.
\end{equation}
Under a symplectomorphism, the image of a Lagrangian (isotropic, coisotropic, symplectic) submanifold is a Lagrangian (resp. isotropic, coisotropic, symplectic) submanifold.
The composition of two symplectic mappings is symplectic and the set $\text{Diff}_{\text{can}}(M)$ of all symplectic diffeomorphisms on a symplectic manifold $M$ is an infinite dimensional Lie group.

\textbf{Symplectic Diffeomorphisms and Lagrangian Submanifolds.}
Consider two symplectic manifolds $\left( M_{1},\Omega _{1}\right) $ and $\left( M_{2},\Omega _{2}\right) $, and let $pr_{1}$ and $pr_{2}$ denote the
canonical projections
\begin{equation}
pr_{1}:M_{2}\times M_{1}\rightarrow M_{1},\qquad pr_{2}:M_{2}\times
M_{1}\rightarrow M_{2}.
\end{equation}%
Let us define a two form on $M_{2}\times M_{1}$ as 
\begin{equation}\label{omega-}
\Omega _{2}\ominus \Omega _{1}:=pr_{2}^{\ast }\Omega _{2}-pr_{1}^{\ast
}\Omega _{1}.
\end{equation}
It is possible to show that $\Omega _{2}\ominus \Omega _{1}$ is a symplectic two-form that is closed and non-degenerate. This makes the product space $M_{2}\times M_{1}$ a symplectic manifold. We define a diffeomorphism $\psi: M_{1}\rightarrow M_{2}$, such that the graph of $\psi$ is the subset of $M_{2}\times M_{1}$ given by
\begin{equation}\label{graph}
graph ( \psi ) =\left\{  ( \psi  ( z )
,z ) \in M_{2}\times M_{1} :\forall z\in M_{1}\right\} .
\end{equation}
A diffeomorphism $\psi: M_{1}\rightarrow M_{2}$ between two symplectic manifolds is a symplectic diffeomorphism if and only if graph$(
\psi)$ defined in \eqref{graph} is a Lagrangian submanifold of the product symplectic manifold $M_{2}\times M_{1}$ equipped with the symplectic two-form $\Omega _{2}\ominus \Omega _{1}$ in \eqref{omega-}.

\textbf{Geometric Hamiltonian Dynamics.}
A (globally) Hamiltonian vector field on a symplectic manifold $( M,\Omega) $ is a unique vector field $X_{H}$ satisfying 
\begin{equation}
\iota_{X_{H}}\Omega =dH,  \label{Hamvf}
\end{equation}
for a real valued function $H$ on $M$ (the function $H$ is the Hamiltonian). Here, $\iota_ {X_H}$ denotes the interior derivative with respect to $X_H$. A Hamiltonian system is determined by the triple $\left( M,\Omega,H\right) $.
In terms of the musical mappings defined in \eqref{bemol} and \eqref{diyez}, the geometric Hamilton equation can be written as \begin{equation}
\Omega^{\flat } (X_{H})=dH, \qquad \Omega^{\sharp }(dH)=X_{H}.
\end{equation}
Since $\Omega^{\sharp}$ is an isomorphism, one can
always find a Hamiltonian vector field $X_{H}$ for a given exact one-form $dH$. It is evident that two functions that differ in a constant give rise to the same Hamiltonian vector field.  We denote the space of Hamiltonian vector fields by $\mathfrak{X}_{ham}(M)$.

\textbf{Canonical Poisson Brackets.}  
Without loss of generalization, we start with the canonical symplectic manifold $(M,\Omega)$, and define the following bracket 
\begin{equation}\label{poissonsymplectic}
\{F,H\}:= \Omega\big(X_{F},X_{H}\big) ,  
\end{equation} 
where $X_{F}$ and $X_{H}$ are Hamiltonian vector fields. This is called a canonical Poisson bracket on the space of smooth functions $\mathcal{F}(M)$. Notice that we can define the Poisson bracket via  
\begin{equation}
X_{H}(F) =\{F,H\}
\end{equation}
as well. 
The pair $(\mathcal{F}(M),\{\bullet,\bullet\})$ is a Lie algebra, the bracket is skew-symmetric and satisfies the Jacobi identity.  
The skew-symmetry of the symplectic two-form $\Omega$ manifests the skew-symmetry of the canonical Poisson bracket, i.e., 
\begin{equation}
\left\{F,H\right\} =-\left\{H,F\right\}.
\end{equation}
A direct calculation shows that
\begin{equation}
\frac{1}{2}d\Omega( X_{F},X_{G},X_{H}) =\{F,\{
H,G\} \} + \{ H,\{G,F\}\} +\{G,\{ F, H \}\}.
\end{equation}
So that the closure of $\Omega$ manifests the fulfillment Jacobi identity
\begin{equation}
\{F,\{H,G\} \} + \{ H,\{G,F\}\} +\{G,\{ F, H \}\} =0
\end{equation}
for the canonical Poisson bracket. Further, it is possible to show that the canonical Poisson bracket satisfies the Leibniz identity 
\begin{equation}
\left\{ FG,H\right\} =\left\{F,H\right\}G+F\left\{ G,H\right\}. 
\end{equation}

\textbf{Space of Hamiltonian Vector Fields (Revisited).}
The following identity manifests that the space of Hamiltonian vector fields $\mathfrak{X}_{ham}(M)$ is closed under Lie bracket
\begin{equation}
\left[ X_{H},X_{F}\right]=-X_{\left\{H,F\right\} },
\label{Poissonham}
\end{equation}
where the bracket on the right hand side is the Lie bracket of vector fields. The identity \eqref{Poissonham} 
restates the fact that $\mathfrak{X}_{ham}(M)$ is a Lie algebra, but it also determines the Hamiltonian function $\{F,H\}$ generating the Lie algebra product $[ X_{H},X_{F}]$ of two Hamiltonian vector fields. Further, \eqref{Poissonham} gives that the mapping  
\begin{equation}
\mathcal{F}(M)\longrightarrow \mathfrak{X}_{ham}(M), \qquad H\mapsto X_H
\end{equation} 
is a Lie algebra anti-homomorphism. This mapping is far from being an isomorphism since its kernel consists of constant functions. 

\textbf{Cotangent Bundle.} 
The cotangent manifold is a symplectic manifold, since there exists a symplectic structure on the cotangent bundle $T^{\ast }%
Q$ that follows from the double vector bundle structure of the tangent of the cotangent bundle $TT^{\ast }Q$. 
For $TT^{\ast }Q$, we have the tangent
fibration and the tangent mapping of the cotangent fibration $\pi _{Q}$ written as
\begin{equation*}
\tau _{T^{\ast }Q}:TT^{\ast }Q\rightarrow T^{\ast }Q,\qquad 
T\pi _{Q}:TT^{\ast }Q\rightarrow TQ,
\end{equation*}
respectively. 
We collect these two fibrations in the following commutative diagram
\begin{equation}
\xymatrix{
&TT^{\ast }Q\ar[dl]_{T\pi _{Q}}\ar[dr]^{\tau _{T^{\ast }Q}}
\\TQ \ar[dr]_{\tau_{Q}}&&T^*Q\ar[dl]^{\pi_{Q}}
\\&Q } \label{T}
\end{equation} 
Accordingly, we define the canonical (Liouville) one-form $\Theta_{Q}$ on $T^{\ast }Q$ as 
\begin{equation}
\Theta_{Q}( \xi) =\left\langle \tau
_{T^{\ast }Q}\left( \xi \right) ,T\pi _{Q}\left( \xi
\right) \right\rangle ,  \label{canonicaloneform}
\end{equation}%
for a vector field $\xi$ on $T^{\ast }Q$ where the pairing on the right hand side is the pairing between $T^*Q$ and $TQ$. Minus of the exterior derivative of the canonical one-form 
\begin{equation}\label{omega-can}
\Omega_{Q}=-d\Theta_{Q}
\end{equation}
is a non-degenerate and closed two-form, so that the pair $(T^*Q,\Omega_{Q})$ becomes a symplectic manifold. In literature, $\Omega_{Q}$ is called a canonical symplectic two-form. There is a distinguished set of coordinates on a cotangent bundle which enables us to write the symplectic two-form as a constant two-form. 

\textbf{The Darboux Coordinates.} Locally, there exists a local coordinate system 
\begin{equation}
(q^i,p_i):T^*Q\longrightarrow \mathbb{R}^{2n}
\end{equation}
on the cotangent bundle, called Darboux coordinates, on cotangent bundle $T^*Q$ such that the canonical one form turns out to be $\Theta_{Q}=p_i dq^i$, whereas  the canonical symplectic two-form in \eqref{omega-can} becomes
\begin{equation}
\Omega_{Q}=dq^i \wedge dp_i.  \label{ss}
\end{equation}
The induced coordinates on the tangent bundle $TT^*Q$ are then given by 
\begin{equation}\label{coord-TT*Q}
(q^i,p_i,\dot{q}^i,\dot{p}_i):TT^*Q\longrightarrow \mathbb{R}^{4n}
\end{equation}
In terms of the Darboux coordinates and for the case of the cotangent bundle, the musical mapping \eqref{bemol}   turns out to be
\begin{equation}\label{beta1}
 \Omega_{Q}^\flat:TT^{*}Q \longrightarrow T^{*}T^{*}Q, \qquad (q^i,p_i,\dot{q}^i,\dot{p}_i)\mapsto(q^i,p_i,-\dot{p}_i,\dot{q}^i). 
 \end{equation} 
 In this local picture, for a Hamiltonian function $H=H(q,p)$ the Hamiltonian vector field is computed to be
 \begin{equation}
 X_H=\frac{\partial H}{\partial p_i}\frac{\partial  }{\partial q^i}-\frac{\partial H}{\partial q^i}\frac{\partial  }{\partial p_i}.
  \end{equation}
Accordingly, the Hamilton equation is 
 \begin{equation}
 \dot{q}^i=\frac{\partial H}{\partial p_i},\qquad \dot{p}_i=-\frac{\partial H}{\partial q^i}.
  \end{equation}
In Darboux coordinates the canonical Poisson bracket is obtained as 
 \begin{equation}
 \{F,H\}=\frac{\partial F}{\partial q^i}\frac{\partial H}{\partial p_i}-\frac{\partial F}{\partial p_i}\frac{\partial H}{\partial q^i}.
  \end{equation}

The Darboux theorem assures that locally all symplectic manifolds are equivalent to a cotangent bundle. 
All symplectic manifolds admit Darboux coordinates, so that the symplectic two-form can be written as a constant two-form as in  (\ref{ss}). 
This generic property of the cotangent bundle and ubiquity of the Darboux chart brings several advantages while performing calculations and proving assertions.

\subsection{HJ for Symplectic Dynamics}\label{Sec-HJ-Symp}
In the time-independent case of the Hamilton-Jacobi theory, we will consider the cotangent bundle case, although the theory can be applied locally to any symplectic manifold using Darboux coordinates. Let us consider the cotangent bundle $T^*Q$ endowed with the canonical symplectic structure $T^*Q$. Let $H$ be a Hamiltonian and $X_H$ be the associated Hamiltonian vector field. We recall that the cotangent bundle is also endowed with the natural projection $\pi_Q$, and that sections of this (vector) bundle are $1$-forms. Consider now one-form $\gamma$ on $Q$, then we can project the Hamiltonian vector field $X_H$ through $\gamma$. We call the projected vector field $X_H^{\gamma}$, the vector field on $Q$ defined by
  \begin{equation}
 X_H^{\gamma}:Q\longrightarrow TQ, \qquad q\mapsto T\pi_Q \circ X_H \circ \gamma(q),
 \end{equation}
 where $T\pi_Q$ is the tangent mappign fo the cotangent bundle projection. 
 The following commutative diagram determines the projected vector field.
 \begin{equation}\label{geomdiagram}
\xymatrix{ T^{*}Q
\ar[dd]^{\pi} \ar[rrr]^{X_H}&   & &TT^{*}Q\ar[dd]^{T\pi}\\
  &  & &\\
 Q\ar@/^2pc/[uu]^{\gamma}\ar[rrr]^{X_H^{\gamma}}&  & & TQ}
\end{equation}

The geometric Hamilton-Jacobi theory relates the projected vector field $X_H^{\gamma}$ and the Hamiltonian vector field $X_H$, see \cite{carinena2006geometric}. 

\begin{theorem}[The geometric HJ thm: the $t$-independent case]\label{th1}
Under the conditions given above, if $\gamma$ is closed $(d\gamma=0)$ then the following conditions are equivalent:
\begin{enumerate}
\item The vector fields $X_H^{\gamma}$ and $X_H$ are $\gamma$-related that is $T\gamma \circ X_H^{\gamma} = X_H^{\gamma}  \circ \gamma$. 
\item The equation $d(H\circ \gamma)=0$ holds.
\end{enumerate}
\end{theorem}
The first item in the theorem says that if $q(t)$
is an integral curve of $X_{H}^{\gamma }$, then $\gamma\circ q(t)$
is an integral curve of the Hamiltonian vector field $X_{H}$, hence a
solution of the Hamilton's equations \eqref{HamEq}. Such a solution of the Hamiltonian equations is
called horizontal since it is on the image of a one-form on $Q$. 
Since a solution $\gamma$ is assumed to a closed one-form, by Poincar\'{e} lemma, there locally exists a function $S$ satisfying $dS=\gamma$. Substitution of this into the second condition $d\left( H\circ \gamma \right)=0$ results with the classical formulation of the Hamilton-Jacobi problem  where the constant $E$ appears as a manifestation of the integration that is
\begin{equation}\label{HJ1}
H\left(q^i,\frac{\partial S}{\partial q^i}\right)=E. 
\end{equation} 

\textbf{Complete Solutions.}
 A complete solution of the Hamilton--Jacobi equation on a symplectic manifold $(T^*Q,\Omega_Q)$ 
 is a diffeomorphism $S:Q\times \mathbb{R}^n\rightarrow T^{*}Q\times \mathbb{R}^n$ such that for a set of
 parameters $\lambda\in \mathbb{R}^n, \lambda=(\lambda_1,\dots,\lambda_n)$, the mapping
  \begin{equation}\label{compsolHJ}
  S_{\lambda}:Q  \longrightarrow    T^{*}Q , \qquad 
  S_{\lambda}(q)  \mapsto   S(q,\gamma(q,t))
 \end{equation}
\noindent
is a solution of the Hamilton--Jacobi equation. We have the following diagram
 \begin{equation}
\xymatrix{ Q\times \mathbb{R}^n
\ar[dd]^{\alpha} \ar[rrr]^{S}&   & &T^{*}Q\ar[dd]^{f_i} \\
  &  & &\\
 \mathbb{R}^n \ar[rrr]^{\pi_i}&  & & \mathbb{R}}
 \end{equation}
with $\pi_i:\mathbb{R}^n\rightarrow \mathbb{R}$ the projection of $(\lambda_1,\dots,\lambda_n)$ to $\lambda_i$.
We define functions $F_i$ on the cotangent bundle as follows 
\begin{equation}\label{functions8}
 F_i:T^*Q\longrightarrow \mathbb{R},\qquad p\mapsto \pi_i\circ \alpha\circ S^{-1}(p).
\end{equation}
whereas  $\alpha:Q \times \mathbb{R}^n\rightarrow \mathbb{R}^n$ is the projection to the second factor.  
 If $S$ is a complete solution of the Hamilton--Jacobi problem on a symplectic manifold, then the functions defined in \eqref{functions8} commute
 with respect to the canonical Poisson bracket, that is,
 \begin{equation}
  \{F_i,F_j\}=0.
 \end{equation}
 
\textbf{The Time-dependent Case.} In the $t$-dependent case, we assume that both the generating function $S$ and the Hamiltonian function $H$ depend explicitly on the time variable $t$. That is, $S$ is a real valued function on $\mathbb{R}\times Q$ and the time-dependent Hamiltonian is defined on $\mathbb{R}\times T^*Q$. In this case, the HJ equation reads
\begin{equation}
    \frac{\partial S}{\partial t}+H\left(t,q^i,\frac{\partial S(t,q^i)}{\partial t}\right)=0.
\end{equation}
To obtain this equation we need the geometric formalism
on the extended space $\mathbb{R}\times T^*Q$, introduced to deal with $t$-dependent Hamiltonians. Note that $\mathbb{R}\times T^*Q$ cannot be a symplectic manifold because it is odd-dimensional. One can fix the problem by extending the space that is
\begin{equation}
(\mathbb{R}\times T^*Q)\times \mathbb{R}\equiv \mathbb{R}^2\times T^*Q\equiv T^*\mathbb{R}\times T^{*}Q\equiv T^*(\mathbb{R}\times Q),
\end{equation}
which is a symplectic manifold and the extra $\mathbb{R}$ does not modify the dynamics significantly. Assume local coordinates $(t, e, q^i, p_i)$ on the cotangent bundle $T^*(\mathbb{R} \times Q)$ where $e$ is the $t$-conjugated momentum. There is a natural mapping 
\begin{equation}
\pi: T^*(\mathbb{R}\times Q) \longrightarrow \mathbb{R}\times T^{*}Q, \qquad (t, e, q^i, p_i) \mapsto (t, q^i, p_i).
\end{equation}
Referring to a time-dependent Hamiltonian function, we define a Hamiltonin function, so-called extended
Hamiltonian, on the cotangent bundle $T^*(\mathbb{R}\times Q) $ as  \begin{equation}
H^e=\pi^*H+e.
\end{equation}
One can easily see that the vector fields $X_{H^e}$
and $\partial/\partial t + X_{H}$ are $\pi$-related.
As a follow up of Theorem \eqref{th1}, we exhibit the geometric  HJ theory for time-dependent case in the following theorem, see \cite{LeonSar2}.  
\begin{theorem}[The geometric HJ theorem: the $t$-dependent case]\label{thm-t-Ham}
Under the previous considerations, let $\gamma$ be a closed $1$-form on $\mathbb{R}\times Q$. Then the following statements are equivalent:
\begin{enumerate}
\item The vector fields $X_{H^e}^{\gamma}$ and $\partial/\partial t+X_{H^e}$ are $\pi\circ \gamma$-related. 
\item The expression 
\begin{equation}\label{thj}
d(H^e\circ\gamma) \in \langle dt \rangle
\end{equation}
holds.
\end{enumerate}
\end{theorem}

Realize that equation \eqref{thj} is more general than \eqref{HJ1}. Again, using Poincar\'e's Lemma, we can locally write $\gamma=dS$ and write the previous expression in terms of $S$. As particular case, \eqref{thj} includes the equations
\begin{equation}
\frac{\partial S}{\partial t}+ H\left(t,q^i,\frac{\partial S}{\partial q^i}\right)=f(t)
\end{equation}
where $f(t)$ is an arbitrary real-valued real function. 

Indeed, if $H(q^i,p_i)$ is a time-independent Hamiltonian, then, the solution of the time-independent equation \eqref{HJ1} and the time-dependent \eqref{thj} are related in the following way. Let us call the solution $S$ of the time-independent equation now $W$. The solution of the time-dependent is $S=Et+W$ is a solution of the time-dependent Hamilton-Jacobi equation. Usually the function $W$ is called the {\it characteristic function} and $S$ the {\it principal function}.

\textbf{Complete Solutions of the Time-dependent Case.}
About the complete solutions for a time-dependent HJ, we start by interpreting the function $S$ as a function on the product manifold $\mathbb{R}\times Q\times  Q$ and so that image space of its exterior derivative $dS$ is a Lagrangian submanifold in $T^{*}(\mathbb{R}\times Q \times Q)$ (here $(q^i)$ are the coordinates on the first factor $Q$, and $(\alpha_i)$ the coordinates on the second factor $Q$). On the other hand, consider the projections 
\begin{equation}
\pi_I:T^{*}(\mathbb{R} \times Q \times Q) \rightarrow \mathbb{R}\times T^{*}Q,\qquad \pi_I(t,e,\alpha^1,\alpha^2)=(t, (-1)^{I+1}\alpha^I ),
\end{equation}
where $I=1,2$.
With these geometric tools, the non-degeneracy condition is equivalent to saying that $\pi_I |_{\text{Im}(dS)}$ is a local diffeomorphism for $I=1,2$. We assume here for simplicity that it is a global diffeomorphism, so we can consider the mapping 
\begin{equation}
\pi_2|_{\text{Im}(dS)}\circ \left(\pi_1|_{\text{Im} (dS)}\right)^{-1}:T^{*}(\mathbb{R}\times Q)\rightarrow T^{*}(\mathbb{R}\times Q).
\end{equation}
This mapping can be easily checked to be the global description of the change of variables introduced in \cite{Vaqtesis}. The Hamilton-Jacobi equation can be understood as the fact that $dS^{*}H^{\text{ext}}=0$, where $H^{\text{ext}}=\pi_1^{*}H+e$. The diagram below helps us have a global picture of the procedure \cite{Vaqtesis}.

\includegraphics[width=\linewidth]{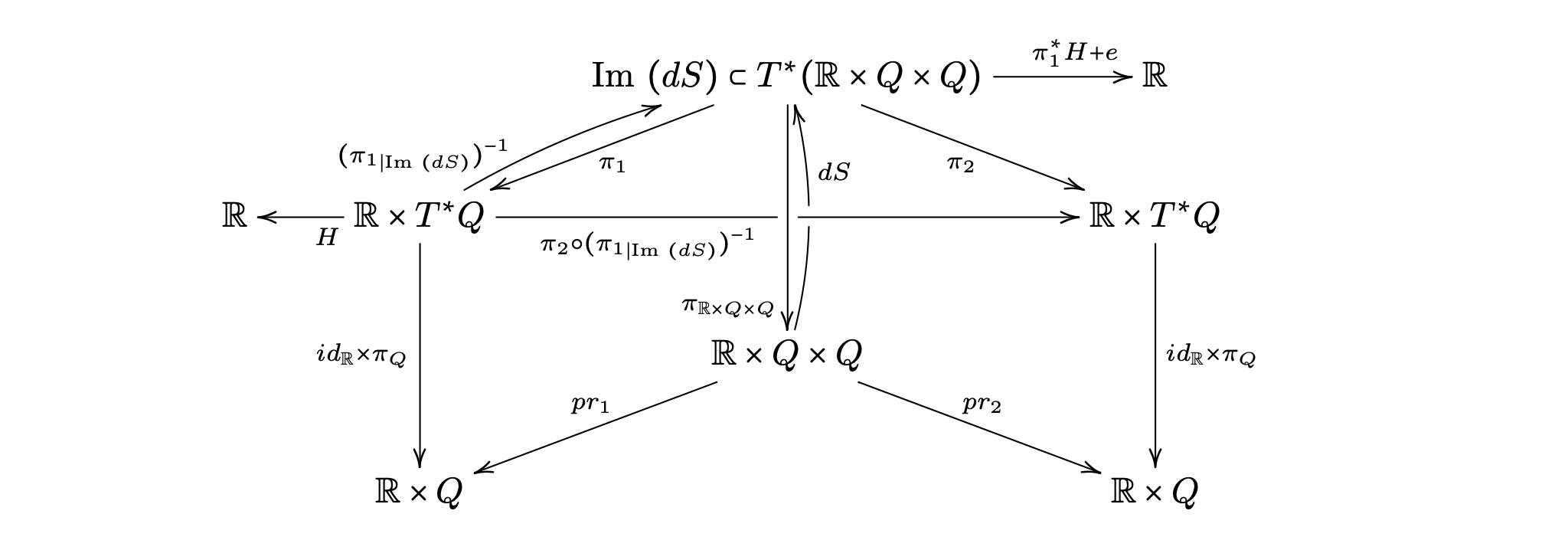}

\subsection{HJ for Conformal Symplectic Dynamics} \label{sec-HJforExt}

We first recall the Hamilton-Jacobi established for Hamiltonian systems admitting external forces. Then we present conformal Hamiltonian dynamics as a particular instance of this construction. Let us now present the geometric formulation of the Hamilton-Jacobi theory for mechanical systems with external forces. Consider Hamiltonian dynamics determined by the triple $(T^*Q,\Omega_Q,H)$. The external forces will be depicted by a semibasic $1-$form $\beta$ (i.e., $\beta$ vanishes over vertical tangent vectors) on $Q$. So, the equation for the dynamics associated with $\beta$ is
\begin{equation}
\iota_{X_{H,\beta}}\Omega_{Q}-dH=\beta
\end{equation} 
Note that $X_{H,\beta}$ is a vector field on $T^*Q$ and it is equal to $X_H$ if $\beta$ vanishes identically. Let $\sigma$ be a closed one-form on $Q$, and introduce a vector field on $Q$ as follows
\begin{equation}
X_{H,\beta}^\sigma:=T\pi_Q\circ X_{H,\beta} \circ \sigma.
\end{equation}
Here is the Hamilton-Jacobi theorem for the present framework.
\begin{theorem} \label{HJforConst}
The following conditions are equivalent:
\begin{itemize}
\item [(i)]$X_{H,\beta}^\sigma$ and $X_{H,\beta}$ are $\sigma$ related.
\item [(ii)]$d\left(H\circ \sigma \right)=-\sigma^*\beta$.
\end{itemize}
\end{theorem}

A $1-$form $\sigma$ on $Q$ satisfying the conditions of  Theorem \ref{HJforConst} will be called a solution for the Hamilton-Jacobi problem given by the Hamiltonian $H$ and the external forces $\beta$. It is evident that if $\beta=0$ then the Hamilton-Jacobi Theorem \ref{HJforConst} reduces to the classical Hamilton-Jacobi Theorem \ref{th1}. We refer \cite{deLeLaLo2021} for a recent study on HJ theory for the systems under external force.

\textbf{HJ for Conformal Hamiltonian Systems.} Consider a symplectic manifold $(M,\Omega)$ A vector field $%
X_c $ on $M$ is called a {\it conformal vector field} \cite{marle2014lie} if it preserves the symplectic two-form up to some scaling factor $c\in%
\mathbb{R}$ that is 
\begin{equation}  \label{conformalvf}
\mathcal{L}_{X_c}\Omega  =c\Omega.
\end{equation}
It is immediate to observe that a conformal vector field reduces to a
Hamiltonian vector field if the scaling factor is zero. The flow $\phi_{t}$
of a conformal vector field preserves the symplectic two-form by the
conformal factor $e^{ct}$, that is 
\begin{equation}
\phi_{t}^{\ast}\Omega  = e^{ct}\Omega. 
\end{equation}

A symplectic manifold $(M,\Omega)$ admits a conformal vector field
with a non-zero scaling factor if and only if the symplectic two-form is the
exterior derivative of a one-form \cite{MaPeSh99}. Assume the cotangent bundle $T^*Q$ equipped with the canonical symplectic two-form $\Omega_Q = -d\Theta_Q$. For a given Hamiltonian
function $H$, the vector field $X_{H}^{c}$ defined through 
\begin{equation}  \label{CHVF}
\iota_{X_{H}^{c}}\Omega_Q = dH - c\Theta_Q
\end{equation}
is conformal. Locally, the inverse of this assertion is also true. That is,
if a vector field is conformal, then there exist a function $H$ and the
conformal vector field can be written in the form of $X_{H}^{c}$ satisfying %
\eqref{CHVF}. More generally, the set of conformal vector fields on $T^{*}Q$ is given by $\{X_H + cZ\}$, where $Z$ is the Liouville vector field
defined by $\iota_Z\Omega_Q = -\Theta_Q$. It is immediate now that generic symplectic structures on cotangent bundles admit conformal Hamiltonian systems, since the canonical symplectic two-form $\Omega_Q$ is minus of the exterior derivative of Lioville one-form. We present this in the following local picture.

In the Darboux coordinates $(q^{i},p_{i})$ on the cotangent bundle $T^*Q$, the Liouville vector
field is $Z=p_{i}\partial /\partial p_{i}$. Hence a
conformal Hamiltonian vector field satisfying (\ref{CHVF}) can be computed to be 
\begin{equation}
X_{H}={\frac{\partial H}{\partial p_{i}}}{\frac{\partial }{\partial q^{i}}}%
+\left(cp_{i}-{\frac{\partial H}{\partial q^{i}}}\right){\frac{\partial }{\partial p_{i}%
}}.  \label{ConfHamVecLoc}
\end{equation}
In this case, the Hamilton equations turn out to be 
\begin{equation} \label{HamEq-con}
\dot{q}^{i}=\frac{\partial H}{\partial p_{i}},\qquad \dot{p}_{i}=cp_{i}-\frac{%
\partial H}{\partial q^{i}}.
\end{equation}
\noindent

\textbf{HJ equation for conformal Hamiltonian systems.} Referring to the symplectic manifold $(T^*Q,\Omega_Q)$, the semi-basic character of the Liouville one-form $\Theta_Q$ enables us to state the following theorem which establishes the Hamilton-Jacobi theorem for conformal systems. 
\begin{theorem}
\label{HJTv2} The following conditions are equivalent:

\begin{enumerate}
\item The vector fields $X_{c}$ and $X_{c}^{\gamma }$ are $\gamma$-related,
that is
\begin{equation}
T\gamma(X_{c}^{\gamma})=X_{c}\circ\gamma.
\end{equation}

\item $
d\left( H\circ \gamma \right) =c\gamma^{\ast }\Theta_{Q}.
$
\end{enumerate}
\end{theorem}
\noindent
{\it Note:} It is easy to see that this is just a  nonholonomic system with a external force, i.e., the case of Chaplygin systems \cite{Balseiro12,deLeMadeDi10}.

\textbf{An Example: A Host-Parasite Model.} Consider $\mathbb{R}^2$ with the coordinates $(x,y)$.
 We consider the host-parasite model described by the following
system of differential equations \cite{EsGu-host}, where $a,b,c,\delta \in \mathbb{R}$
\begin{equation}
\begin{cases}
\dot{x} & =ax-byx \\ 
\dot{y} & =cy-\delta y^{2}/x%
\end{cases}
\label{hp.1}
\end{equation}
If we set $c=0$ in (\ref{hp.1}), we obtain the following reduced system 
\begin{equation}\label{hp.2}
\begin{cases} 
\dot{x} & =ax-byx \\ 
\dot{y} & =-\delta y^{2}/x%
\end{cases}%
\end{equation}
The Hamiltonian function is
computed to be 
\begin{equation}
H(x,y)=-\frac{a}{y}-b\ln y-{\frac{\delta }{x}}, \label{HamHP2-}
\end{equation}
provided that $x$ and $y$ are not zero. 
We write the corresponding symplectic two-form for this case as follows 
\begin{equation}
\Omega =xy^{2}dx\wedge dy, \label{HamHP2--}
\end{equation}
It is now a matter of direct calculation to show that the reduced system (%
\ref{hp.2}) is Hamiltonian with Hamiltonian function (\ref%
{HamHP2-}) and the symplectic two-form $\Omega$ is (\ref{HamHP2--}). We need to add the term $cy$ to the right hand side of the second equation
in (\ref{hp.2}) to study the original system that we propose (\ref{hp.1}). The introduction of the term is performed in a pure geometrical way, by
introducing the Liouville vector field $Z=\frac{1}{3}y\frac{\partial}{\partial y}$. Then, the
complete system (\ref{hp.1}) corresponds with a conformal Hamiltonian system whose dynamics is described by the conformal vector field $X_c$ such that
\[
\iota_{X_c}\Omega =dH-3c\theta ,
\]%
where $X_c=X_{H}+3cZ$, $\Omega$ is the symplectic form in (\ref{HamHP2--}) and that $\theta $
is the associated potential one-form $\theta =\frac{1}{3}y^{3}xdx,$
satisfying $\Omega =-d\theta$.
In the particular case of the host-parasite model, the Hamiltonian vector field and its projection are:

\begin{equation}
X_c=\left(\frac{a}{y^2}-\frac{b}{y}\right)\frac{\partial}{\partial x}+\left(cy-\frac{\delta}{x^2}\right)\frac{\partial}{\partial y},\qquad X_c^{\gamma}=\left(\frac{a}{y^2}-\frac{b}{y}\right)\frac{\partial}{\partial x}
\end{equation}

\noindent
Now, using condition 1 and 2 in Theorem \eqref{HJTv2}, we obtain an expression for $\gamma$
that reads:

\begin{equation}
\frac{\partial \gamma(x)}{\partial x}=-\frac{\delta}{x^2}+c\gamma(x)
\end{equation}

This equation is a linear equation and we can integrate it easily. The general solution is:

\begin{equation}
  \gamma(x)=  \delta  c\left(\frac{1}{cx}-\text{Ei}_1(x)e^{cx}\right)
\end{equation}

where $Ei_1$ is the exponential integral.

\subsection{Cosymplectic Manifolds and Dynamics}\label{Sec-Alm-cos}

To study $t$-dependent (or parameter dependent Hamiltonians)
we need to introduce {almost cosymplectic structures} and their application in classical Hamiltonian mechanics.
By an {\it almost cosymplectic structure} we understand a $2n+1$-dimensional manifold equipped with a one-form $\eta$ and a two-form $\Omega$ such that $\eta\wedge \Omega^n$ is a volume form.
In particular, we will study two cases of almost cosymplectic manifolds. On one hand, the case of {cosymplectic manifolds} \cite{Blair, LeonTuyn}, and on the other hand, the case of
{contact manifolds} \cite{Blair,boothwang,Godbillon}. Cosymplectic manifolds are a particular type of Poisson manifolds, and contact manifolds are a particular case of Jacobi manifolds. We will focus on cosymplectic manifolds in this section, and we will come back to contact manifolds in the third part of this manuscript.
Cosymplectic manifolds have shown their usefulness in theoretical physics, as in gauge theories of gravity, branes and string theory \cite{Becker,Cham,Hitchin}. 
Among the early studies of cosymplectic manifolds we mention A. Lichnerowicz \cite{Lich,Lichnerowicz-Jacobi}, who studied the Lie algebra of infinitesimal automorphisms of a 
cosymplectic manifold, in analogy with the symplectic case. 
Posteriously, some works have endowed cosymplectic manifolds with a Riemannian metric, the so-called {coK\"ahler manifolds} \cite{oku}. These
are the odd dimensional counterpart of K\"ahler manifolds.
Another important role of the cosymplectic theory is the reduction theory to reduce time-dependent Hamiltonians by symmetry groups \cite{Albert}.
But since very foundational papers by P. Libermann, very sporadic papers have appeared on cosymplectic settings \cite{Libermann-LCS,Takizawa}. It is our future intention to provide surveys on cosymplectic geometry due to their lack \cite{Cape}.
Our particular interest in cosymplectic structures resides in their use in the description of time-dependent mechanics. 

\textbf{(Almost) Cosymplectic Dynamics.} 
An almost cosymplectic structure on a manifold $M$ of odd dimension $2n + 1$ equipped with a one-form $\eta$ and a two-form $\Omega$ satisfying that $\eta\wedge\Omega^n$ is a non-vanishing top-form. This permits us to define an isomorphism from the space of vector fields $\mathfrak{X}(M)$ to the space of one-forms $\Gamma^1(M)$ given by
\begin{equation}\label{flat-map}
    \flat: \mathfrak{X}(M)\longrightarrow \Gamma^1(M),\qquad X\mapsto \flat(X)=\iota_X\Omega + \eta(X)\eta.
\end{equation}
We denote the inverse of this isomorphism by $\sharp$. The inverse image of the one-form $\eta$ is a vector field denoted by $\mathcal{R}$ and called as Reeb field. It is straightforward to see that the Reeb field satisfies the following contractions
\begin{equation}\label{Reeb-cosymp}
\iota_{\mathcal{R}}\eta=1,\qquad \iota_{\mathcal{R}}\Omega=0.
\end{equation}
The inverse of this assertion is also true. That is, consider a manifold $M$ equipped with a pair $(\eta,\Omega)$ establishing the flat map $\flat$ in \eqref{flat-map} as an isomorphism. Assume also a distinguished vector field $\flat$ satisfying \eqref{Reeb-cosymp}, then  $M$ turns out to be $2n + 1$-dimensional, and the pair $(\eta,\Omega)$ determines $M$ as an almost cosymplectic manifold. In this case, the almost cosymplectic manifold admits $\mathcal{R}$ as the Reeb field.

Consider an almost cosymplectic manifold $(M,\eta,\Omega)$. Referring to the the flat map $\flat$, for a real valued function $H$ on $M$, we assign a vector field, called gradient vector field, $\operatorname{grad} H$ on $M$ as follows
\begin{equation}
\flat\circ \operatorname{grad} H=dH.
\end{equation}
It is evident that the gradient field $\operatorname{grad} H$ is unique for a given $dH$ and vice versa. Accordingly, we define a bracket of functions on $M$ as 
\begin{equation}\label{Poisson-bra-coi}
\{F,H\}:=\Omega(\operatorname{grad} F,\operatorname{grad} H).
\end{equation}
It is evident to see that this bracket is skew-symmetric and satisfies the Leibniz identity. But, in general it fails to satisfy the Jacobi identity. So that an almost cosymplectic manifold is an almost Poisson manifold. We can endow the Hamiltonian vector field definition in \eqref{Hamvf}. In accordance with this, for a Hamiltonian function $H$ on $M$ the associated Hamiltonian vector field $X_H$ is computed to be
\begin{equation}
\iota_{X_H}\Omega =dH- \mathcal{R}(H)\eta,\qquad \iota_{X_H}\eta=0.
\end{equation}
The difference between the gradient vector field and the Hamiltonian vector field is computed to be 
\begin{equation}
\operatorname{grad} H-X_H= \sharp(\mathcal{R}(H)\eta).
\end{equation}
Additionally one can define evolutionary vector field as the linear combination of the Hamiltonian vector field and the Reeb field that is 
\begin{equation}
E_H=\mathcal{R}+X_H.
\end{equation}

An almost cosymplectic manifold $(M,\eta,\Omega)$ is called cosymplectic if both $\Omega$ and $\eta$ are closed. In this case, it becomes to a matter of a direct calculation to show that the bracket \eqref{Poisson-bra-coi} satisfies the Jacobi identity as well. So, we can claim that it is a Poisson bracket.

%
%
\subsection{HJ for Cosymplectic Dynamics}\label{Sec-HJ-cos}
Consider the extended phase space $T^{*}Q\times \mathbb{R}$ and its canonical projections of the first and second factor, $\rho:T^{*}Q\times \mathbb{R}\rightarrow T^{*}Q$ and  $t:T^{*}Q\times \mathbb{R}\rightarrow \mathbb{R}$, respectively
 and a time-dependent Hamiltonian $H:T^{*}Q\times \mathbb{R}\rightarrow \mathbb{R}$. Let us depict the problem with a diagram
\[
\xymatrix{ T^{*}Q\times \mathbb{R}
\ar[dd]^{\rho}  \ar[ddrr]^{t}  \ar@/^2pc/[ddrr]^{H}  &   &\\
  &  & &\\
T^{*}Q & & \mathbb{R}}
\]
We have canonical coordinates $\{q^i,p_i,t\}$ with $i=1,\dots,n$, where $(q^i,p_i)$ are fibered coordinates in $T^{*}Q$ and $t\in \mathbb{R}$. In Draboux coordinates, the pair $\rho^*\Omega_Q$ and $dt$ determines a cosymplectic pair. In this case the Reeb field is  $\mathcal{R}=\partial/\partial t$ whereas the Hamiltonian and gradient fiedls are  computed to be
\begin{equation}
X_H=\frac{\partial H}{\partial p_i}\frac{\partial}{\partial q^i}-  \frac{\partial H}{\partial q^i}\frac{\partial}{\partial p_i},\qquad \operatorname{grad} H=\frac{\partial H}{\partial t}\frac{\partial}{\partial t}+X_H,
\end{equation}
respectively. A direct calculation gives the evolutionary vector field as
\begin{equation}\label{evo-field-cos}
E_H=\frac{\partial}{\partial t}+\frac{\partial H}{\partial p_i}\frac{\partial}{\partial q^i}-  \frac{\partial H}{\partial q^i}\frac{\partial}{\partial p_i}.
\end{equation}

To have an alternative coymplectic pair on the extended cotangent bundle, we consider the two-form on $T^{*}Q\times \mathbb{R}$ as $\Omega_H=-d\theta_H$
and 
\begin{equation}\label{thetah}
\theta_H=\Theta_Q-Hdt
\end{equation}
where $\Theta_Q$ is the canonical Liouville one-form. 
We abuse notation by identifying the pullbacks of the one-forms with the one-forms themselves. That is, $\rho^{*}(\Theta_Q)=\Theta_Q$.
Hence, 
\begin{equation}\label{oh}
\Omega_H=\sum_{i=1}^n dq^i\wedge dp_i+dH\wedge dt.
\end{equation}
Let us consider the cosymplectic structure $(\eta=dt,\Omega_H).$
The corresponding Reeb vector field needs to satisfy 
\begin{equation}\label{reebc}
\iota_{\mathcal{R}_H}dt=1,\quad \iota_{\mathcal{R}_H}\Omega_H=0.
\end{equation}
The unique Reeb vector field that satisfies \eqref{reebc} has the following expression in coordinates
\begin{equation}\label{reebcosym}
 \mathcal{R}_H=\frac{\partial}{\partial t}+ \frac{\partial H}{\partial p_i}\frac{\partial}{\partial q^i}-  \frac{\partial H}{\partial q^i}\frac{\partial}{\partial p_i}.
\end{equation}
A direct comparison shows that the evolutionary field $E_H$ in \eqref{evo-field-cos} and the Reeb field $\mathcal{R}_H$ are coinciding. For both of these vector fields, the corresponding Hamilton's equations are
\begin{equation}\label{hamileq22} 
 {\dot q}^i =\frac{\partial H}{\partial p_i},\qquad 
 {\dot p}_i =-\frac{\partial H}{\partial q^i},\qquad  {\dot t}=1. 
 \end{equation}
\noindent
Since $\dot{t}=1$, we can consider $t$ a time-parameter (up to an affine change).

\noindent
We consider the fibration
$\pi: T^{*}Q\times \mathbb{R}\rightarrow Q\times \mathbb{R}$ and a section $\gamma$ of $\pi:T^{*}Q\times \mathbb{R} \rightarrow Q\times \mathbb{R}$, i.e., $\pi\circ \gamma=\text{id}_{Q\times \mathbb{R}}$. 
Also, we assume that $\text{Im}(\gamma_t)$ with $\gamma_t:Q\rightarrow T^{*}Q\times \mathbb{R}$ such that
$\gamma_t(q^i)$ in coordinates $(q^i,\gamma^i(q^i,t),t)$ is a lagrangian submanifold for a fixed time $t$ of the cosymplectic manifold $(T^{*}Q\times \mathbb{R},dt,\Omega_H)$ for a fixed time, that is $d\gamma_t=0$.
%

\[
\xymatrix{
 T^{*}Q\times \mathbb{R}\ar[rr]^{\rho} & &  T^{*}Q\ar[rr]^{\pi} & &  Q\times \mathbb{R}\ar[rr] & &  Q\ar@/^2pc/[llllll]^{\gamma_t}}
\]

We can use $\gamma$ to project $\mathcal{R}_H$ on $Q\times \mathbb{R}$
just defining a vector field $\mathcal{R}^{\gamma}_H$, the denominated {\it projected} vector field on $Q\times \mathbb{R}$ by
\begin{equation}\label{hjr}
 \mathcal{R}^{\gamma}_H=T_{\pi}\circ \mathcal{R}_H\circ \gamma
\end{equation}
The following diagram summarizes the above construction
\[
\xymatrix{ T^{*}Q\times \mathbb{R}
\ar[dd]^{\pi} \ar[rrr]^{\mathcal{R}_H}&   & &T(T^{*}Q\times \mathbb{R})\ar[dd]^{T_{\pi}}\\
  &  & &\\
Q\times \mathbb{R} \ar@/^2pc/[uu]^{\gamma}\ar[rrr]^{\mathcal{R}^{\gamma}_H}&  & & T(Q\times \mathbb{R})}
\]
\begin{definition}
 If $\alpha$ is a one-form, locally expressed as $\alpha=\sum_{i=1}^n\alpha_i dq^i$, we designate by $\alpha^{V}$ the {\it vertical lift} \cite{yanoishi}
or vector fields associated with $\alpha$, defined by
\begin{equation*}
 \iota_{\alpha^V}\Omega_Q=\alpha
\end{equation*}
Hence, the vector field $\alpha^V$ has the local expression
\begin{equation}
\alpha^{V}=-\sum_{i=1}^n\alpha_i\frac{\partial}{\partial p_i}.
\end{equation}
\end{definition}

\begin{theorem}
 The vector fields $\mathcal{R}_H$ and $\mathcal{R}^{\gamma}_H$ are $\gamma$-related if and only if the following equation is satisfied
\begin{equation}\label{eqtheorem2}
[d(H\circ \gamma_t)]^{V}=\dot{\gamma}_q
\end{equation}
where $[\dots]^{V}$ denotes the vertical lift of a one-form on $Q$ to $T^{*}Q$.  Now $\dot{\gamma}_q$ is the tangent vector in a point $q$ associated with the curve
\[
\xymatrix{
 \mathbb{R}\ar[rr]\ar@/^2pc/[rrrrrr]^{\gamma_q} & &  Q\times \mathbb{R}\ar[rr] & &  T^{*}Q\times \mathbb{R}\ar[rr]^{\rho} & &  T^{*}Q}
\]
Notice that these applications are given for a fixed point $q\rightarrow (q,t,\gamma)$.

\end{theorem}

Equation \eqref{eqtheorem2} is known as a {\it Hamilton--Jacobi equation on a cosymplectic manifold.}
In local coordinates $\{q^i,p_i,t\}$, we have
\begin{equation}\label{lc}
 \frac{\partial \gamma^j}{\partial t}+  \frac{\partial H}{\partial p_i}\frac{\partial \gamma^j}{\partial q^i}+\frac{\partial H}{\partial q^j}=0.
\end{equation}

\textbf{An Example: the damped Hamiltonian.} 
A damped Hamiltonian comes from a Lagrangian for a harmonic oscillator that has a time-dependent mass of the from $m=m_0e^{\Gamma t}$ accreting with time in order to mimic exponential energy dissipation. This Hamiltonian reads:
\begin{equation}\label{dampham}
    H=\frac{p^2}{2m}e^{-\Gamma t}+\frac{m}{2}\Omega_0^2e^{\Gamma t}x^2
\end{equation}
For this example the Reeb vector field takes the form:
\begin{equation}
    \mathcal{R}_H=\frac{\partial}{\partial t}+\frac{p}{m}e^{-\Gamma t}\frac{\partial}{\partial x}-m\Omega_0^2e^{\Gamma t}x\frac{\partial}{\partial p}
\end{equation}
and the solution of the HJ will be a section $\gamma$ that in fibered coordinates reads $\gamma=(t,x,\gamma(x))$ defined on the cosymplectic manifold $(\mathbb{R}^{+}\times T^{*}Q, \theta_H,\Omega_H)$ given in \eqref{thetah} and \eqref{oh}. 
So, the HJ equation is:
\begin{equation}
    \gamma \frac{\partial \gamma}{\partial x}=-m^2\Omega_0^2e^{2\Gamma t}x
\end{equation}
and its solution by partial integration reads
\begin{equation}
    \gamma^2=-m^2\Omega_0^2e^{2\Gamma t}x^2+F(t)
\end{equation}

\section{HJ on Jacobi Manifolds}\label{Sec-Jacobi}

In this section we present the HJ theory for dynamical systems whose  evolution is described by one particular case of Jacobi structure, as it is the case of locally conformally symplectic and contact structures. Systems that accord to these descriptions are: systems that locally behave as symplectic dynamical systems but not globally (for example, intertwined dynamical systems) and systems with dissipation, respectively. These two structures have in common that they obey the Jacobi identity for integrability. Let us start by describing generally a Jacobi structure.

\subsection{Jacobi Manifolds and Dynamics}\label{Sec-Jac-Dyn}

    A Jacobi manifold~\cite{Kirillov,Lichnerowicz-Poi,Lichnerowicz-Jacobi,Marle-Jacobi} is a triple $(P,\Lambda,Z)$, where $\Lambda$ is a bivector field and $Z$ is a vector field, such that the following identities are satisfied
\begin{equation}\label{ident-Jac}
        [\Lambda,\Lambda] = 2 Z \wedge \Lambda, \qquad 
       [Z,\Lambda] = 0,
\end{equation}
    where $[\bullet,\bullet ]$ is the Schouten-Nijenhuis bracket.

    Given a Jacobi manifold $(P,\Lambda,Z)$, we define the {Jacobi bracket} of two functions $H,F\in C^{\infty}(P)$ as
   \begin{equation}\label{bra-Jac}
        \{H,F\} = \Lambda(dH, dF) + H\cdot Z(F) - F\cdot Z (H).
 \end{equation}
This bracket is bilinear, antisymmetric, and satisfies the Jacobi identity. The bilinearity and the antisymmetry properties of the bracket \eqref{bra-Jac} follow from its definition. The satisfaction of the Jacobi identity is a manifestation of the identities in \eqref{ident-Jac}.  Furthermore, it fulfills the weak Leibniz rule:
\begin{equation}
        \operatorname{supp}(\{F,G\}) \subseteq \operatorname{supp} (F) \cap \text{supp} (G).
\end{equation}
This characteristic of the Jacobi bracket reads the Lie algebra $(C^\infty(M), \{\bullet,\bullet\})$ as a local Lie algebra in the sense of Kirillov, \cite{Kirillov}. The inverse of this assertion is also true, that is, a local Lie algebra determines a Jacobi structure. 
Consider two Jacobi manifolds, say $P_1$ and $P_2$. A smooth mapping $\phi:P_1\mapsto P_2$ is said to be a Jacobi map if it preserves the Jacobi bracket introduced in \eqref{bra-Jac}. 

Notice that the Jacobi bracket reads the product of two functions in one of its entries as 
\begin{equation}
\{F, GH\} = G \{F, H\} + H \{F, G\} + G\cdot H\cdot Z(F).
\end{equation}
It is evident that for the satisfaction of the Leibniz identity the vector field $Z$ must vanish identically. This is an if and only if argument. See that, in this case, the first identity in \eqref{ident-Jac} reduces to the Poisson condition in \eqref{Poisson-cond}, whereas the second one is trivially satisfied. To sum up, we argue that $(P,\Lambda,Z)$ is a Poisson manifold if and only if $Z=0$.

\textbf{Hamiltonian Vector Field.} Consider a Jacobi manifold  $(P,\Lambda,Z)$ and assume a function $H$ on $P$. The vector field 
\begin{equation}\label{Ham-v-f-Jac}
X_H=\Lambda^\sharp(dH)+HZ
\end{equation}
is called a Hamiltonian vector field of $H$. Here, $\Lambda^\sharp$ stands for the musical mapping generated by the bivector field $\Lambda$, whereas $dH$ is the de Rham (exterior) derivative of $H$. In this picture, $H$ is called the Hamiltonian function. It is interesting here to note that the Hamiltonian vector field corresponding with the constant function $1$ is precisely $Z$. A direct calculation proves that the mapping taking a Hamiltonian function $H$ to the Hamiltonian vector field $X_H$ is a Lie algebra homomorphism. 

\textbf{Foliation of Jacobi Manifolds and Subcases.} The characteristic distribution of a Jacobi manifold $P$ is the subbundle of $TP$ generated by the Hamiltonian vector fields. Referring to the definition \eqref{Ham-v-f-Jac}, for a point $p$ in $P$, the characteristic distribution is determined by the following sum 
\begin{equation}
\mathfrak{C}_p=\Lambda^\sharp(T^*_p P)+\langle Z_p \rangle,
\end{equation}
where $\langle Z_p \rangle$ is a one-dimensional space spanned by the value of $ Z$ at $p$. 
A Jacobi manifold is called transitive if its characteristic distribution is equal to $TP$. An even dimensional transitive Jacobi manifold is a locally conformally symplectic manifold whereas an odd dimensional transitive Jacobi manifold is a contact manifold. We shall examine these two cases in detail in the upcoming sections. The characteristic distribution of a Jacobi manifold is integrable. Therefore, an even dimensional Jacobi manifold foliates into a collection of locally conformally symplectic (LCS) manifolds, whilst an odd dimensional Jacobi manifold foliates into a collection of contact manifolds. 
We draw the following table summarizing these observations. 

\begin{table}[H]{\footnotesize
  \noindent
\caption{{\small {\bf Jacobi Manifolds}. $M$ is a manifold, $Z$ is a vector field and $\Lambda$ is a bivector field.  $\Lambda^\sharp$ is the musical map induced by $\Lambda$. $[\bullet,\bullet]$ is the Schouten bracket. $\Omega^n$ is the $n$-th wedge power of $\Omega$. }}
\label{table2}
\medskip
\noindent\hfill
\resizebox{\textwidth}{!}
{\begin{minipage}{\textwidth}
\begin{tabular}{ l l l l }
 \hline
 &&\\[-1.5ex]
 Structure&  Characterization &Bracket and h.v.f.&  (Bi)vector fields   \\[+1.0ex]
\hline
 &  & \\[-1.5ex]
 {\bf Jacobi}& $[Z,\Lambda] = 0$   & $\{F,G\}=\Lambda(dF,dG)+F  Z(G)-G  Z(F) $ &    \\[+1.0ex] $(M,\Lambda, Z)$& $[\Lambda,\Lambda]=2Z\wedge \Lambda$ &   $X_H=\Lambda^\sharp (dH)+HZ$  & \\[+1.0ex]  &&$Z=0\quad $ Poisson  &  \\[+2.0ex] 
\hline 
 &  & \\[-1.5ex]
 {\bf LCS}&  $d\theta=0,\quad \Omega^n\neq 0$ & Same as Jacobi &  $\Lambda(\alpha,\beta)=\Omega(\sharp(\alpha),\sharp(\beta))$  \\[+1.0ex] $(M,\Omega, \theta)$ & $d\Omega=\theta \wedge \Omega$ &     & $Z_\theta=\Omega^\sharp(\theta)$   \\[+1.0ex]even $2n-$dim && &  \\[+2.0ex]
\hline\\[0.5ex]
 {\bf Contact} & $\eta\wedge {(d\eta)}^n\neq 0$  & Same as Jacobi   &   $\Lambda(\alpha,\beta)=-d\eta(\sharp(\alpha),\sharp(\beta))$ \\[+1.0ex]$(M,\eta)$  &   & &  $Z=-\mathcal{R}=-\sharp(\eta)$ \\[+1.0ex]odd $(2n+1)-$dim &&&   \\[+2.0ex]
\hline\\[0.5ex]
\end{tabular}
\end{minipage}}
\hfill}
\end{table}

\subsection{LCS Manifolds and Dynamics}

Let $M$ be a manifold. 
A non-degenerate two-form $\Omega$ on $M$ is said to be locally conformally symplectic abbreviated as LCS if the two-form is closed locally up to a conformal parameter i.e., if around each point $x$ in $M$ there exists an open neighborhood, say $U_\alpha$, , and a function $\sigma_{\alpha}:U_\alpha\to\mathbb R$ such that the exterior derivative $d(e^{-\sigma_{\alpha}}\Omega\vert_{\alpha})$ vanishes identically on $U_\alpha$, see \cite{Vaisman85}.  
Here, $\Omega\vert_{\alpha}$ denotes the restriction of the two-form $\Omega$ to $U_\alpha$. We refer also some more recent studies on LCS and the dynamics on this geometry, \cite{Banyaga2002,Bazzoni2018,HallerRybicki,Stanciu2019,WojtkowskiLiverani1998}

Collecting the local conformal factors, we define the Lee one-form $\theta$ on $M$ such that $\theta\vert_{\alpha}=d\sigma_{\alpha}$ \cite{LeeHC}. This realization permits us to denote a LCS manifold by a triple $(M,\Omega,\theta)$.   A direct calculation gives that $d\Omega=\theta\wedge \Omega$.  Since $\theta$ is locally exact, it is also closed.  A LCS manifold $(M,\Omega,\theta)$ is a globally conformally symplectic GCS manifold if the Lee form $\theta$ is an exact one-form.  

\noindent 

\noindent \textbf{Musical Mappings.} Consider a LCS manifold $(M,\Omega,\theta)$. The non-degeneracy of the two-form $\Omega$ leads us to define a musical isomorphism 
\begin{equation} \label{mus-iso}
\Omega^\flat:\mathfrak{X}(M)\longrightarrow \Gamma^1(M): X \mapsto \iota_X\Omega,
\end{equation}
where $\iota_X$ is the interior derivative. Here, $\mathfrak{X}(M)$ is the space of vector fields on $M$, whilst $\Gamma^1(M)$ is the space of one-form sections on $M$. 
The image of the Lee form $\theta$ under this isomorphism 
determines the Lee vector field 
\begin{equation} \label{Lee-v-f}
Z_\theta:=\Omega^\sharp(\theta), \qquad \iota_{Z_\theta}\Omega=\theta.
\end{equation}
Notice that $\mathcal{L}_{Z_\theta}\theta=0$ and that $\mathcal{L}_{Z_\theta}\Omega=0$.

\noindent \textbf{The Lichnerowicz-de Rham Differential.}
Referring to a one-form $\theta$ on a manifold $M$, the Lichnerowicz-de Rham differential (LdR) on the space of differential forms is defined to be \cite{GuLi84}
\begin{equation} \label{LdR-Diff}
d_\theta: \Gamma^k(M) \rightarrow \Gamma^{k+1}(M) : \beta \mapsto d\beta-\theta\wedge\beta.
\end{equation}
See that, $d_\theta$ is a differential operator of order $1$ that transforms a $k$-form into a $k+1$-form. If $\theta$ is closed, then one has $d_{\theta}^2=0$. This way, one can study the cohomology of $d_\theta$ \cite{HaRy99}. It is immediate to see that a triple $(M,\Omega,\theta)$ is a LCS manifold if and only if  $d_\theta \Omega = 0$.

\textbf{LCS on Cotangent Bundles.} 
Consider  the canonical symplectic manifold $(T^*Q,\Omega_Q)$, and assume a closed one-form  $\vartheta$ on the base manifold $Q$. The pull back of $\vartheta$ to $T^*Q$ by means of the cotangent bundle projection $\pi_Q$ determines a one-form $\theta=\pi_Q^*(\vartheta)$. Notice that $\theta$ is  closed and semi-basic. Referring to the Lichnerowicz-de Rham differential, we define a two-form 
\begin{equation} \label{omega_theta}
\Omega_\theta=-d_\theta(\Theta_Q)= -d\Theta_Q+\theta\wedge \Theta_Q=\Omega_Q+\theta\wedge \Theta_Q
\end{equation}
on the cotangent bundle $T^*Q$. It is evident that
\begin{equation}
d\Omega_\theta=\theta\wedge \Omega_\theta.
\end{equation} 
We denote this LCS structure by the triple $(T^*Q,\Omega_\theta,\theta)$ or shortly by $T^*_\theta Q$. 
Let us note that $T^*_\theta Q
$ is a generic example for LCS manifolds, which means that all LCS manifolds locally look like $T^*_\theta Q$ with a closed one-form $\vartheta$, see \cite{ChantraineMurphy2019,OtimanStanciu2017}. 
\bigskip

\noindent 

\textbf{Lagrangian Submanifolds.} Consider $(T^*Q,\Omega_\theta,\theta)$ and a one-form $\gamma$ on $Q$. A direct computation shows that the pull-back of the LCS structure is $d_\theta$ exact, that is 
\begin{equation} \label{LagSubT*Q}
\gamma^* \Omega_\theta = - d_\vartheta \gamma,
\end{equation} 
where $d_\vartheta$ denotes the LdR differential defined by the one-form $\vartheta$ on $Q$. This implies that the image space of $\gamma$ is a Lagrangian submanifold of $T_\theta^*Q$ if and only if $d_\vartheta \gamma=0$. Since $d_\vartheta^2$ is identically zero, the image space of the one-form $d_\vartheta f$ is a
 Lagrangian submanifold of $ (T^*Q,\Omega_\theta,\theta)$ for some function $f$ defined on $Q$. 
 
\textbf{Dynamics on LCS Manifolds } Consider a LCS manifold $(M,\Omega,\theta)$. For a Hamiltonian function $H$, the Hamiltonian vector field is defined to be 
\begin{equation}\label{semiglobal3}
\iota_{X_H} \Omega = d_\theta H,
\end{equation}
where $d_\theta$ is the Lichnerowicz-de Rham differential in (\ref{LdR-Diff}). In terms of the Lee vector field $Z_\theta$, the Hamiltonian vector field is computed to be 
\begin{equation}\label{vflcs}
X_H=\Omega^\sharp(dH)+HZ_\theta. 
\end{equation}
Notice that if $H = 1$, then the Hamiltonian vector field is $Z_\theta$, which is far from being zero. Now, we can write $Z_\theta=X_1$. More generally, a vector field $X$ is called a locally Hamiltonian vector field if 
\begin{equation}
d _ { \theta } ( \iota_{X}\Omega)=0.
\end{equation}
A Hamiltonian vector field $X_H$ is locally Hamiltonian since $d^2_\theta =0$.

\subsection{HJ for LCS Dynamics}\label{Sec-lcs-man}

Consider a LCS manifold $ (T^*Q,\Omega_\theta,\theta)$ and a Hamiltonian vector field $X_H$ defined by (\ref{semiglobal3}). For a section $\gamma$ of the cotangent bundle, we define a vector field $ X_{H}^{\gamma}$ on $Q$ by the following commutative diagram
\[
\xymatrix{ T_{\theta}^*Q
\ar[dd]^{\pi_Q} \ar[rrr]^{X_{H}}&   & &TT_{\theta}^*Q  \ar[dd]^{T\pi_Q}\\
  &  & &\\
Q \ar@/^2pc/[uu]^{\gamma}\ar[rrr]^{X_{H}^{\gamma}}&  & & TQ}
\]
as $X_{H}^{\gamma}:=T\pi\circ X_{H}\circ\gamma$. 
\begin{theorem}
\label{HJT} Consider a one-form $\gamma$ whose image is a Lagrangian submanifold of the locally conformally symplectic manifold $T^*_\theta Q$ with respect to the almost symplectic two-form $\Omega_\theta$, that is $d_\vartheta \gamma=0$. Then, the following conditions are equivalent:
\begin{enumerate}
\item The vector fields $X_{H}$ and $X_{H}^{\gamma}$ are $\gamma$-related,
that is
\begin{equation}
T\gamma \circ X_{H}^{\gamma}=X_H\circ\gamma.
\end{equation}
\item  The following equation holds
\[
d_\vartheta (H\circ\gamma)=0.
\]
\end{enumerate}
\end{theorem} 

Extensions of HJ theory for LCS geometry to the case of k-symplectic formalism and fiber bundles have been studied in  \cite{EsLeSaZa-k-sympl} and \cite{EsLeSaZa-Cauchy}, respectively. These works are establishing locally conformal approach to the field theoretical problems.  

\textbf{Example: Gaussian isokinetic dynamics.}


In some cases the LCS structure appears naturally when we consider systems of equations restricted to submanifolds of the phase space \cite{WojtkowskiLiverani1998}. In physics an example of such a situation is provided by Gauss isokinetic dynamics. Let us consider an $n$-dimensional Euclidean space $Q=\mathbb{R}^n$ with a standard flat metric. We will denote by $|\vec v|$ a square of a vector $\vec v=(v^1,...,v^n)$, i.e. $|\vec v|^2=\delta_{ij}v^jv^i$, where $\delta_{ij}$ is a Kronecker delta. The Gaussian isokinetic dynamics is described by the system of equations 
\begin{equation}\label{gas-equation}
\delta_{ij}\dot q^j=p_i,  \quad \dot p_i= F_i-\alpha p_i, 
\end{equation}
where $\alpha=\frac{F_jp_j}{|\vec p|^2}$. We will assume that the force field $F$ is potential, i.e., there exists a local potential function $U:\mathbb{R}^n\rightarrow\mathbb{R}$ such that $F_i=-\frac{\partial U}{\partial q^i}$. Let us introduce a one-form $\widetilde\theta=-\frac{F_i}{|\vec p|^2}\d q^i$ and the Hamiltonian $h=\frac{1}{2}|\vec p|^2-c$. Notice that $\d\widetilde\theta\neq 0$. We define now a two-form
\begin{equation} \label{Ex2-Omega}
\widetilde\Omega_\theta=\Omega_Q-\frac{F_i}{|\vec p|^2} \d q^i\wedge \Theta_Q=\d q^i\wedge\d p_i-p_j\frac{F_i}{|\vec p|^2}\d q^i\wedge\d q^j  .
\end{equation}
It is a matter of straightforward calculations to see that for $X=\dot q^i\partial_{q^i}+\dot p_j\partial_{p_j}$, the equation $\widetilde\Omega_\theta(X,\cdot)=\d h$ is equivalent to (\ref{gas-equation}).
However, the form $\widetilde\Omega_\theta$ is not a locally conformal form since $\d\widetilde\Omega_\theta\neq\widetilde\theta\wedge\widetilde\Omega_\theta$.

In physics, we are often interested in the restrictions of the Hamiltonian equations to a level-set of the Hamiltonian. Let us consider then the subset \mbox{$M:=h^{-1}(0)\subset T^*\mathbb{R}^n$}. Furthermore, we define new forms $\theta$, $\Omega_\theta$ on $Q$ given by
\begin{equation}
\theta=-\frac{F_i}{2c}\d q^i, \quad \quad \Omega_\theta= \d q^i\wedge\d p_i-p_j\frac{F_i}{2c}\d q^i\wedge\d q^j. 
\end{equation}
Notice that $\theta=\d\left(\frac{U}{2c}\right)$ and $d\Omega_\theta=\theta\wedge\Omega_\theta$. This means that the triple $(T^*Q,\Omega_\theta,\theta)$ is a globally conformal symplectic manifold. Moreover, on $M$ we have $\Omega_\theta(X,\cdot)=\d_\theta h$ which means that we can analyze the restriction of equations (\ref{gas-equation}) to $M$ in terms of a LCS structure $(T^*Q,\Omega_\theta,\theta)$. One can also see that the form $\theta$ is a pull-back of the form $\vartheta=-\frac{F_i}{2c}\d q^i$ on $Q$.

Let us apply now the global HJ Theorem \ref{HJT} to the present example. We start with a one-form $\gamma(q^i)=\beta_k(q^i)dq^k$ satisfying the condition $d_\vartheta \gamma=0$. In coordinates, it reads
$$ d_\vartheta \gamma=\Big( \frac{\partial\beta_k}{\partial q^i} + \frac{F_i}{2c}\beta_k(q^i) \Big) \d q^i\wedge dq^k $$
and the condition $d_\vartheta \gamma=0$ is equivalent to
\begin{equation}\label{gammacondition}
\frac{\partial\beta_k}{\partial q^i}+  \frac{F_i}{2c}\beta_k= \frac{\partial\beta_i}{\partial q^k}+  \frac{F_k}{2c}\beta_i.
\end{equation}
Now we consider the vector fields
$$ X_H\circ\gamma= \beta_i\partial_{\dot q^i}+(F_i-\alpha \beta_i)\partial_{p_i}   $$
and 
$$  X_{H}^{\gamma}= \beta_i(q)\partial_{\dot q^i},  \qquad \gamma_* X_{H}^{\gamma}=\beta_i\partial_{\dot q^i}+\frac{\partial \beta_i}{\partial q^k}\beta_k\partial_{p_i}. $$
From the above, we have that the first condition in Theorem \ref{HJT} reads
$$F_i-\alpha \beta_i =  \frac{\partial \beta_i}{\partial q^k}\beta_k, 
\qquad  \alpha=\frac{F_j\beta_j}{2c}.  $$
Referring to the second condition in Theorem \ref{HJT}, we write HJ equation as 
\begin{equation}\label{diffGauss}
d_\vartheta(h\circ\gamma)=\beta_k\d\beta_k=\beta_k\frac{\partial\beta_k}{\partial q^i}\d q^i=0. 
\end{equation}
Finally, we obtain that the first and second equation in Theorem \ref{HJT} are equivalent to equations 
\begin{equation}\label{label1}
 F_i-\frac{F_k\beta_k\beta_i}{2c}  =  \frac{\partial \beta_i}{\partial q^k}\beta_k,\quad \quad \beta_k\frac{\partial \beta_k}{\partial q^i}=0. 
\end{equation}
It is trivial to show that the equations \eqref{label1} are equivalent if (\ref{gammacondition}) holds. 


\subsection{Contact manifolds and Dynamics}\label{Sec-Cont-Man}

The interest in contact structures roots in their applications in thermodynamics, geometric optics, geometric quantization and applications in low dimensional topology \cite{Rajeev}. Also, the theory of contact structures is linked to many other geometric backgrounds, as it is the case of symplectic geometry, Riemannian and complex geometry, analysis and dynamics \cite{Blair,Bravetti2017,deLeon2019a}. Developing a HJ on a contact manifold will be very useful for many of these applications. Let us introduce the fundamentals to develop our theory.


A $2n+1$ dimensional manifold $M$ is called contact  if it admits a one form $\eta$, known as the contact form, satisfying $\eta \wedge d\eta^n \not= 0$, see \cite{Goldstein-book,Liber87}. The contact structure is the kernel of the contact form $\eta$, and it is a nonintegrable distribution in $TM$. We denote a contact manifold by a pair $(M, \eta)$. Then, there exists a unique vector field $\mathcal R$ (called Reeb vector field) such that
\begin{equation}
\iota_{\mathcal R} \, d\eta = 0, \qquad \iota_{\mathcal R}\, \eta = 1. 
\end{equation}
The nondegenerate character of the top-form $\eta \wedge d\eta^n$ induces an isomorphism from the space of vector fields to the space of one-form sections given by
\begin{equation}\label{bar-flat}
\flat: \mathfrak{X}(M)\longrightarrow \Gamma^1(M), \qquad X \mapsto \iota_X d \eta  + \eta(X) \eta.
\end{equation}
One can establish that $\flat(\mathcal{R}) = \eta$,
so that, in this sense, we may argue that the Reeb field $\mathcal{R}$ is the \textit{dual} element of $\eta$.

Given a contact $2n+1$ dimensional manifold $(M, \eta)$, we can consider the following two distributions on $M$, that we define as horizontal distribution $ \mathfrak{H} \mathcal{D}$, and vertical distribution $ \mathfrak{V} M $ so that their Whitney sum decomposition reads
\begin{equation}
        \mathfrak{H} M = \ker \eta, \qquad          \mathfrak{V} M = \ker d \eta, \qquad TM =  \mathfrak{H}M \oplus \mathfrak{V}M,
\end{equation}
We notice that $\dim   \mathfrak{H}M = 2 n$ and $\dim  \mathfrak{V}M = 1$, and that $(d\eta)_{|_{ \mathfrak{H}M}}$ is nondegenerate and $ \mathfrak{V}M$ is generated by the Reeb vector field $\mathcal{R}$.

\textbf{Extended Cotangent Bundle.}
Consider a trivial line bundle over a manifold given by $Q\times \mathbb{R}\mapsto Q$. The first jet bundle is diffeomorphic to the extended cotangent bundle $T^*Q\times \mathbb{R}$. One can consider $T^*Q\times \mathbb{R}$ as a line bundle over the canonical symplectic manifold $(T^*Q,\Omega_{Q})$. Referring to this fibration, we pull the canonical one-form $\Theta_{Q}$ back to $T^*Q\times \mathbb{R}$ and define a contact one-form 
\begin{equation}\label{eta-Q}
\eta_{Q}:=dz-\Theta_{Q}
\end{equation}
on $T^*Q\times \mathbb{R}$ . Here, $z$ is a real coordinate on $\mathbb{R}$. There exist Darboux coordinates $(q^i,p_i,z)$ on $T^*Q\times \mathbb{R}$, where $i=1,\dots,n$. In these coordinates, the contact one-form and the Reeb vector field are computed to be
\begin{equation}\label{cont-one-form}
\eta _{Q}= d z -  p_i d q^i, \qquad \mathcal{R}=\frac{\partial}{\partial z},
\end{equation}
respectively. Notice that, in this realization, the horizontal bundle is generated by the vector fields
\begin{equation}\label{Horizontal-space}
\mathfrak{H} (T^*Q\times \mathbb{R})=span\{\xi_i,\xi^i\},\qquad \xi_i=\frac{\partial}{\partial q^i} + p_i \frac{\partial}{\partial z},~\xi^i=\frac{\partial}{\partial p_i}.
\end{equation}
We denote the inverse of the is isomorphism $\flat$ by $\sharp$. In Darboux coordinates these mappings are computed to be 
\begin{equation}
\begin{split}
\flat&: X^i\frac{\partial}{\partial q^i}+X_i\frac{\partial}{\partial p_i}+v\frac{\partial}{\partial z} \mapsto -(X_i+p_iv)dq^i+X^idp_i+vdz
\\
\sharp&: \alpha_idq^i+\alpha^idp_i+u dz \mapsto \alpha^i\frac{\partial}{\partial q^i}-(\alpha_i+p_iu)\frac{\partial}{\partial p_i}+u\frac{\partial}{\partial z},
\end{split}
\end{equation}
respectively. The nonintegrable character of the horizontal bundle can easily be observed by the simple calculation $[\xi^{i},\xi_{j}]=\delta^i_j\mathcal{R}$ 
where $\delta^i_j$ stands for the Kronecker delta. The extended cotangent bundle is a generic example for contact manifolds. That is, every contact manifold takes this form in a local level. So, for every contact manifold there exist Darboux coordinates $(q^i,p_i,z)$ admitting the contact one-form in \eqref{cont-one-form}.

\textit{Remark.} A contact manifold is an example of Jacobi manifolds. In the following part we shall examine this formulation and provide a HJ theorem for the Hamiltonian flow. On the other hand, a contact manifold can be recasted as an example of almost cosymplectic manifolds. We shall discuss this realization in Subsection \ref{Sec-Evo-Cont}.

\textbf{Contact Manifolds as Jacobi Manifolds and Hamiltonian Flow.}


Consider a contact manifold $(M,\eta)$. We define a bivector field $\Lambda$ and a vector field $E$ as follows
\begin{equation} \label{con-bivect}
    \Lambda(\alpha,\beta) = - d \eta (\sharp\alpha, \sharp\beta), \qquad
    E = - \mathcal{R},
\end{equation}
where $\sharp$ is the inverse of $\flat$ in \eqref{bar-flat} whereas ${\mathcal R}$ is the Reeb field for the contact one-form.  Direct calculations verifies that the pair $(\Lambda,E)$ satisfies the equalities in \eqref{ident-Jac}. Accordingly, we can argue that a contact manifold is a Jacobi manifold. Note that the musical isomorphism $\sharp_\Lambda$ for the bivector field $\Lambda$ in \eqref{con-bivect} is related with with the isomorphism $\sharp$ as
\begin{equation} 
\sharp_\Lambda(\alpha)=\sharp(\alpha)-\alpha({\mathcal R}){\mathcal R}. 
\end{equation}
Referring to the Hamiltonian vector field definition \eqref{Ham-v-f-Jac} given for the case of Jacobi Manifold, we define Hamiltonian vector field for a Hamiltonian function $H$ as 
\begin{equation}
X_H = \sharp_\Lambda( dH) - H {\mathcal R}.
\end{equation}
Notice that, in terms of the isomorphims $\flat $ in \eqref{bar-flat}
the Hamiltonian vector field $X_H$ is given by 
\begin{equation}
\flat (X_H) = dH - (\mathcal R (H) + H) \, \eta.
\end{equation} 
In an alternative way, the contact Hamiltonian vector field 
can be defined in terms of the contact one-form $\eta$ as follows
\begin{equation}
\iota_{X_{H}}\eta =-H,\qquad \iota_{X_{H}}d\eta =dH-\mathcal{R}(H) \eta,   \label{contact}
\end{equation}%
where $\mathcal{R}$ is the Reeb vector field.  A direct computation determines the conformal factor $g=\mathcal{R}(H)$ for a given Hamiltonian vector field
\begin{equation}\label{L-X-eta}
\mathcal{L}_{X_{H}}\eta =
d\iota_{X_{H}}\eta+\iota_{X_{H}}d\eta= -\mathcal{R}(H)\eta.
\end{equation}

In this realization, the contact Jacobi bracket of two smooth functions on $M
$ is defined by
\begin{equation}\label{cont-bracket}
\{F,H\}^c=\iota_{[X_F,X_H]}\eta, 
\end{equation}
 where $X_F$ and $X_H$ are Hamiltonian vectors fields determined through \eqref{contact}.  Here, $%
\left[ \bullet,\bullet \right]$ is the Lie bracket of vector
fields. Then, the identity $ -\left[
X_{K},X_{H}\right]=X_{\left\{ K,H\right\}^c}$  establishes the
isomorphism 
\begin{equation}
\left( \mathfrak{X}_{con}\left( M\right) ,-\left[\bullet ,\bullet
\right] \right) \longleftrightarrow \left( \mathcal{F}\left( M%
\right) ,\left\{ \bullet,\bullet \right\} ^c\right)  \label{iso1}
\end{equation}%
between the Lie algebras of real smooth functions and contact vector fields. According to \eqref{L-X-eta}, the flow of a contact Hamiltonian system preserves the contact structure, but it does not preserve neither the contact one-form nor the Hamiltonian function. Instead we obtain
\begin{equation}
{\mathcal{L}}_{X_H} \, H = - \mathcal{R}(H) H.
\end{equation}
Being a non-vanishing top-form we can consider $d\eta^n
\wedge \eta$ as a volume form on $M$.  
Hamiltonian motion does not preserve the volume form since
\begin{equation}
{\mathcal{L}}_{X_H}  \, (d\eta^n
\wedge \eta) = - (n+1)  \mathcal{R}(H) d\eta^n
\wedge \eta.
\end{equation}%
However, it is immediate to see that, for a nowhere vanishing Hamiltonian function $H$, the quantity $ {H}^{-(n+1)}   (d\eta)^n \wedge\eta$ 
is preserved along the motion.

In terms of the Darboux coordinates $(q^i,p_i,z)$, we compute the image of a one-form in $T^*M$ by $\sharp_\Lambda$ as
\begin{equation}\label{zap}
\sharp_\Lambda:\alpha_i dq^i + \alpha^i dp^i + 
u dz\mapsto 
\alpha^i \frac{\partial}{\partial q^i}-(\alpha_i + p_i u)\frac{\partial}{\partial p_i}  + \alpha^ip_i  \frac{\partial}{\partial z}. 
\end{equation}
Therefore, the Hamiltonian vector field is in the local form
\begin{equation}\label{hcont2}
X_H = \frac{\partial H}{\partial p_i} \frac{\partial}{\partial q^i} - 
\left(\frac{\partial H}{\partial q^i} + p_i \frac{\partial H}{\partial z} \right)\,  \frac{\partial}{\partial p_i} + 
\left(p_i \frac{\partial H}{\partial p_i} - H\right) \, \frac{\partial}{\partial z}
\end{equation}
Thus, an integral curve $(q^i(t), p_i(t), z(t))$ of $X_H$ satisfies the 
contact Hamilton equations 
\begin{equation}\label{hcont3} 
\frac{dq^i}{dt}   =  \frac{\partial H}{\partial p_i}, \qquad 
\frac{dp_i}{dt}   =  - \ \frac{\partial H}{\partial q^i}-  p_i \frac{\partial H}{\partial z}, \qquad 
\frac{dz}{dt}   =  p_i \frac{\partial H}{\partial p_i} - H. 
\end{equation}

\textbf{Legendrian Submanifolds.}
Let $(M,\eta)$ be a contact manifold. Recall the associated bivector field $\Lambda$. Consider a linear subbundle $ \Xi$ of the tangent bundle $TM$ (that is, a distribution on $M$). We define the contact complement of $\Xi$ as $\Xi^\perp : = \sharp_ \Lambda(\Xi^o)$, where $\Xi^o$ is the annihilator of  $\Xi$.
Let $N$ be a submanifold of $M$. We say that $N$ is:
    \begin{itemize}
        \item \emph{Isotropic} if $TN\subseteq {TN}^{\perp }$.
        \item \emph{Coisotropic} if $TN\supseteq {TN}^{\perp }$.
        \item \emph{Legendrian} if $TN= {TN}^{\perp }$.
    \end{itemize}

Consider the first order jet bundle  $\mathcal{T}^*Q$ endowed with the contact structure given in  \eqref{eta-Q}. 
Let $F$ be a real valued function on the base manifold $Q$. Its first prolongation is 
\begin{equation}\label{j1F}
\mathcal{T}^* F:Q\longrightarrow \mathcal{T}^*Q=T^*Q\times \mathbb{R},\qquad q\mapsto (dF(q),F(q)).
\end{equation}
A simple computation characterizes the local structure of Legendrian submanifolds. 
\begin{proposition}\label{legendrian_submanifold_jet}
A section $\gamma$ of $\mathcal{T}^*Q \to Q$ is a Legendrian submanifold if and only if it is the first jet prolongation of a function on $Q$. 
\end{proposition}

\subsection{HJ for Contact Hamiltonian Dynamics}\label{Sec-HJ-con-Ham}

In this section our goal is to write a geometric Hamilton-Jacobi theorem for contact Hamiltonian dynamics on the extended cotangent bundle $T^*Q\times \mathbb{R}$. This has also been studied in~\cite{Grillo2020} and in~\cite{Cannarsa2020,Jin2020} from a variational perspective in relationship to Herglotz's principle. One of the steps in this geometrization is to project the dynamics on $T^*Q\times \mathbb{R}$ to a base manifold. Due to the product structure of the total space $T^*Q\times \mathbb{R}$, we provide two Hamilton-Jacobi theories, one is obtained by projecting into $Q\times \mathbb{R}$ and the other projects onto $Q$. We shall label these formulations as Approach I and Approach II.  Both approaches and their proof are described with further detail in~\cite{de2021hamilton}.

\textbf{HJ Approach I.}
We consider the extended phase space $T^{*}Q\times \mathbb{R}$, and a Hamiltonian function
$H:T^{*}Q\times \mathbb{R} \rightarrow \mathbb{R}$ (see the diagram below). 
\[
\xymatrix{ T^{*}Q\times \mathbb{R}
\ar[dd]^{\rho} \ar[ddrr]^{z}\ar@/^2pc/[ddrr]^{H}\\
  &  & &\\
T^{*}Q &  & \mathbb{R}}
\]

Consider $\gamma$ a section of $\pi:T^{*}Q\times \mathbb{R} \rightarrow Q\times \mathbb{R}$, i.e., $\pi\circ \gamma=\text{id}_{Q\times \mathbb{R}}$. We can use $\gamma$ to project $X_H$ on $Q\times \mathbb{R}$
just defining a vector field $X_{H}^{\gamma}$ on $Q\times \mathbb{R}$ by
\begin{equation}\label{hjpar}
 X_H^{\gamma}=T_{\pi}\circ X_{H}\circ \gamma.
\end{equation}
The following diagram summarizes the above construction
\[
\xymatrix{ T^{*}Q\times \mathbb{R}
\ar[dd]^{\pi} \ar[rrr]^{X_H}&   & &T(T^{*}Q\times \mathbb{R})\ar[dd]^{T{\pi}}\\
  &  & &\\
Q\times \mathbb{R} \ar@/^2pc/[uu]^{\gamma}\ar[rrr]^{X^{\gamma}_H}&  & & T(Q\times \mathbb{R})}
\]

In this context we have the following Hamilton-Jacobi equation.
\begin{theorem}\label{hj_contact_I}
 Assume that a section $\gamma$ of the projection $T^*Q \times \mathbb{R} \longrightarrow Q \times \mathbb{R}$
is such that $\gamma(Q\times \mathbb{R})$ is a coisotropic submanifold of 
$(T^{*}Q\times \mathbb{R}, \eta_Q)$, and $\gamma_z (Q)$ is a Lagrangian submanifold of $(T^{*}Q, \Omega_Q)$, for any $z \in \mathbb{R}$. 
Then, the vector fields $X_H$ and $X_H^{\gamma}$ are $\gamma$-related if and only if 
\begin{equation}\label{hjlocal2}
\frac{\partial H}{\partial q^j} + 
\frac{\partial H}{\partial p_i} \frac{\partial \gamma_i}{\partial q^j} +
\gamma_j \left( \frac{\partial H}{\partial z} + \frac{\partial H}{\partial p_i} \frac{\partial \gamma_i}{\partial z} \right) - H \frac{\partial \gamma_j}{\partial z} = 0.
\end{equation}
Equivalently, we can rewrite this equation geometrically
\begin{equation}\label{hjglobal}
d (H \circ \gamma_z) + \gamma_o (\gamma^* \Theta_Q) - (H\circ \gamma)\cdot  \iota_{\frac{\partial}{\partial z}} d(\gamma^* \Theta_Q) = 0.
\end{equation}
where $\gamma_o =  \frac{\partial H}{\partial z} + \frac{\partial H}{\partial p_i} \frac{\partial \gamma_i}{\partial z}$.

\end{theorem}

Equations (\ref{hjlocal2}) and (\ref{hjglobal}) are indistinctly referred to as the {\it Hamilton-Jacobi equation on a contact manifold}. A section $\gamma$ fulfilling the assumptions of the theorem and the Hamilton-Jacobi equation will
be called a {\it solution} of the Hamilton-Jacobi problem for $H$.
 
Notice that if $\gamma$ is a solution of the Hamilton-Jacobi problem for $H$, then
$X_H$ is tangent to the coisotropic submanifold $\gamma(Q \times \mathbb{R})$, but it is
not necessarily tangent to the Lagrangian submanifolds $\gamma_z(Q)$, $z \in \mathbb{R}$. This occurs when $X_H(z -z_0) = 0$ for any $z_0$, that is, if and only if
$$
H \circ \gamma_{z_0} = \gamma_i \frac{\partial H}{\partial p_i}
$$
In such a case, we call $\gamma$ an {\it strong solution} of the Hamilton-Jacobi problem.

A characterization of conditions on the submanifolds $\gamma(TQ \times \mathbb{R}), \gamma_z(TQ)$ can be given as follows. Let $\sigma: Q \times \mathbb{R} \to \Lambda^k(T^* Q)$ be a $z$-dependent $k$-form on $Q$. Let $d_Q \sigma$ be the exterior derivative at fixed $z$, that is
\begin{equation}
    d_Q \sigma(q^i,z) = d \sigma_z(q^i),
\end{equation}
where $\sigma_z= \sigma(\cdot,z)$. In local coordinates, we have
\begin{equation}
    d_Q f= \frac{\partial f}{\partial  q^i} d q^i,\qquad d_Q (\alpha_i d q^i)= \frac{\partial \alpha_j}{\partial  q^i} d q^i \wedge d q^j,
\end{equation}
where $f:Q \times \mathbb{R} \to \mathbb{R}$ is a function and $\alpha = \alpha_i dq^i: Q \times \mathbb{R} \to \Gamma^1(T^* Q)$ is a $z$-dependent $1$-form.

\begin{theorem}\label{thm:coisotropic_lagrangian_section}
    Let $\gamma$ be a section of $T^*Q \times \mathbb{R}$ over $Q \times \mathbb{R}$. Then $\gamma(Q \times \mathbb{R})$ is a coisotropic submanifold and $\gamma_{z_0}(TQ)$ are Lagrangian submanifolds for all $z_0$ if and only if $d_Q \gamma = 0$ and $\mathcal{L}_{\partial/\partial z} \gamma = \sigma \gamma$ for some function $\sigma:Q\times \mathbb{R} \to \mathbb{R}$. That is, there exists locally a function $f:Q \times \mathbb{R} \to \mathbb{R}$ such that $d_Q f = \gamma$ and $d_Q \frac{\partial f}{\partial z} = \sigma d_Q f$.
\end{theorem}
The proof of this theorem can also be found in~\cite{de2021hamilton}.

\noindent
\textbf{HJ Approach II.}
Instead of considering sections of $\pi : T^*Q \times \mathbb R \longrightarrow Q \times \mathbb{R}$ as we have performed before, 
we could consider a section of the canonical projection 
$\pi : T^*Q \times \mathbb R \longrightarrow Q$, say
$\gamma :  Q \to T^*Q \times \mathbb R$. In local coordinates, we have
$(q^i) \mapsto \gamma(q^i) = (q^i, \gamma_j(q^i), \gamma_z(q^i))$. We want $\gamma$ to fulfill
\begin{equation}\label{eq:HJ}
X_H\circ\gamma=T\gamma \circ X_H^\gamma,
\end{equation}
where $X_H^\gamma=T\pi\circ X_H\circ\gamma$. Now, notice that $\tilde{\gamma} = \rho \circ \gamma$ is a 1-form on $Q$. Then, we locally have $\tilde{\gamma} = \gamma_i(q) \, dq^i$.
Next, we assume that $\gamma(Q)$ is a Legendrian submanifold of $(T^*Q \times \mathbb{R}, \eta_Q)$. This implies that $\tilde{\gamma}(Q)$ is a Lagrangian submanifold of $(T^*Q, \Omega_Q)$.
By Proposition~\ref{legendrian_submanifold_jet}, $\gamma(Q)$ is a Legendrian submanifold if and only if it is the 1-jet of a function, namely $\gamma=j^1 \gamma_z$, where we consider $\gamma_z$ as a function from $Q$ to $\mathbb R$. In other words, we have $\gamma_i = \partial \gamma_z/\partial q^i$.
If we assume that the section $\gamma$ fulfills the above condition, we can see that the Hamilon--Jacobi equation reduces to $H \circ \gamma = 0$ \cite{de2021hamilton}.
\begin{theorem}\label{hj_contact_II}
 Assume that a section $\gamma$ of the projection $T^*Q \times \mathbb{R} \longrightarrow Q$
is such that  and $\gamma (Q)$ is a Legendrian submanifold of $(T^*Q \times \mathbb{R}, \eta_Q)$. 
Then, the vector fields $X_H$ and $X_H^{\gamma}$ are $\gamma$-related if and only if $H\circ \gamma=0$. 
\end{theorem}

If so, we say that $\gamma$ is a solution of the Hamilton-Jacobi equation on a contact Hamiltonian manifold.

\textbf{An Example: The Parachute Equation.}
Consider the Hamiltonian $H$~\cite{gaset2020new}
\begin{equation}\label{eq:parachute_hamiltoniam}
    H(q,p,z) = \frac{1}{2 m}{(p+2\lambda z)}^2 + \frac{mg}{2 \lambda} (e^{2 \lambda q} -1),
\end{equation}
where $\lambda, g \in \mathbb{R}$ are constant. The extended phase space is $T^*Q \times \mathbb{R} \simeq \mathbb{R}^3$.

The Hamiltonian  field is given by
\begin{align}
    X_H = \frac{p + 2 \gamma z}{m} \frac{\partial }{\partial q} 
    &-\left(mg e^{2 \lambda q } + 2 \gamma \frac{p}{m}(p + 2 \lambda z)\right) \frac{\partial }{\partial p} + \nonumber\\ +
    &\left(  \frac{1}{2 m}{(p+2\lambda z)}^2 - \frac{mg}{2 \lambda} (e^{2 \lambda q} -1) + 2 \frac{\lambda}{m} (p + 2 \lambda z)  \right) \frac{\partial }{\partial z}.
\end{align}

\textbf{Hamiltonian vector field. Approach I}

Assume that $\gamma : Q \times \mathbb{R} \to T^* Q \times \mathbb{R}$ is a section of the canonical projection $T^* Q \times \mathbb{R} \to Q \times \mathbb{R}$, that is,
$\gamma(q,z) = (q, \gamma_p(q,z), z)$. Assume that $\gamma(Q \times \mathbb{R})$ is a coisotropic submanifold and $\gamma_{z_0}(TQ)$ are Lagrangian submanifolds for all $z_0$. Then, since $\gamma$ only has one non-trivial component, Theorem~\ref{thm:coisotropic_lagrangian_section} reduces to state that $\gamma_p = \partial{f}/ \partial q$ for some $f$ on $Q \times \mathbb{R}$.

On this situation, by Theorem~\ref{hj_contact_I}, the vector fields $X_H$ and $X_H^{\gamma}$ are $\gamma$-related if and only if
\begin{equation*}
2 \, \lambda^{3} {z}^{2} \frac{\partial f}{\partial {q}} \frac{\partial^{2} f}{\partial {q}\partial {z}} - \, g \lambda m^{2} e^{2 \, \lambda {q}} - 2 \, \lambda^{2} \left(\frac{\partial f}{\partial {q}}\right)^{2} - \, \lambda \frac{\partial f}{\partial {q}} \frac{\partial^{2} f}{\partial {q}^{2}}+ \frac{G_1}{2}\frac{\partial f}{\partial {q}} \frac{\partial^{2} f}{\partial {q}\partial {z}}-2G_2z=0
\end{equation*}
where
\begin{equation*}
G_1=g m^{2} e^{2 \, \lambda {q}} - g m^{2} + \lambda \left(\frac{\partial f}{\partial {q}}\right)^{2} - 2 \, \lambda \frac{\partial f}{\partial {q}},\quad
G_2=\lambda^2\left(2 \, \lambda \frac{\partial f}{\partial {q}} +  \frac{\partial^{2}f}{\partial {q}^{2}} - {\left( \frac{\partial f}{\partial {q}} - 1\right)} \frac{\partial f}{\partial {q}} \frac{\partial^{2} f}{\partial {q}\partial {z}}\right)
\end{equation*}
\textbf{Hamiltonian vector field. Approach II}

Assume that $\gamma : Q \to T^* Q \times \mathbb{R}$ is a section of the canonical projection $T^* Q \times \mathbb{R} \to Q$, that is, $\gamma(q) = (q,\gamma_p(q), \gamma_z(q))$.
We assume that $\gamma(Q)$ is a Legendrian submanifold of $T^*Q \times \mathbb{R}$; then, $\gamma_p(q) = \frac{\partial \gamma_z}{\partial q}$ and $X_H$ and $X^\gamma_H$ are $\gamma$-related if and only if
\begin{equation}
    H \circ \gamma = k,
\end{equation}
for a constant $k \in \mathbb{R}$. Then, the Hamilton–Jacobi equation becomes

\begin{equation}
     H(\gamma(q)) =\frac{1}{2 m}{(\gamma_p
     +2\lambda \gamma_z)}^2 + \frac{mg}{2 \lambda} (e^{2 \lambda q} -1) = k,
 \end{equation}
or, equivalently,
\begin{equation*}\label{eq:hjex1}
\frac{1}{2 m}{ \left(\frac{\partial \gamma_z}{ \partial q}
     +2\lambda \gamma_z \right)}^2 + \frac{mg}{2 \lambda} (e^{2 \lambda q} -1) = k,
\end{equation*}
This equation can be solved for $k\leq 0$. We can take the square root and obtain
\begin{equation*}
    \frac{\partial \gamma_z}{ \partial q}
     +2\lambda \gamma_z = \pm m \sqrt{\frac{g}{\lambda} e^{2 \lambda q} - 1 - k}
\end{equation*}
This equation is a linear ODE, which can be integrated using standard techniques, obtaining
\begin{equation*}
    \gamma_z(q) = \frac{e^{-2 \lambda q}}{3g\lambda} \left(3 c_1 g \lambda \pm \sqrt{\frac{g e^{2 \lambda q}}{\lambda }-k-1} \left(-e^{2 \lambda q}g + (k+1) \lambda \right) \right),
\end{equation*}
where $c_1\in \mathbb{R}$ is a constant of integration.

\section{HJ on Almost Poisson Manifolds}\label{Sec-Almost}

In this section we review the HJ theory for dynamical systems that develop their dynamics in manifolds where their compatible geometric structure does not satisfy the Jacobi identity but it fulfills the Leibniz rule. 

We start by describing a general HJ on almost Poisson manifolds, to apply it to two particular cases: nonholonomic systems and evolutionary dynamics. So, let us summarize our theory on almost Poisson manifolds to start.

\subsection{Dynamics on Almost Poisson Manifolds}\label{Sec-Almost-Pois}

   Let $P$ be a manifold equipped with a bilinear operation $\{\bullet,\bullet\}$ on the smooth functions $\mathcal{F}(P)$ and that satisfies the Leibniz identities 
\begin{equation}\label{LeibId-2}
\{F\cdot G,H\}=\{F,H\}\cdot G+F\cdot \{G,H\}, \qquad \{F,G\cdot H\}=G\cdot \{F,H\}+\{F,G\}\cdot H.
\end{equation}
for all $F$, $G$ and $H$ in $\mathcal{F}(P)$. In this case, $P$ is called a Leibniz manifold, and the bracket is called a Leibniz bracket \cite{guha2007metriplectic,Ortega-Leibniz}. 

In addition to the Leibniz identities, if further the bracket $\{\bullet,\bullet\}$ is skew-symmetric, then it is called an almost Poisson bracket \cite{Mackenzie-book,Marle-Lie-algebroid,Weinstein-Groupoid2001}. In this case the manifold $P$ is called an almost Poisson manifold. 

Consider an almost Poisson manifold $(P,\{\bullet,\bullet\})$ and the Jacobiator, 
\begin{equation}\label{Jacobiator}
\mathfrak{J}:\mathcal{F}(P)\times \mathcal{F}(P) \times \mathcal{F}(P)\longrightarrow \mathcal{F}(P), \qquad (F,G,H)\mapsto ~ \circlearrowright \{F,\{H,G\}\},
\end{equation}
where $\circlearrowright $ denotes the cyclic sum. If for an almost Poisson bracket $\{\bullet,\bullet\}$ the corresponding Jacobiator vanishes identically, then, the bracket satisfies the Jacobi identity \eqref{Jac-ident-Poiss}. An almost Poisson bracket satisfying the Jacobi identity is as Poisson bracket. In this case, $P$ turns out to be a Poisson manifold.

\textbf{Hamiltonian Vector Fields and Characteristic Distribution.} In this formalism, the definition of a Hamiltonian vector field is the same as the one in the Poisson case. We define a Hamiltonian vector field as
\begin{equation}\label{Hamvf-almost-Poiss}
X_H(F):=\{F,H\}.
\end{equation}
Accordingly, following the formalism given in Poisson case, the Hamilton equation is determined as in \eqref{HamEq}. Since we assume the skew-symmetry, the Hamiltonian function is conserved along the motion.
Further, the definition of a bivector field $\Lambda$ is exactly the same as the one in the Poisson case. So, we refer to \eqref{bivec-PoissonBra} for this definition. The rank of an almost Poisson manifold at a point $p$ is the dimension of the characteristic distribution 
\begin{equation}
\mathfrak{C}_p=\Lambda^\sharp(T^*_pP).
\end{equation}
See that the rank is not necessarily the same for every point and in this regard, it is only a generalized distribution. Further, we can state that $\mathfrak{C}$  is not integrable in general.  

\textbf{From Almost Symplectic to Almost Poisson.} Consider an almost symplectic manifold $P$ equipped with a two-form $\Omega$ which is non-degenerate but necessarily closed \cite{FassoSansonetto}. Then the identity 
\begin{equation}\label{sf-pt}
(\Lambda^\sharp)^{-1}:=-\Omega^\flat
\end{equation}
 determines a bivector field, but this time the bivector does not necessarily satisfy the Jacobi identity \eqref{Poisson-cond}. Here, $\Omega^\flat$ is the musical map associated with the symplectic two-form, as defined in \eqref{bemol}. One computes that
\begin{equation}\label{Jac-almost-sym}
\frac{1}{2}[\Lambda,\Lambda]=\Lambda^\sharp(d\Omega)
\end{equation} 
manifesting that $\Lambda$ is an almost Poisson manifold. A direct computation shows that the Poisson bracket in \eqref{sf-pt} is 
\begin{equation}
\{F,H\}:=\Omega(X_F,X_H).
\end{equation}
where $X_F$ and $X_H$ are Hamiltonian vector fields wrt. an almost symplectic structure. It is evident that if $\Omega$ is a symplectic two-form (that is additionally closed) then, the right hand side of the identity \eqref{Jac-almost-sym} vanishes. Thus, by satisfying the Jacobi identity \eqref{Poisson-cond}, one arrives at a Poisson structure. In either of the cases, the Poisson structure is non-degenerate so that there exists no non-trivial Casimir functions. 

\textbf{From Non-degenerate Two-forms to Almost Poisson.} 
It is possible to carry the discussion done in the previous paragraph to a more general geometry as follows. Consider a (regular) distribution $\mathfrak{F}$ on $P$, assume the existence of a two-form 
\begin{equation}
\Omega_\mathfrak{F}\in\Gamma(\bigwedge^2\mathfrak{F}^*)
\end{equation}
which is non-degenerate on $\mathfrak{F}$. We define now a bivector $\Lambda$  by determining the musical map $\Lambda^\sharp$ as follows. Let $\alpha$ be a one-form on $P$ and $X\in\mathfrak{X}(\mathfrak{F})$ be a vector field on $\mathfrak{F}$. Then, 
\begin{equation}\label{omega-F-Lambda}
\iota_X \Omega_\mathfrak{F}=\alpha\vert_\mathfrak{F} \quad \Leftrightarrow \quad \Lambda^\sharp(\alpha)=-X.
\end{equation}
Here, the second equation reads that $\mathfrak{F}$ is the characteristic distribution of the bivector $\Lambda$. Further, if $\mathfrak{F}$ is an integrable distribution and $\Omega_\mathfrak{F}$ is additionally closed, then, one arrives at a Poisson manifold. 
This example is the generic case of (almost) Poisson manifolds. If $\Lambda$ is an almost Poisson bivector with a regular characteristic distribution, then there exits a non-degenerate two form $\Omega_\mathfrak{F}$ satisfying \eqref{omega-F-Lambda}. If, particularly, $\Lambda$ is a Poisson bivector, then one can find an integrable distribution $\mathfrak{F}$ and a symplectic two-form $\Omega_\mathfrak{F}$ on $\mathfrak{F}$.

\begin{table}[H]{\footnotesize
  \noindent
\caption{{\small {\bf Almost Poisson Manifolds} The triple $(\tau, [\cdot,\cdot]_{D},\rho_D)$ is an almost Lie algebroid. The dual bundle $\tau^*:D\rightarrow Q$ where $D^*$ is the total space of the bundle and $\tau^*$ the projection. The almost Poisson bracket on the dual is $\{\}_{D^*}$. The dual coordinates on $D^*$ are $(q^i,p_{\alpha})$ and $\rho_{\beta}^i$ are the components of the anchor. Here $\Lambda_{D^*}$ is an almost Poisson bivector and $\{\cdot,\cdot\}_{D^*}$ is an almost Poisson bracket. Here $F^0$ is an arbitrary codistribution playing the role of reaction forces.}}
\label{table3}
\medskip
\noindent\hfill
\resizebox{\textwidth}{!}{\begin{minipage}{\textwidth}
\begin{tabular}{ l l l l }
 
\hline
 &&\\[-1.5ex]
 Structure&  Characterization &Bracket and h.v.f.&  \\[+1.0ex]
\hline
 &  & \\[-1.5ex]

 {\bf Lin Almost Pois. }&   & 
$\left\{ p_{\alpha },p_{\beta }\right\} _{D^{\ast }}=-C_{\alpha \beta
}^{\gamma }p_{\gamma }$ &   $
\Lambda _{D^{\ast }}=\rho _{\alpha }^{i}\frac{\partial }{\partial q^{i}}%
\wedge \frac{\partial }{\partial p_{\alpha }}-$  
\\[+1.0ex] $(\tau, \{\cdot,\cdot\}_{D},\rho_D)$&  $\left[ X_{\alpha },X_{\beta }\right] _{D}=C_{\alpha \beta }^{\gamma
}X_{\gamma }$&   $\left\{ q^{i},p_{\alpha }\right\}
_{D^{\ast }}=\rho _{\alpha }^{i}$  & \quad $ -\frac{1}{2}C_{\alpha \beta
}^{\gamma }p_{\gamma }\frac{\partial }{\partial p_{\alpha }}\wedge \frac{%
\partial }{\partial p_{\beta }}$ 
\\[+1.0ex] $\tau:D\rightarrow Q$ &   $\rho _{D}\left( X_{\alpha }\right) =\rho _{\alpha
}^{i}\frac{\partial }{\partial q^{i}}$ 
& $\left\{ q^{i},q^{j}\right\}
_{D^{\ast }}=0$  
& $X_{H}^{\Lambda _{D^{\ast }}}=\rho _{\beta }^{i}\frac{\partial H}{\partial
p_{\beta }}\frac{\partial }{\partial q^{i}}-$  \\[+1.0ex] 
&   &$ X_{H}^{\Lambda _{D^{\ast }}}=-\iota_{dH}\Lambda _{D^{\ast }}$   &   $\quad -\big( \rho _{\alpha }^{i}\frac{%
\partial H}{\partial q^{i}}+C_{\alpha \beta }^{\gamma }p_{\gamma }\frac{%
\partial H}{\partial p_{\beta }}\big) \frac{\partial }{\partial p_{\alpha }%
}$
   \\[+1.0ex] 
\hline\\[0.5ex]
 {\bf Nonholonomic}& $\text{dim}(T_xM)=\text{dim}(F_x),\forall x$  & $X_{H,M}|_M$ &  $\Lambda_{nh}(\alpha,\beta)=\Lambda_Q(\mathcal{P}^*(\alpha),\mathcal{P}^*(\beta))$  \\[+1.0ex]  $(T^*Q,\Omega_Q,H)$&  $TM\cap F^{\perp}=\{0\}$ & $X_c=X_H+\lambda_a Z^a$    & $\{\cdot,\cdot\}_{nh}=\Lambda_{nh}(d\cdot, d\cdot)$ \\[+1.0ex] $F^0$ & $TT^*Q|_M=TM\oplus F^{\perp}$& $(\iota_{X_c}\Omega_Q-dH)|_M\in F^0$ & 
 \\[+1.0ex] $M,\quad \{\Psi^a\}$ & & $X_c\Psi^a=0$ &   \\[+1.0ex] $\mathcal{P}:TT^*Q|_{M}\rightarrow M$ & & $X_c=\mathcal{P}(X_H)$  &   \\[+2.0ex] 
\hline\\[0.5ex]

{\bf Evolutionary Dyn.}&   & $\{f,g\}= \frac{\partial f}{\partial p_i}\frac{\partial g}{\partial q^i}-
		\frac{\partial f}{\partial q^i}\frac{\partial g}{\partial p_i}
	$ &  $\varepsilon_H=\frac{\partial H}{\partial p_i}\frac{\partial}{\partial q^i}$  \\[+1.0ex] 
	$(M,\eta)$ & $\qquad \eta\wedge\eta^n\neq 0$ & $	-\frac{\partial f}{\partial S}\big(p_i\frac{\partial g}{\partial p_i}\big)+\frac{\partial g}{\partial S}\big(p_i\frac{\partial f}{\partial p_i}\big)$ &  \quad $- \left (\frac{\partial H}{\partial q^i} + \frac{\partial H}{\partial z} p_i \right)
	\frac{\partial}{\partial p_i}$ 
	\\[+1.0ex] $\mathcal{L}_{\epsilon_H}\eta=dH-\mathcal{R}(H)\eta$ & $\qquad \eta(\varepsilon_H)=0$ & $\varepsilon_H=X_H^c+H\mathcal{R}$ & \quad $+ p_i\frac{\partial H}{\partial p_i} \frac{\partial}{\partial z}$ \\[+2.0ex] 
\hline\\[0.5ex]

    
\end{tabular}
  \end{minipage}}
\hfill}
\end{table} 

\subsection{HJ for Linear Almost Poisson Manifolds}\label{Sec-lin-almost}

To introduce our HJ on linear almost Poisson spaces we introduce some fundamental concepts: almost Lie algebroids, linear almost Poisson structures, the concept of an almost differential, see \cite{Mackenzie-book,Marle-Lie-algebroid,Weinstein-Groupoid2001}. For studies directly relating with our further sections, we refer \cite{BaMaMaPa2010,CantrijnLeonMartin1999,LeMaDi2010}. A recent publication \cite{GrilloPadron} gives a general framework in terms of fibered manifolds. 

\textbf{Almost Lie Algebroids.}
An almost Lie algebroid is a triple $(\tau,[\bullet,\bullet]_D,\rho_D)$ that consists of a vector bundle $\tau :D\rightarrow Q$ over a manifold $Q$, a skew-symmetric bilinear operation $[\bullet,\bullet]_D$ on the space $\Gamma(\tau)$ of smooth sections of $\tau$, and the anchor map $\rho_D$ from $\Gamma(\tau)$ to the space $\mathfrak{X}(Q)$ of vector fields satisfying the identity
$$ [X,fY]_D=f[X,Y]+\rho_D(X)(f)Y$$
for all real valued functions $f$ on $Q$. 
Note that, we are not necessarily assuming the satisfaction of the Jacobi identity for the bracket $[\bullet,\bullet]_D$. If the Jacobi identity is satisfied, then an almost algebroid turns out to be a Lie algebroid. 

Let us assume that the rank of the vector bundle $\tau$ is $m$, and the dimension of the base manifold $Q$ is $n$. Assume that  $\left( q^{i}\right) $ is a local coordinate system on $Q$, and $(X_{\alpha })$
is a basis of sections $\Gamma (\tau )$. Here, the Latin index $i$ runs from $1$ to $n$ whereas the Greek index $\alpha$ runs from $1$ to $m$. In these local coordinates, the almost Lie algebroid bracket $[\bullet,\bullet]_D$, and the anchor map $\rho_D$ can be written as
\begin{equation} \label{local-sa}
\left[ X_{\alpha },X_{\beta }\right] _{D}=C_{\alpha \beta }^{\gamma
}X_{\gamma }, \qquad \rho _{D}\left( X_{\alpha }\right) =\rho _{\alpha
}^{i}\frac{\partial }{\partial q^{i}}, 
\end{equation}%
where $C_{\alpha \beta }^{\gamma
}$ are the structure constant of the algebroid, whereas $\rho _{\alpha
}^{i}$ is the local representative of the anchor. 

\textbf{Linear Almost Poisson Structures.} Consider an almost Lie algebroid $(\tau,[\bullet,\bullet]_D,\rho_D)$, and the dual bundle $\tau ^{\ast }:D\rightarrow Q$ where $D^{\ast }$ is the total space of the bundle, and $%
\tau ^{\ast }$ is the projection. There exists a one-to-one correspondence between the space $\Gamma (\tau )$ of sections of
the vector bundle and the space of linear functions on the dual bundle $D^{\ast }$. This enables us to arrive at an almost Poisson bracket $\left\{ \bullet,\bullet \right\} _{D^{\ast }}$ on $D^{\ast }$. Consider the dual coordinates $(q^{i},p_{\alpha })$ on $D^{\ast }$. In this picture, the almost Poisson bracket can be written as 
\begin{equation} \label{hat}
\left\{ p_{\alpha },p_{\beta }\right\} _{D^{\ast }}=-C_{\alpha \beta
}^{\gamma }p_{\gamma }, \qquad\left\{ q^{i},p_{\alpha }\right\}
_{D^{\ast }}=\rho _{\alpha }^{i}, \qquad \left\{ q^{i},q^{j}\right\}
_{D^{\ast }}=0 
\end{equation}
where $C_{\alpha \beta }^{\gamma }$ are the structure constants of the almost Lie algebroid bracket $[\bullet,\bullet]_D$, and $\rho _{\beta }^{i}$ are the local representatives of the anchor map $\rho_D$ as determined in the equations (\ref{local-sa}). See that the bracket of two linear functions on $D^{\ast }$ is a linear function. This property labels the bracket $\left\{ \bullet,\bullet \right\} _{D^{\ast }}$ as linear almost Poisson bracket.
The associated almost Poisson $2$-vector $\Lambda _{D^{\ast }}$ can be
computed to be 
\begin{equation}
\Lambda _{D^{\ast }}=\rho _{\alpha }^{i}\frac{\partial }{\partial q^{i}}%
\wedge \frac{\partial }{\partial p_{\alpha }}-\frac{1}{2}C_{\alpha \beta
}^{\gamma }p_{\gamma }\frac{\partial }{\partial p_{\alpha }}\wedge \frac{%
\partial }{\partial p_{\beta }}.
\end{equation}
In the realm of an almost Lie algebroid $(\tau,[\bullet,\bullet]_D,\rho_D)$, the Hamiltonian dynamics defined on the dual space $(D^*,\Lambda _{D^{\ast }})$ is generated by a real valued Hamiltonian function $H$ defined by 
\begin{equation} \label{nhHamD*}
X_{H}^{\Lambda _{D^{\ast }}}=-\iota_{dH}\Lambda _{D^{\ast }}, \qquad \dot{z}=\{H,z\}_{D^*}
\end{equation}
where $\Lambda _{D^{\ast }}$ is the almost Poisson bivector and $\{\bullet,\bullet\}_{D^*}$ is the almost Poisson bracket in (\ref{hat}). In the local coordinates, the Hamiltonian vector field is computed to be
\[
X_{H}^{\Lambda _{D^{\ast }}}=\rho _{\beta }^{i}\frac{\partial H}{\partial
p_{\beta }}\frac{\partial }{\partial q^{i}}-\left( \rho _{\alpha }^{i}\frac{%
\partial H}{\partial q^{i}}+C_{\alpha \beta }^{\gamma }p_{\gamma }\frac{%
\partial H}{\partial p_{\beta }}\right) \frac{\partial }{\partial p_{\alpha }%
} 
\]%
whereas the Hamilton's equations become
\[
\dot{q}^{i}=\rho _{\beta }^{i}\frac{\partial H}{\partial p_{\beta }}, \qquad \dot{p}_{\alpha }=-\left( \rho _{\alpha }^{i}\frac{\partial H}{\partial
q^{i}}+C_{\alpha \beta }^{\gamma }p_{\gamma }\frac{\partial H}{\partial
p_{\beta }}\right) . 
\]
We shall denote a Hamiltonian system (\ref{nhHamD*}) on a linear almost Poisson space $D^*$ by the triple $(D,\{\bullet,\bullet\},H)$ where $H$ is the Hamiltonian function on $D^*$.

\textbf{Almost Differential.}
If de Rham exterior derivative of a function is zero then it is constant. This is a manifestation of integrability, that is the result of the Jacobi identity. As discussed in the previous section, in the realm of constrained systems, the brackets do not necessarily satisfy the Jacobi identity. This implies that the first cohomological space may not be equal to the field of real numbers. That is, there may exist some non-constant functions that vanish under the derivative. In this subsection, we elaborate this in the framework of almost Lie algebroids, as follows.

Assume that the triplet $(\tau,[\bullet,\bullet]_D,\rho_D)$ is a almost Lie algebroid or, equivalently, $(D^*,\{\bullet,\bullet,\}_{D^*})$ is an almost Poisson manifold. Then there exists an almost differential $d^{D}$ on the vector bundle $\tau$ defined by 
\begin{eqnarray*}
d^{D}\Psi \left( X_{0},X_{1},...,X_{k}\right)
&=&\sum\limits_{i=0}^{k}\left( -1\right) ^{i}\rho _{D}\left( X_{i}\right)
\left( \Psi \left( X_{0},X_{1},...\hat{X}_{i},...,X_{k}\right) \right) \\
&&+\sum\limits_{i<j}\left( -1\right) ^{i+j}\Psi \left( \left[ X_{i},X_{j}%
\right] _{D},X_{0},X_{1},...\hat{X}_{i},...\hat{X}_{j},...,X_{k}\right),
\end{eqnarray*}
for a $k$-form $\Psi$. 
The almost differential of a function $f$ and a one-form $\Phi$ are given by  
\begin{eqnarray*}
\left( d^{D}f\right) \left(
X\right)&=&\rho _{D}\left( X\right) (f)
\\ d^{D} \Phi \left( X,Y\right)&=&d^{D}\left(\Phi \left(
Y\right) \right) \left( X\right) -d^{D}\left( \Phi \left( X\right) \right)
\left( Y\right) -
\Phi \left( \left[ X,Y\right] _{D}\right)
\end{eqnarray*}
respectively. 

Consider the local coordinates $(q^i)$ on $Q$, a basis $(X_\alpha)$ for the sections in $\Gamma(\tau)$ and the dual sections $(X^\alpha)$ of the dual fibration $\tau^*:D^*\mapsto Q$. In these local realizations, the almost differential reads 
\[
d^{D}q^{i}=\rho _{\alpha }^{i}X^{\alpha },\qquad d^{D}X^{\gamma }=-%
\frac{1}{2}C_{\alpha \beta }^{\gamma }X^{\alpha }\wedge X^{\beta },
\]
where $\rho _{\alpha }^{i}$ are the coordinate representation of the anchor map, and $C_{\alpha \beta }^{\gamma }$ are the structure constants of the algebroid in (\ref{local-sa}). 

\begin{remark}
It is important to note that the set of linear almost Poisson structures on $D^{\ast }$, the set of skew-symmetric algebroids on $D$, and the set of almost differential are in 
one-to-one correspondence. 
\end{remark}
\begin{remark}
If an almost Lie algebroid bracket satisfies the Jacobi identity, so does the associated almost bracket. This means that, if an almost Lie algebroid is a Lie algebroid, then the almost Poisson bracket becomes a Poisson bracket. In this case the almost differential satisfies $(d^D)^2=0$. 
\end{remark}

So, we are ready to introduce our HJ on linear almost Poisson manifolds. Let $(D,\tau,Q)$ be a almost Lie algebroid, and $(D,\{\bullet,\bullet\}_{D^*},H)$ be a Hamiltonian system 
(\ref{nhHamD*}). We denote the Hamiltonian vector field by $X_{H}^{\Lambda _{D^{\ast }}}$. Consider now a section $\psi$ of the dual bundle. By commutation of the following diagram, we define a vector field $X_{H,\psi }^{\Lambda _{D^{\ast }}}$ on the base manifold $Q$ as follows. 
\begin{equation}
  \xymatrix{ D^{\ast }
\ar[dd]^{\tau^*} \ar[rrr]^{X_{H}^{\Lambda _{D^{\ast }}}}&   & &TD^{\ast }\ar[dd]^{T\tau^*}\\
  &  & &\\
 Q\ar@/^2pc/[uu]^{\psi}\ar[rrr]^{X_{H,\psi }^{\Lambda _{D^{\ast }}}}&  & & TQ }
\end{equation}
Here, $T\tau^*$ is the tangent lift of the dual projection $\tau^*$. 
It is interesting to note that the image space of the reduced vector field $X_{H,\Psi }^{\Lambda _{D^{\ast }}}$ at every point $q$ in $Q$ is an element of the image space $\rho_{D}\left( D_{q}\right)$ of the anchor map at that point. For the details and proof of this assertion, we cite \cite{LeMaDi2010,deLedeDiVa14}.
\begin{theorem}\label{HJ-skew}
Let $(D,\{$\textperiodcentered $,$\textperiodcentered $\}_{D^{\ast }},H)$ be
a Hamiltonian system. Assume that $\psi$ be a section of the dual bundle $(D^*,\tau^*,Q)$ such that $%
d^{D}\Psi =0$. The following conditions are equivalent.
\begin{enumerate}
\item If $c\left( t\right) $ is an integral curve of $X_{H,\psi
}^{\Lambda _{D^{\ast }}}$ then $\psi \circ c(t)$ is a solution of $%
X_{H}^{\Lambda _{D^{\ast }}}$.

\item $\psi $ satisfies the Hamilton-Jacobi equation%
\begin{equation}
d^{D}(h\circ \psi )=0.
\end{equation}
\end{enumerate}
\end{theorem}

In the local 
picture, the previous theorem takes the following particular form. Let $\psi (q) = ( q^{i},\psi_{\gamma }( q ))$ then the derivative of $\psi$ with $d^D$ is zero if
\begin{equation}
d^{D}\psi =0\Longleftrightarrow C_{\delta \beta }^{\gamma }\psi _{\gamma
}=\rho _{\delta }^{i}\frac{\partial \psi_{\beta }}{\partial q^{i}}-\rho
_{\beta }^{i}\frac{\partial \psi _{\delta }}{\partial q^{i}},
\end{equation}
whereas the vanishing of the derivative of the composition $h\circ \psi$ can be written as
\begin{equation}
d^{D}\left( H\circ \psi \right) =0\Longleftrightarrow \rho _{\gamma
}^{i}\left( \frac{\partial H}{\partial q^{i}}+\frac{\partial \psi _{\beta }%
}{\partial q^{i}}\frac{\partial H}{\partial p_{\beta }}\right) =0.
\end{equation}

It is immediate to observe that when $D$ is the tangent bundle, the anchor map is the identity mapping and the Lie bracket on the sections of the almost Lie algebroid is the Jacobi Lie bracket of vector fields on $Q$. Then, the theorem reduces to the one in (\ref{HJT}). Note that, in this case, the derivative $d^D$ turns out to be the de Rham exterior derivative. In this sense, this theorem is a generalization of the geometric Hamilton-Jacobi theorem (\ref{HJT}).

\textbf{Totally Nonholonomic Case.} As we shall discuss in the following section, a submanifold of the tangent manifold is called nonholonomic if it is not integrable. Since the Jacobi identity is not assumed in the present discussion, we can claim that Theorem \ref{HJ-skew} is proper for so called nonholonomic dynamics. Let us consider the following particular realization of an almost Lie algebroid valid for many physical theories. Consider a generalized distribution $\tilde{D}$ on a maniold $Q$ whose
characteristic space is $\tilde{D}_{q}=\rho _{D}\left( D_{q}\right) $.
Denote by $Lie^{\infty }\left( \tilde{D}\right) $ the smallest Lie subalgebra of $%
\mathfrak{X}(Q)$ containing $\tilde{D}$. Consider 
\[
Lie_{q}^{\infty }\left( \tilde{D}\right) =\left\{ \tilde{X}\left( q\right)
\in T_{q}Q:\tilde{X}\in Lie^{\infty }\left( \tilde{D}\right) \right\} . 
\]%
then 
\[
q\rightarrow Lie_{q}^{\infty }\left( \tilde{D}\right) \subset T_{q}Q 
\]%
defines a generalized foliation on $Q$. 
The skew-symmetric algebroid $(D,[$ $,$ $]_{D},\rho _{D})$ over Q is said to
be completely nonholonomic if 
\[
Lie_{q}^{\infty }\left( \rho _{D}\left( D\right) \right) =Lie_{q}^{\infty
}\left( \tilde{D}\right) =T_{q}Q 
\]%
for all $q\in Q$. If the skew-symmetric algebroid $(D,[$ $,$ $]_{D},\rho _{D})$ over $Q$ is
completely nonholonomic and $Q$ is connected, then $H^{0}\left( d^{D}\right) =%
\mathbb{R}
$. However, the condition $H^{0}\left( d^{D}\right) =%
\mathbb{R}
$ does not imply, in general, that the almost Lie algebroid $(D,[$ $,$ $%
]_{D},\rho _{D})$ is completely nonholonomic. Note that, if $Q$ is connected and $D$ is a transitive almost Lie algebroid ($\rho
_{D}\left( D_{q}\right) =T_{q}Q$), then $H^{0}\left( d^{D}\right)=\mathbb{R}$ is so. If $H^{0}\left( d^{D}\right) =\mathbb{R}$ or that the almost Lie  algebroid is completely nonholonomic and $Q$ is
connected, then the second condition in Theorem \ref{HJ-skew} turns out to be $h\circ \psi =cst$. Let us record this case in the following theorem.
 
\begin{theorem}\label{Thm-HJ-skew}
Let $(D,\{$\textperiodcentered $,$\textperiodcentered $\}_{D^{\ast }},H)$ be
a Hamiltonian system. Assume that $H^{0}\left( d^{D}\right) =\mathbb{R}$ or that the almost Lie  algebroid is completely nonholonomic and $Q$ is
connected. Consider also that $\psi$ is a section of the dual bundle $(D^*,\tau^*,Q)$ such that $%
d^{D}\Psi =0$. The following conditions are equivalent.
\begin{enumerate}
\item If $c\left( t\right) $ is an integral curve of $X_{H,\psi
}^{\Lambda _{D^{\ast }}}$ then $\psi \circ c(t)$ is a solution of $%
X_{H}^{\Lambda _{D^{\ast }}}$.

\item $\psi $ satisfies the Hamilton-Jacobi equation%
\begin{equation}
H\circ \psi =cst.
\end{equation}
\end{enumerate}
\end{theorem}

\subsection{Nonholonomic Hamiltonian Dynamics}

In the literature, the nonholonomic dynamics first examined in the realm  of Lagrangian dynamics. In this setting, the nonholonomic constraints are defined to be a nonintegrable submanifold of velocity phase space. In general, it is assumed that the projection of this submanifold to the configuration space is assumed to be surjective. This physically corresponds to the idiom that "all positions are permitted but not all velocities". In this survey, we approach the nonholonomic formulation from the point of view of Hamiltonian dynamics. 

Hamiltonian dynamics under nonholonomic constraints has been widely discussed in various studies \cite{BatesSniatycki93,KoonMarsden97,Marle1998,SchaftMaschke94}. The lack of integrability of the constraints manifests the lack of fulfillment of the Jacobi identity. So, one needs to work on an almost Poisson manifold \cite{Balseiro12,BlochMarsdenZenkov2005,CantrijnLeonMartin1999,Sniatycki2001}. 
To develop a HJ theory for nonholonomic systems we need to introduce first some fundamentals of constrained dynamics. Here we start then by introducing constrained Hamiltonian dynamics, a nonholonomic bracket, etc.

{\bf{Constrained Hamiltonian Dynamics.}}
Consider a Hamiltonian system $(T^*Q,\Omega_Q,H)$. 
Let $M$ be a submanifold of $T^*Q$ of codimension $k$ defining constraints. Then, $M$ may be described in terms of a set $\{\Psi^a \}$ of independent constraint functions in the following way
\begin{equation}
M=\left\{z\in T^*Q: \Psi^a(z)=0\right\},
\end{equation}
where $a$ runs from $1$ to $k$. 
The equation of motion associated to the constrained Hamiltonian system is defined to be 
\begin{equation} \label{nHELag}
\left. \left(\iota_{X_c}\Omega_{Q} -dH\right) \right\vert _{M}\in F^{o}, \qquad
\left. X_{H,M} \right\vert _{M} \in TM,
\end{equation}
where $F^{o}$ is just an arbitrary codistribution on $T^{*}Q$ playing the role of reaction forces. One can define a distribution $F$ on $T^{*}Q$ such that $F^{o}$ is the annihilator of $F$. 
Here, the latter condition in \eqref{nHELag} suggests that the image space of the constraint vector field $X_{H,M}$ is in the tangent bundle $TM$ of the constraint manifold.

To guarantee the existence and uniqueness of a vector field $X_{H,M}$ satisfying  \eqref{nHELag} we need to impose two conditions:
\begin{itemize}
\item[(i)] Admissibility condition: 
\begin{equation}\label{admissibility condition}
dim \left( T_{x} M \right) = dim \left( F_{x} \right), \ \forall x \in M
\end{equation}
\item[(ii)] Compatibility condition: 
\begin{equation}\label{compatibility condition}
TM \cap F^{\perp} = \{0 \}.
\end{equation}
\end{itemize}
Then, the tangent space at each point $z$ in the constrained submanifold $M$ can be written as a direct sum of $T_{z}M$ and $F^{\bot}_{z}$, that is 
\begin{equation}\label{projections}
\left. TT^*Q\right\vert _{M}=TM\oplus F^{\bot}. 
 \end{equation}
and we can define two complementary projectors as well
\begin{equation}
\mathcal{P}:TT^*Q \vert_M \mapsto TM, \qquad Q:TT^*Q\vert_M \mapsto F^{\bot}.
 \end{equation}

\textbf{Constrained Dynamics as a Projection.} Let us try to exhibit some of the features given in more concrete terms by using the one-forms and vector fields.
 Assume that the codistribution $F^0$ is generated by the set of one-forms $(\sigma^{a})$ where $a$ runs from $1$ to the dimension of $F^0$, say $k$. By using the musical isomorphism $\Omega^{\sharp}_{Q}$, we compute the symplectic gradients of these generators and obtain a basis $\{Z^a\}$ for the symplectic orthogonal $F^\bot$ of the distribution $F$. If the admissibility condition is assumed, then the submanifold $M$ will be defined in terms of the set of constraints
$\{\Psi^a\}$ where the index $a$ runs from $1$ to the codimension of $M$.  If we denote $X_{H}$ as the unconstrained Hamiltonian vector field, a solution of the constrained system takes the form 
\begin{equation} \label{-soln-nHHe}
X_c=X_{H}+\lambda_{a}Z^a, 
\end{equation}
where $\{\beta_{a}\}$ are Lagrange multipliers which can
be determined as follows. We have the tangency condition that $X\left(\Psi^a\right)$ vanishes identically. So that, we have a linear system of
equations%
\[
X_{H}\left(\Psi^a\right) +\lambda_{b}Z^{b}\left(\Psi^a\right)=0. 
\]%
The Lagrange multipliers can be determined up to some extent depending on
the degeneracy level of the matrix 
\begin{equation} \label{Cab}
[C^{ab}]=[Z^a(\Psi^b)].
\end{equation}
If the compatibility condition holds, then the rank of the matrix $[C^{ab}]$ is the maximum. In this case, one can compute the Lagrange multipliers $\lambda_{a}$ uniquely and write a solution as
\begin{equation} \label{Proj-exp}
X_{c}=\mathcal{P}(X_H)=X_H-C_{ab}X_H(\Psi^{b})Z^a
\end{equation}
where $C_{ab}$ are the components of the inverse of the matrix $[C^{ab}]$ in (\ref{Cab}). 
 
\textbf{Nonholonomic Bracket.}
Consider a constrained system \eqref{nHELag} holding the compatibility condition \eqref{compatibility condition}. As previously stated, in this case, we arrive at a well-defined projection (\ref{projections}) from the total space $TT^*Q$ to the tangent bundle $TM$ of the constrained submanifold. Using this projection, we push the canonical Poisson bracket on $T^*Q$ to a bracket on $M$ \cite{CantrijnLeonMartin1999}. To see this bracket, consider two functions $F_1$ and $F_2$ in $C^{\infty}(M)$, then 
\begin{equation} \label{nh-bracket}
\{F_1,F_2\}_{nh}=\Omega_{Q}(\mathcal{P}(X_{F_1}),\mathcal{P}(X_{F_2})).
\end{equation}
This bracket is called a Dirac bracket \cite{IbLeMa}. Nowadays, we call it the nonholonomic bracket. 
It is easy to see that the bracket $\{\bullet,\bullet\}_{nh}$ is skew-symmetric and verifies the Leibniz rule. However, in general, the Jacobi identity is not satisfied. So  $\{\bullet,\bullet\}_{nh}$ is an almost Poisson structure. 

Recall that, by assuming the admissibility condition (\ref{admissibility condition}), we can write the projection $\mathcal{P}$ in terms of the symplectic gradients $(Z_j)$ of the generators of the codistribution, and constraint functions $(\Psi_a)$ as in (\ref{Proj-exp}). In this case, the nonholonomic bracket (\ref{nh-bracket}) can be computed to be  
\begin{eqnarray} \label{nhbraclocal}
\{F_1,F_2\}_{nh}=\{F_1,F_2\}&+&C_{ab}Z^a(F_2)\{F_1,\Psi^b \}-C_{ab}Z^a(F_1)\{F_2,\Psi^b\}\notag \\&+&C_{ab}C_{cd}
\{\Psi^b,\Psi^d\}Z^a(F_1)Z^d(F_2),
\end{eqnarray}
where the bracket on the right hand side $\{\bullet,\bullet\}$ is the canonical Poisson bracket on $T^*Q$, and $C_{ab}$ are the components of the inverse matrix of $[C^{ab}]=[Z^a(\Psi^b)]$. Note that the constrained functions $\Psi^a=\Psi^a(q,p)$ are Casimir functions. Now, the dynamics can be written in terms of the nonholonomic bracket as 
   \begin{eqnarray}
   \dot{z}=\{H,z\}_{nh}.
   \end{eqnarray}
Here, $z=z(t)$ is a curve in $M$.
As a result, we deduce that if a unique solution of the constrained dynamics exists, then there exists an almost Poisson structure (\ref{nhbraclocal}) on the constrained manifold $M$. In this case, the constrained dynamics is generated by the Hamiltonian function $H$ as well.
Consider the canonical Poisson $\Lambda_Q$ induced by the canonical symplectic structure.
For any two sections $\alpha$ and $\beta$ of $T^{*}_M T^{*}Q$ we put
$$\Lambda_{nh}(\alpha, \beta) = \Lambda_Q\left(\mathcal{P}^{*}(\alpha), \mathcal{P}^{*}(\beta)\right).$$
\noindent
with $\mathcal{P}$ the projector onto $TM$ defined previously. It is easy to see that the operator $\Lambda_{nh}$ can be extended to a genuine two-contravariant skew-symmetric tensor field on an open neighbourhood of $M$ in $T^{*}Q$. Indeed, it suffices to observe that the projector $\mathcal{P}$ can always be extended to a type $(1, 1)$- tensor field on a neighbourhood of $M$. This $\Lambda_{nh}$ is obviously related to the nonholonomic bracket as: $\{\bullet, \bullet \}_{nh}= \Lambda_{nh}\left(d\bullet, d\bullet\right)$.

\textbf{Almost Lie Algebroid Realization.}
Start by defining a set of linearly independent sections $(X_\alpha,Y_a)$, where $a$ runs from $1$ to $k$, and $\alpha$ runs from $1$ to $m$ with $m+k$ equal to $n$ of the tangent bundle $TQ$ as follows, for arbitrary $ a=1,...k$ and $\alpha=1,...,m$, we have that
\begin{equation}
\langle Y^a, X_\alpha \rangle =0, \qquad \langle Y^a, Y_b \rangle =\delta^{a}_b,
\end{equation}
where $\delta^{b}_a$ is the Kronecker delta. In a local chart, we have that 
\begin{equation} \label{genD}
X_\alpha= X_\alpha^i\frac{\partial}{\partial q^i}, \qquad Y_a= Y_a^i\frac{\partial}{\partial q^i},
\end{equation}
where the coefficients $X_\alpha^i$ and $Y_a^i$ are real valued functions on $Q$. Note that, at every $q$ in $Q$, the tangent space $T_qQ$ is spanned by the set of vectors $\{X_\alpha(q),Y_a(q)\}$. By considering the subset $\{X_\alpha(q)\}$, one arrives at a subspace of $T_qQ$, denoted by $\bar{D}_q$. The union results in the total space $\bar{D}$ of a vector bundle $\tau:\bar{D}\rightarrow Q$ which is the geometric framework of the skew algebroid.

The vector fields $(X_\alpha,Y_a)$ can be considered as linear functions on $T^*Q$ as follows 
 \begin{equation}
 X_\alpha(z)= \langle X_\alpha(\pi_Q(z)), z \rangle,
 \end{equation}
where the pairing $\langle \dot, \dot \rangle$ is the linear algebraic duality between $T_qQ$ and $T^*_qQ$, with $q=\pi_Q(z)$. Using this identification, we define a new set of coordinates $(q^i,p_\alpha,p_a)$ on the cotangent bundle $T^*Q$. Explicitly, the dual coordinates are defined as
\begin{equation}
p_\alpha= X_\alpha (z), \qquad p_a=Y_a(z).
\end{equation}
Considering the Darboux coordinates $z=(q^i,p_j)$  we have that 
\begin{equation}
p_i=X^\alpha_i p_\alpha+Y^a_{i}p_a, \quad\qquad p_\alpha= X_\alpha^i p_i,\quad p_a= Y_a^i p_i,
\end{equation}
where $X^\alpha(X_\beta)=\delta^\alpha_\beta$, $X^\alpha(Y_a)=0$ for all $\alpha$ and $a$. In other words, the set $(X^\alpha,Y^a)$ forms a linear independent basis for the sections of the cotangent bundle fibration.  

In the coordinate system $(q^i,p_\alpha,p_a)$, the constrained submanifold $M$ presented in (\ref{LT}) turns out to be 
\begin{equation} \label{LT-}
M=\left \{(q^i,p_\alpha,p_a)\in T^*Q: \frac{\partial H}{\partial p_a}(q^i,p_\alpha,p_a)=0, \quad a=1,..,k \right \},
\end{equation}
We have assumed that the Lagrangian function $L$ is non-degenerate. This implies the non-degeneracy of the Hamiltonian function $H$. 
In light of this fact, a straightforward calculation proves that the rank of the matrix $(\frac{\partial^2 H}{\partial p_a\partial p_b})$ is full. 
So that, in a local open neighborhood, one may solve $ p_a$ in terms of $(q^i,p_\alpha)$. 
This observation result with a local coordinate system $(q^i,p_\alpha)$ for the constrained submanifold $M$, and local sections $\{X^\alpha\}$ from the base manifold $Q$ to the constrained manifold $M$. 
Using these sections, we define a vector space $\bar{D}^*_q$ at each point $q$ in $Q$, and we can consider $M$ as the total space of the vector bundle with projection $\pi_Q \rvert_M$ and the base manifold $Q$. Here, the fibration $\pi_Q \rvert_M$ is the restriction of the cotangent bundle projection $\pi_Q$ to $M$. Note that this bundle is dual to the bundle $\bar{D}\rightarrow Q$ generated by the sections $\{X_\alpha\}$ presented in (\ref{genD}). This leads us to denote $M$ by $\bar{D}^*$ as well. The canonical Poisson bracket on the cotangent bundle $T^*Q$ reduces to an almost Poisson bracket on $M=\bar{D}^*$. It is computed to be 
\begin{eqnarray}
\{q^i,q^j\}_{nh}&=&0, 
\qquad 
\{q^i,p_\alpha\}_{nh}=X^i_\alpha, \qquad 
\{p_\alpha,p_\beta\}_{nh} =  \left (\frac{\partial X^i_\alpha}{\partial q^j} X^j_\beta - 
\frac{\partial X^i_\beta}{\partial q^j} X^j_\alpha 
  \right)p_i \bigg\rvert_M.
\end{eqnarray}
Here, the first two brackets depend on the induced coordinates $(q^i,p_\alpha)$ on $M$, whereas the last one is defined in terms of $(p_i)$. For the latter one, one should substitute the constraints $p_a=p_a(q^i,p_\alpha)$ to arrive at a final representation. The constrained dynamics in (\ref{nHELag}) is given by
\begin{equation}
\dot{q}^i=\{q^i,H_M\}_{nh}, \qquad \dot{p}_\alpha=\{p_\alpha,H_M\}_{nh}
\end{equation}
where $H_M$ is the restriction of $H$ to the submanifold $M$.

\subsection{HJ for Nonholonomic Hamiltonian Dynamics}\label{Sec-HJ-nonholo}

Consider a submanifold $N \subseteq TQ$ of codimension $k$ in $TQ$ locally described in terms of independent constraint functions $\{\psi^a \}_{a=1, \dots , k}$ in the following way
\begin{equation} \label{Npre}
N=\left\{(q^i,\dot{q}^j)\in TQ: \psi^a(q,\dot{q})=0\right\}.
\end{equation}
We assume that $\tau_{Q}  ( N  ) = Q$ or, equivalently, the constraints are purely kinematic.
A particularly important example of  constraint manifold $N$ is a linear subbundle of the tangent bundle $TQ$. In this case, $N$ is locally generated by linear constraint functions in the velocities given by
\begin{equation} \label{constLin}
\psi^a(q,\dot{q})=\psi_i^a(q)\dot{q}^i.
\end{equation}
So, we can write $N$ in this case
\begin{equation} \label{N2}
N =\left\{(q^i,\dot{q}^j)\in TQ: \psi_i^a(q)\dot{q}^i=0\right\}.
\end{equation}
Accordingly, we can define a family of one-forms
\begin{equation} \label{barpsi}
\overline{\psi}^a(q)=\psi_i^a(q) dq^i \in \Gamma^1(Q)
\end{equation}
It is easy to realize that the codistribution $N ^{o}$ on $Q$ given by the annihilator of $N $ is (locally) generated by the $1-$forms $\overline{\psi}^{a}$. Then, by pulling the one-forms $\overline{\psi}^a$ back to $TQ$ with the help of the tangent bundle projection $\tau_Q$, we arrive at one-forms
\begin{equation} \label{barpsilift}
\tau_Q^*\overline{\psi}^a(q)=\psi_i^a(q) dq^i = S^{*}\left( d \psi^a \right)\in \Gamma^1(TQ).
\end{equation} 
Here, $S$ is the almost tangent structure on $TQ$.
Hence,
$$\tau_{Q}^{*}\left( N ^{o}\right) = S^{*}\left( TN ^{o}\right).$$
As a particular case, let us now apply the 
(inverse) Legendre transformation and arrive at
\begin{equation} \label{LT}
M=\left \{(q,p)\in T^*Q: \Psi^a(q,p)=\psi^a_{i}(q)\frac{\partial H}{\partial p_i}(q,p) =0 \right \}
\end{equation}
of the cotangent bundle $T^*Q$. With these definitions, $M$ and 
\begin{equation} \label{F^{o}Lag}
F^{o}=\left\langle\sigma^a:=\pi_Q^*\overline{\psi}^a= \psi^a_{i}(q) dq^i \right\rangle.
\end{equation}
satisfy the admissibility condition. In Darboux coordinates, the constraint Hamiltonian system (\ref{nHELag}) is computed to be
\begin{equation} \label{nhHamEqcoor}
\dot{q}^i=\frac{\partial H}{\partial p_i}, \qquad \dot{p}_i=-\frac{\partial H}{\partial q^i}+\beta_a\psi^a_i(q), \qquad  \psi^a_i(q)\frac{\partial H}{\partial p_i}=0.
\end{equation}
Here, the Hamiltonian function has the form (\ref{LT}). 
It is possible to determine the Lagrange multipliers $\beta_a$ in the equations of motion by simply taking the derivative of the constraint with respect to time. By taking the isomorphic image of the space $F^{o}$ in Eq. (\ref{F^{o}Lag}), we compute the symplectic orthogonal distribution 
\begin{equation}
F^\bot=\left\langle Z^a=\psi^a_{i}(q)\partial / \partial p_i \right \rangle
\end{equation}
over the cotangent bundle $T^*Q$.

Consider the vector subbundle $N$ given in \eqref{N2}. As discussed previously, see Eq. (\ref{barpsi}), the constraint functions $\psi ^a_i(q)$ lead to the determination of the set of differential one-forms $\overline{\psi^a}$ on $Q$. Note that, the image space of these one-forms determine a subbundle $N ^o$ of the cotangent bundle $T^*Q$ annihilating $N $, and  
define an ideal 
\begin{equation}
{\mathcal I}(N ^o)=\left\{\beta_a\wedge \overline{\psi ^a}: \beta_a\in \Lambda^s(Q)  \right\}
\end{equation}
 of the exterior algebra $\Lambda(Q)$. In this setting, a vector field $X$ on $Q$ is called a characteristic vector field of the ideal satisfying $\iota_X(\overline{\psi^a})=0$ for all $\overline{\psi^a}$. A characteristic vector field $X$ of ${\mathcal I}(N ^o)$ preserves the ideal, that is, $\iota_X{\mathcal I}(N ^0)\subset {\mathcal I}(N ^o)$. Notice that a vector field $X$ taking values in the constraint subbundle $N $ is a characteristic vector field of the ideal. 
 
 As above, for a fixed $1-$form $\sigma$ on $Q$ such that $\sigma \left( Q \right) \subseteq M$ we will define a vector field $X_{H,M}^{\sigma}$ on $Q$ by satisfying the following diagram,
\begin{equation}\label{Xg22}
  \xymatrix{ M\subset T^* Q
\ar[dd]^{\pi_{Q}} \ar[rrr]^{X_{H,M}}&   & &TT^* Q\ar[dd]^{T\pi_{Q}}\\
  &  & &\\
 Q\ar@/^2pc/[uu]^{\sigma}\ar[rrr]^{X_{H,M}^{\sigma}}&  & & TQ }
\end{equation}

 \begin{theorem}\label{nhhj1}
Let $\sigma$ be a 1-form on $Q$ such that $\sigma(Q) \subset
{M}$ and $d\sigma\in {\mathcal I} (N ^0)$. Then the following
conditions are equivalent:

\begin{itemize}
\item[(i)] $X_{H,M}^{\sigma}$ and $X_{H,M}$ are $\sigma$-related

\item[(ii)] $d(H \circ \sigma) \in N ^0$.
\end{itemize}
\end{theorem}

A distribution $\mathcal{D} \leq TQ$ is said to be completely nonholonomic (or bracket-generating) if
$\mathcal{D}$ along with all of its iterated Lie brackets $\left[\mathcal{D}, \mathcal{D}\right], \left[ \mathcal{D}, \left[ \mathcal{D} , \mathcal{D} \right] \right], \hdots$ spans the tangent bundle $TQ$. 
An important result related with this kind of distribution is the following
\begin{theorem}[Chow-Rashevskii's theorem]\label{LordChowTheorem}
Let $Q$ be a connected differentiable manifold. If a distribution $\mathcal{D} \leq TQ$ is completely nonholonomic, then any two points on $Q$ can be joined by a horizontal path $\gamma$, i.e., the derivative of $\gamma$ is tangent to $\mathcal{D}$.

\end{theorem}
A detailed proof of this theorem for regular distributions can be found in \cite{Montg}. For singular distributions see \cite{Harms}. A particularly important consequence of this theorem is given by the following result.\\
\begin{proposition}\label{Importantconsequuence2134}
Let $Q$ be a connected differentiable manifold and $\mathcal{D} \leq TQ$ be a completely nonholonomic distribution. Then there is no non-zero exact one-form in the annihilator $\mathcal{D}^{o} \leq T^{*}Q$.
\end{proposition}
It is important to remark that this proposition can be analogously proved for the case in which $\mathcal{D}$ is a singular distribution by taking into account the Chow-Rashevskii's theorem for singular distributions.

Let us now present this same theorem in its Hamiltonian counterpart \cite{NonholIglesiasLeonDiego}, notice that the constraint submanifold is defined to be $M=\mathbb{F}L(N_\ell)$ as in Eq. (\ref{LT}).

\begin{theorem}\label{nhhj1corollary4again}
Assume that the distribution defined by the constraint vector bundle $N $ is completely nonholonomic. Let $\sigma$ be a 1-form on $Q$ such that $\sigma(Q) \subset
{M}$ and $d\sigma\in {\mathcal I} (N ^0)$. Then, the following
conditions are equivalent:

\begin{itemize}
\item[(i)] $\Gamma_{H,M}^{\sigma}$ and $\Gamma_{H,M}$ are $\sigma$-related

\item[(ii)] $d(H \circ \sigma) =0$.
\end{itemize}
\end{theorem}

{\bf Example: Nonholonomic free particle}

As an example of application of a nonholonomic HJ theory let us consider a nonholonomic free particle in a three-dimensional space  \cite{BatesSniatycki93,CaGrMaMaMuMiRa,deLedeDiVa14,EsJiLeSa19,GracMar}. We have $Q=\mathbb R^3$ and the Lagrangian of this system reads
$$ L= \frac{1}{2}m(\dot x^2+\dot y^2+\dot z^2)-V(x,y,z)  $$
with constraints given by the equation
\begin{equation}\label{eqn:distbn}
\dot z-y\dot x=0.  
\end{equation}
From (\ref{eqn:distbn}) it is straightforward to see that the distribution describing the constraints is the subspace
$$ N=span\bigg\{  \frac{\partial }{\partial x}+y\frac{\partial }{\partial z}    , \frac{\partial }{\partial y} \bigg\}  .$$
The Legendre transform for this system is a map
$$ \lambda:TQ\to T^*Q, \qquad \lambda(x,y,z,\dot x,\dot y,\dot z)= (x,y,z,m\dot x,m\dot y,m\dot z), $$
which means that in the Hamiltonian picture the system is described by a Hamiltonian
\begin{equation}\label{nonholham}
 H= \frac{1}{2m}(p_x^2+p_y^2+p_z^2)+V(x,y,z). 
\end{equation}
Using Legendre map we also obtain a constraint submanifold $M:=\lambda(N)$ in $T^*Q$, which is described by the equation $ p_z-yp_x=0$ and therefore it is spanned by the one form 
$$ M=span\bigg\{   dx+ydz    \bigg\}  .$$
Using the above results we can find the decomposition (\ref{projections}), which in this case reads
$$  TT^*Q|_M=TM\oplus F^\perp, $$ 
$$F^\perp=span\bigg\{ \frac{\partial}{\partial p_z} - y\frac{\partial}{\partial p_x}  \bigg\}. $$
From the above one can derive that the projection on $TM$ is a map
$$ P = id_{TT^*Q_{|M}}-\frac{1}{1+y^2}\bigg( \frac{\partial}{\partial p_z}-y\frac{\partial}{\partial p_x}  \bigg)\otimes (dp_z-ydp_x-p_xdy).    $$ 
One can use the above projection to write down the nonholonomic Hamiltonian vector field of (\ref{nonholham}), which by definition reads
$$ X_{H,M}=P(X_H)=\frac{p_x}{m}\frac{\partial}{\partial x}+\frac{p_y}{m}\frac{\partial}{\partial y}+\frac{p_z}{m}\frac{\partial}{\partial z} -\frac{1}{y^2+1}\Big( \frac{\partial V}{\partial x}+y\frac{\partial V}{\partial z}+\frac{yp_xp_y}{m^2}  \Big)\frac{\partial}{\partial p_x}   $$
$$ \qquad -\frac{\partial V}{\partial y}\frac{\partial }{\partial p_y}  -\frac{1}{y^2+1}\Big(   y\frac{\partial V}{\partial x}+y^2\frac{\partial V}{\partial z}-\frac{p_xp_y}{m^2}  \Big)\frac{\partial}{\partial p_z}.    $$
Let us move to the HJ theorem. We start with a one-form $\sigma=\sigma_1dx+\sigma_2dy+\sigma_3dz$. To enforce the condition $d\sigma_{|_{N\times N}}=0$ we have to calculate
$$ d\sigma\Big(   \frac{\partial }{\partial x}+y\frac{\partial }{\partial z}    , \frac{\partial }{\partial y}  \Big)=0, $$
which in coordinates leads to the equation
$$ (1+y^2)\frac{\partial\sigma_1}{\partial y}+y\sigma_1-\frac{\partial\sigma_2}{\partial x}-y\frac{\partial\sigma_2}{\partial z}=0.  $$
Notice that $\{ X_1,X_2,[X_1,X_2] \}=T^*Q$, which means that the distribution $N$ is completely nonholonomic. Therefore, the equation $d(H\circ\sigma)\in N^o$ is equivalent to the equation  $H\circ\sigma=E$ and the nonholonomic HJ equation reads
$$ \frac{1}{2m}\Big( (1+y^2)\sigma_1^2+\sigma_2^2 \Big)+V(x,y,z)  = E, $$
where $E$ is a constant. For a discussion of the solution of the above equation for $m=1$ and $V=0$ see \cite{deLedeDiVa14}.

\subsection{HJ for Evolution Dynamics on Contact Manifolds}\label{Sec-Evo-Cont}

Another vector field can be defined from a Hamiltonian function $H$ on a contact manifold $(M,\mathcal{\eta})$, this is the so-called the \emph{evolution vector field} of $H$~\cite{simoes2020contact}.
The use of the evolution vector field is interesting when one describes a thermodynamics problem geometrically (see \cite{simoes2020geometry}).

It is denoted as $\varepsilon_H$ and it satisfies
\begin{equation}\label{evo-def} 
\mathcal{L}_{\varepsilon_H}\eta=dH-\mathcal{R}(H)\eta,\qquad \eta(\varepsilon_H)=0.
\end{equation} 
In local coordinates it is given by
\begin{equation}\label{evo-dyn}
	\varepsilon_H=\frac{\partial H}{\partial p_i}\frac{\partial}{\partial q^i}  - \left (\frac{\partial H}{\partial q^i} + \frac{\partial H}{\partial z} p_i \right)
	\frac{\partial}{\partial p_i} + p_i\frac{\partial H}{\partial p_i} \frac{\partial}{\partial z},
\end{equation}
so that the integral curves satisfy the evolution equations
\begin{equation}\label{evo-eq}
	\dot{q}^i= \frac{\partial H}{\partial p_i}, \qquad \dot{p}_i = -\frac{\partial H}{\partial q^i}- 
	p_i\frac{\partial H}{\partial z}, \qquad \dot{z} = p_i\frac{\partial H}{\partial p_i}.
\end{equation}

The evolution and Hamiltonian vector field are related by
$\varepsilon_H = X^c_H + H \mathcal{R}.$ There are two possible approaches to review the dynamics explained in terms of evolutionary vector fields. 

Assume that $\gamma$ is a section of the canonical projection 
$\pi : T^*Q \times \mathbb R \longrightarrow Q \times \mathbb{R}$, say
$\gamma :  Q \times \mathbb{R} \longrightarrow T^*Q \times \mathbb R$. In local coordinates we have $(q^i, z) \mapsto \gamma(q^i) = (q^i, \gamma_j(q^i), z)$. Therefore, we can define the projected evolution vector field $\mathcal{E}_H^\gamma = T\pi \circ \mathcal{E}_H \circ \gamma$.

\begin{theorem}\label{hj_evolution_I}
 Assume that a section $\gamma$ of the projection $T^*Q \times \mathbb{R} \longrightarrow Q \times \mathbb{R}$
is such that $\gamma(Q\times \mathbb{R})$ is a coisotropic submanifold of 
$(T^{*}Q\times \mathbb{R}, \eta_Q)$, and $\gamma_z (Q)$ is a Legendrian submanifold of $(T^{*}Q\times \mathbb{R}, \eta_Q)$, for any $z \in \mathbb{R}$. 
Then, the vector fields $\mathcal{E}_H$ and $\mathcal{E}_H^{\gamma}$ are $\gamma$-related if and only if 
\begin{equation}\label{nuevo2}
d( H \circ \gamma) + \gamma_o \, \gamma^*(\Theta_Q) = 0 \; ,
\end{equation}
where $\gamma_o =  \frac{\partial H}{\partial z} + \frac{\partial H}{\partial p_i} \frac{\partial \gamma_i}{\partial z}$.
\end{theorem}

Equation (\ref{nuevo2}) is referred to as the {\it Hamilton-Jacobi equation for the evolution vector field}.
A section $\gamma$ fulfilling the assumptions of the theorem and the Hamilton-Jacobi equation will
be called a {\it solution} of the Hamilton-Jacobi problem for the evolution vector field of $H$.

\textbf{Approach II}
We will keep the notations of the previous subsection, although now $\gamma$ is a section
of the canonical projection 
$\pi : T^*Q \times \mathbb R \longrightarrow Q$, say
$\gamma :  Q \to T^*Q \times \mathbb R$. In local coordinates we have $(q^i) \mapsto \gamma(q^i) = (q^i, \gamma_j(q^i), \gamma_z(q^i))$. As in previous sections, we define the projected evolution vector field $\mathcal{E}_H^\gamma = T\pi \circ \mathcal{E}_H \circ \gamma $.

In this context, we have the following.

\begin{theorem}
 Assume that a section $\gamma$ of the projection $T^*Q \times \mathbb{R} \to Q$
is such that $\gamma(Q)$ is a Legendrian submanifold of $(T^{*}Q\times \mathbb{R}, \eta_Q)$. 
Then, the vector fields $\mathcal{E}_H$ and $\mathcal{E}_H^{\gamma}$ are $\gamma$-related if and only if
\begin{equation}\label{nuevo4}
d (H \circ \gamma) = 0,
\end{equation}
where $\gamma = j^1 f,$ for some $f:Q \to \mathbb{R}$.

\end{theorem}

Equation (\ref{nuevo4}) is referred to as the {\it Hamilton-Jacobi equation for the evolution vector field}.
A section $\gamma$ fullfilling the assumptions of the theorem and the Hamilton-Jacobi equation will
be called a {\it solution} of the Hamilton-Jacobi problem for the evolution vector field of $H$.

\textbf{An Example}

We again consider the Hamiltonian $H$ of the parachute equation, as in~\eqref{eq:parachute_hamiltoniam}:
\begin{equation}
    H(q,p,z) = \frac{1}{2 m}{(p+2\lambda z)}^2 + \frac{mg}{2 \lambda} (e^{2 \lambda q} -1),
\end{equation}
where $\lambda, g \in \mathbb{R}$ are constants. The extended phase space is $T^*Q \times \mathbb{R} \simeq \mathbb{R}^3$.

Its evolution vector field is given by
\begin{equation*}
        \varepsilon_H = \frac{p + 2 \gamma z}{m} \frac{\partial }{\partial q} 
    -\left(mg e^{2 \lambda q } + 2 \gamma \frac{p}{m}(p + 2 \lambda z)\right) \frac{\partial }{\partial p} +
    \left(  \frac{p}{m}(p + 2 \lambda z)  \right) \frac{\partial }{\partial z}.
\end{equation*}

\textbf{Evolution vector field. Approach I}

We choose a section $\gamma : Q \times \mathbb{R} \to T^* Q \times \mathbb{R}$ as in the approach I for the Hamiltonian vector field. Again, we assume that $\gamma(Q \times \mathbb{R})$ is a coisotropic submanifold and $\gamma_{z_0}(TQ)$ are Lagrangian submanifolds for all $z_0$. Using Theorem~\ref{hj_evolution_I},  $\mathcal{E}_H$ and $\mathcal{E}^\gamma_H$ are $\gamma$-related if and only if
\begin{equation*}
2\, \lambda^{3} {z}^{2} \frac{\partial f}{\partial {q}} \frac{\partial^{2} f}{\partial {q}\partial {z}} -  \, g \lambda m^{2} e^{2 \, \lambda {q}} - \, \lambda \frac{\partial f}{\partial {q}} \frac{\partial^{2}f}{\partial {q}^{2}} + \frac{F_1}{2} \frac{\partial f}{\partial {q}} \frac{\partial^{2} f}{\partial {q}\partial {z}}+2F_2z=0,
\end{equation*}
where
\begin{equation*}
F_1=g m^{2} e^{2 \, \lambda {q}} - g m^{2} + \lambda \left(\frac{\partial f}{\partial {q}}\right)^{2} - 2 \, \lambda \frac{\partial f}{\partial {q}},\quad
    F_2=\lambda^2\left(\frac{\partial^{2}f}{\partial {q}^{2}} - {\left(\frac{\partial f}{\partial {q}} - 1\right)} \frac{\partial f}{\partial {q}} \frac{\partial^{2} f}{\partial {q}\partial {z}}\right)
\end{equation*}
\textbf{Evolution vector field. Approach II}

Now we choose a section $\gamma: Q \to TQ \times \mathbb{R}$ and assume that is image is Legendrian. Equivalently $\gamma = j^1 f$.  $\mathcal{E}_H$ and $\mathcal{E}^\gamma_H$ are $\gamma$-related if and only if
\begin{equation}
    d(H \circ j^1 f) = 0.
\end{equation}

That is,
\begin{equation}
    \frac{1}{m}{ \left(\frac{\partial \gamma_z}{ \partial q}
     +2\lambda {\gamma_z} \right)
     \left(\frac{\partial^2 \gamma_z}{ \partial q^2}
     +2\lambda \frac{\partial \gamma_z}{ \partial q} \right)} + {mg} e^{2 \lambda q} = 0.
\end{equation}

\section{Conclusions}

This review paper presents an enhancement of some former results published by the authors. The theories have been categorized, merged and enhanced in such a way that our scattered results are gathered here as subcases of more general structures that classify them. So, in this review we have presented a new classification of HJ theories according to two main identities: the Leibniz and the Jacobi identity. This classification depends on the character of the geometric structure under study: either fulfills the Leibniz identity, or the Jacobi identity, or both.
All the examples included in this review paper are new examples, hence, we are providing new content alongside the review of our former content.
This review also includes a HJ theory for conformal Hamiltonian systems that had not been included in any of our former publications. This is why the introduction of a HJ for conformal Hamiltonian systems is a novelty in this review. The HJ for conformal Hamiltonian systems is introduced as a particular case of the HJ for systems with external forces.

Our starting point in Section \ref{Sec-Poisson} is the geometric HJ Theorem \ref{th1} recast in the symplectic framework and then we have presented some immediate generalizations of this formulation within the Poisson category. We have presented a geometric HJ Theory for conformal Hamiltonian dynamics in Subsection \ref{sec-HJforExt}. This novel HJ formulation is obtained as a particular case of the HJ theory including external forces. After cosymplectic Hamiltonian dynamics is exhibited, a HJ formulation for this geometric background is given in Subsection \ref{Sec-HJ-cos}, constituting the proper setting for time dependent Hamiltonian dynamics. The dynamical equations in this section are in the realm of Poisson manifolds. Table \ref{Table-1} briefly summarizes the ingredients of this section.

As mentioned previously, in order to determine a compendium of various HJ theories, we have relaxed the Leibniz and Jacobi identity fulfillment (one or the other) to arrive at Jacobi and almost Poisson geometries, respectively. We plot the following diagram to summarize the geometric frameworks discussed in this survey. 

\medskip

\begin{center}
\fbox{ 
 \includegraphics[width=0.6\textwidth]{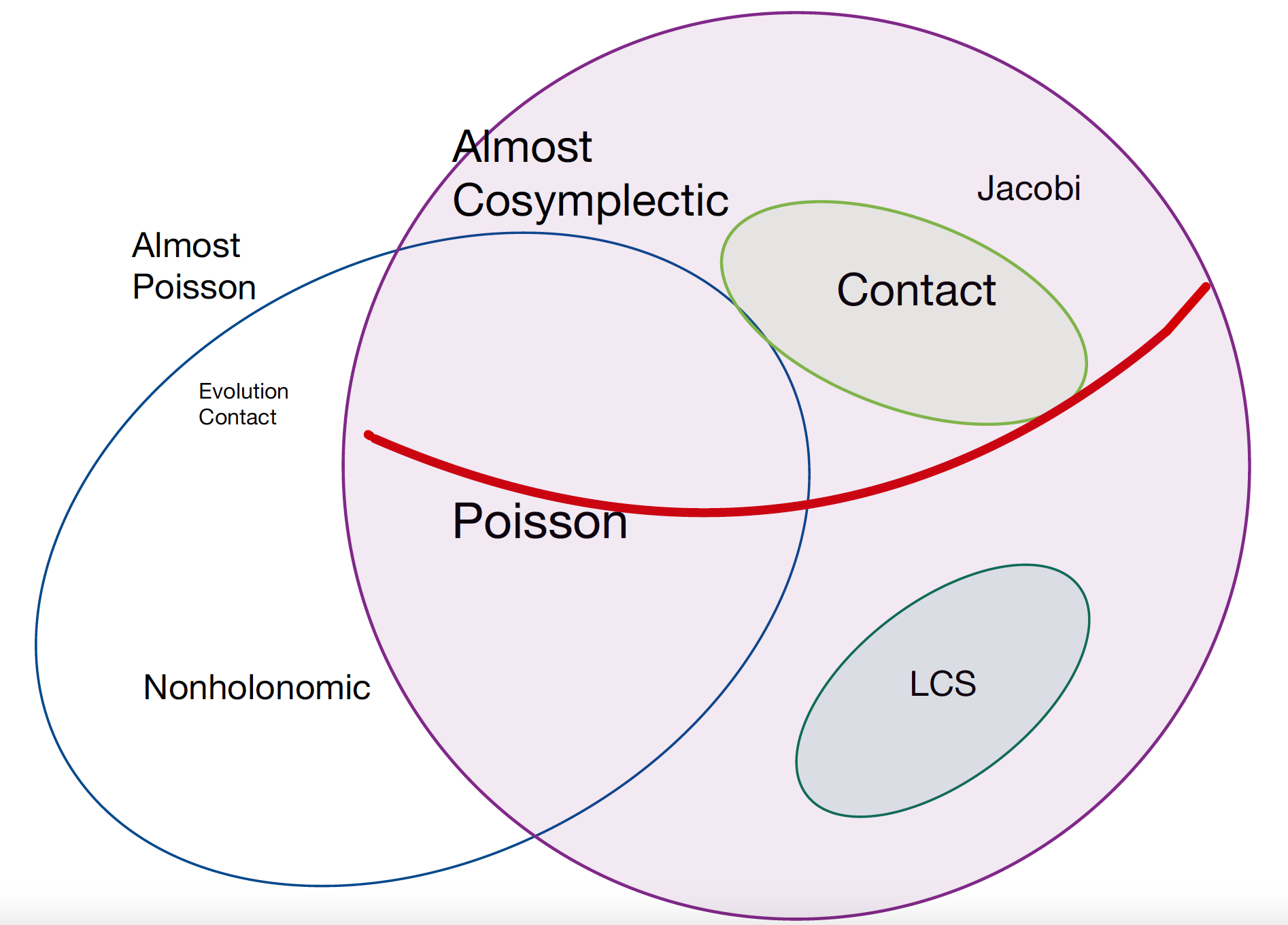}   
 }
\end{center}

\medskip

It is important to note that Jacobi manifolds can be interpreted as a foliation of LCS manifolds for even dimensions and contact manifolds for odd dimensions. According to the classification, in Section \ref{Sec-Jacobi}, we have 
examined the HJ theory for these two subcases of Jacobi manifolds, namely LCS and contact manifolds. In Subsection \ref{Sec-lcs-man} a geometric Hamilton-Jacobi theory is given for Hamiltonian dynamics on LCS manifolds, whereas in Subsection \ref{Sec-HJ-con-Ham} a geometric Hamilton-Jacobi theory is given for classical Hamiltonian dynamics on  contact manifolds. Such a direct classification of almost Poisson manifolds does not exist in the literature. Section \ref{Sec-Almost} is reserved for generalizations of the HJ theory for several particular instances of almost Poisson manifolds. In Subsection \ref{Sec-lin-almost}, we started the section by examining linear almost Poisson dynamics and the HJ theory for this geometry. Then, in Subsection \ref{Sec-HJ-nonholo}, we have focused on the HJ theory of nonholonomic Hamiltonian dynamics. HJ for evolution Hamiltonian dynamics on contact manifolds is given in Subsection \ref{Sec-Evo-Cont}. We present the following picture to exhibit several particular geometries in the realm of almost Poisson frameworks.  

\medskip

\begin{center}
\fbox{ 
 \includegraphics[width=0.6\textwidth]{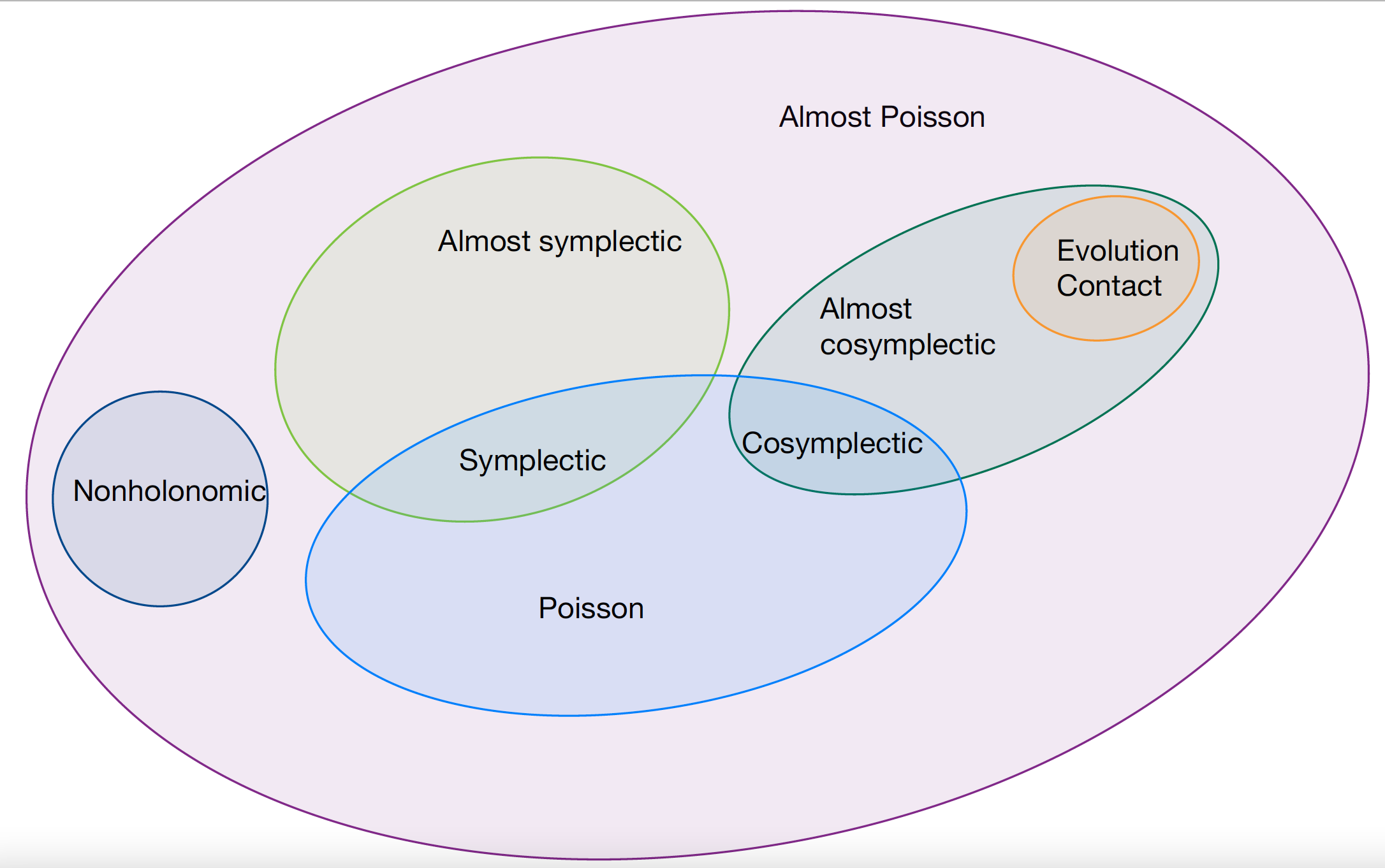}   
 }
\end{center}

\medskip

In the picture we also exhibit almost cosymplectic manifolds, for which contact manifolds and cosymplectic manifolds are  particular instances. Further, we comment that on a contact manifold two different  brackets can be determined. One is the classical contact bracket which is Jacobi, the other one is the evolution contact bracket, which is almost Poisson.

Two generalizations or/and complements of this survey can be given by discrete versions and field theoretical extensions of the classical formulations presented. As mentioned previously, there is extensive literature for these generalizations. We plan to pursue these directions in our future works in order to collect a  more complete picture of the geometric HJ theory. 
Let us finish this paper stating two open problems in the classical setting. 

\medskip
(1) We cite \cite{Contact-nonholonomic}, where nonholonomic constraints are imposed on both the Lagrangian and Hamiltonian picture on contact manifolds. We find interesting to pursue a geometric Hamilton-Jacobi theory for these constrained dynamics. 

\medskip
(2) As studied in Subsection \ref{Sec-lin-almost}, linear almost Poisson geometry can be recast on the dual of an almost Lie algebroid. We find interesting to discuss the HJ theory for the case of contact geometry in the realm of almost Lie algebroids, so we can arrive at a theory in the intersection area of almost Lie algebroids, contact manifolds and geometric Hamilton-Jacobi theory.

\bibliographystyle{abbrv}
\bibliography{references}
\end{document}